\documentclass[a4paper,fleqn]{cas-sc}

\usepackage[authoryear,longnamesfirst]{natbib}
\usepackage{amsmath,amssymb,mathtools,bm}
\usepackage{graphicx}
\usepackage{booktabs}
\usepackage{xcolor}
\usepackage{microtype}

\usepackage{tikz}
\usetikzlibrary{arrows.meta,calc,fit,positioning,shapes.geometric}

\newcommand{\R}{\mathbb R}
\newcommand{\ii}{\mathrm i}
\newcommand{\dd}{\mathrm d}
\newcommand{\Dt}{\Delta t}
\newcommand{\Dx}{\Delta x}
\newcommand{\Dy}{\Delta y}
\newcommand{\calR}{\mathcal R}

\newcommand{\calL}{\mathcal L}

\newcommand{\avg}[1]{\overline{#1}}

\newcommand{\ubar}{\bm u}
\newcommand{\W}{\bm W}
\newcommand{\U}{\bm U}
\newcommand{\F}{\bm F}
\newcommand{\G}{\bm G}
\newcommand{\x}{\bm x}

\newcommand{\bs}[1]{\boldsymbol{#1}}

\newcommand{\ww}{\bs w}
\newcommand{\WW}{\bs W}

\newcommand{\Ih}{\mathcal I_h}
\newcommand{\Rh}{\mathcal R_h}
\newcommand{\Ocal}{\mathcal O}
\newcommand{\OO}{\mathcal O}

\newcommand{\w}{\bm w}

\newcommand{\q}{\widehat{\bm q}}
\newcommand{\Gm}{\mathsf G}

\newcommand{\A}{\mathcal A}
\newcommand{\D}{\mathcal D}
\newcommand{\Pp}{\mathcal P}
\newcommand{\Ssel}{\mathcal S}
\newcommand{\E}{\mathcal E}

\newcommand{\e}{\mathrm e}

\begin{document}
\let\WriteBookmarks\relax
\def\floatpagepagefraction{1}
\def\textpagefraction{.001}
\shorttitle{Active Flux with Transported Acoustic Increments}
\shortauthors{K. Duraisamy}

\title [mode = title]{Fully Discrete Active Flux Method based on Transported Acoustic Increments for the Compressible Euler Equations}                      
\tnotemark[1]

\tnotetext[1]{This document is the result of a research
   project funded by Los Alamos National Laboratory.}

\author[1]{Karthik Duraisamy}[type=editor,
                        auid=000,bioid=1,
                        prefix=,
                        orcid=0000-0002-3519-8147]
\ead{kdur@umich.edu}
\ead[url]{caslab.engin.umich.edu}

\credit{Conceptualization, Methodology, Results, Writing}

\affiliation[1]{organization={Department of Aerospace Engineering \& Michigan Institute for Computational Discovery and Engineering  University of Michigan},
                city={Ann Arbor},
                state={MI},
                country={USA}}

\cortext[cor1]{Corresponding author}

\begin{abstract}
A fully discrete Active Flux method is proposed for the two-dimensional compressible Euler equations. The method builds on the evolution-operator formulation of multidimensional Active Flux proposed by Roe in which conservative cell averages are updated by unsplit flux quadrature while primitive point values are evolved by acoustic and advective subsolvers. This has recently been extended to the Euler equations by Barsukow, who developed an efficient approach by  adding acoustic increments at the Eulerian target point. It is shown that, in a frozen linearized setting, this misses the leading symmetric advection–acoustic cross interaction. The proposed method instead reconstructs the acoustic increment as a cellwise $Q_2$ field, using vertex, edge-midpoint, and cell-center acoustic increments, and evaluates this field at the convective foot of the target point. For constant frozen coefficients, the resulting point update reduces to the transported composition, eliminating the additive split defect and yielding the exact unsplit frozen evolution when the acoustic and advective generators commute. The resulting method preserves the exact locally linearized acoustic evolution operator of Barsukow, the compact stencil, and the conservative one-stage average update. Numerical experiments probe several facets of the numerical method. A mixed Fourier wave packet with {the linearized Euler equations} isolates the split error and shows third-order point accuracy for the transported update, compared with second-order behavior for the additive update. Isentropic vortex convection confirms reduced error constants, and an enlarged empirical CFL range. Nonlinear Gaussian acoustic pulse evolution demonstrates preservation of radial symmetry and near-third-order decay of the symmetry error. Low-Mach shear-layer tests show coherent vorticity evolution, {low} entropy dissipation, and absence of the coarse-grid secondary vortices seen in displayed DG/CG comparisons. Finally, a compressible under-resolved Kelvin–Helmholtz test demonstrates robust no-limiter evolution to late time with consistent entropy dissipation. { Fourier analysis of the fully coupled and point linearzed operators support the observed improvements in accuracy, acoustic phase and amplification behavior.}
\end{abstract}

%


\maketitle

\section{Introduction}
\label{sec:intro}

Despite many advances in high order methods~\citep{wang2013high,ferrer2023high,mani2023perspective} over the past two decades, effective computation of compressible flows  remains an outstanding challenge. This is particularly true when viewed via the lens of practical computations of complex  problems such as high Reynolds number aerodynamics, radiation hydrodynamics where one cannot afford  spatio-temporal resolutions corresponding to  asymptotically convergent regimes. Thus, there is a need for accurate, robust and efficient numerical methods that preserve important flow features at coarse resolutions, while providing the order of accuracy with refinement of the discretization.

\cite{roe2025musings} reminds us that the fundamental questions concerning a numerical method for hyperbolic conservation laws are easy to state but rarely respected:
{what information should be used to find the solution?}, and is
that information being propagated through the right geometric channels?
In one space dimension the question is benign: The exact domain of
dependence is a segment, and standard upwinding can be made to coincide
with it. {For the Euler equations in multiple dimensions, however,  the relevant signals propagate along a
range of directions,} and dimensionally split or one-dimensional
Riemann-solver-based discretizations inevitably approximate this cone
by anisotropic, axis-aligned proxies. Such schemes can be formally
high-order in the dispersive sense yet, in Roe's diagnosis, deficient
in \emph{bandwidth}~\citep{roe2021bandwidth}: they fail to transmit
short-wavelength multi-dimensional information cleanly, regardless of
formal order. A genuinely satisfactory numerical method should
therefore (i) carry a \emph{compact, symmetric} stencil whose discrete
domain of dependence approximates the actual Mach cone rather than its
axis-aligned envelope; (ii) propagate information along the true
multi-dimensional characteristics of the system rather than along the
coordinate axes~\citep{roe2017multidim}; (iii) be \emph{fully discrete}
and \emph{single stage}, so as not to import the ad-hoc damping and
effective stencil enlargement of multistage time integrators; and
(iv) deliver these properties uniformly across regimes, without
recourse to solution-dependent fixes such as low-Mach
preconditioning~\citep{guillard1999,dellacherie2010}, asymptotic-
preserving patches, or limiter machinery tuned per regime.
A scheme that propagates the genuinely multi-dimensional signals correctly, in
particular, the divergence-free, vorticity, and stationarity content
of the linear acoustic system, inherits the correct low-Mach
behaviour as a by-product, rather than as a fix.

The Active Flux method is one of the
few CFD frameworks designed from the outset around all four of the
principles above. It originates in the seminal `Scheme V' of
\cite{vanleer1977}, which proposed a method  for one-dimensional 
conservation laws, carrying alongside the cell average, an additional
pointwise degree of freedom and updates it by tracing the foot of the
characteristic. Around 15 years ago, \cite{eymann2011,eymann2013} revived this idea and extended the construction to
nonlinear systems and to genuinely two-dimensional linear problems. The defining feature of the framework is the simultaneous
storage of conservative cell averages {\em and}  { appropriately chosen pointwise 
degrees of freedom located on cell interfaces (vertices and edge
midpoints),} tied together by a globally continuous reconstruction. Two distinct mechanisms operate on
these unknowns: a conservative average update, which uses the point
values directly as quadrature data for the unsplit fluxes and
therefore requires no Riemann solver; and a non-conservative point
update, which is responsible for upwinding and is naturally formulated
as a (possibly approximate) evolution operator. When that evolution
operator is exact, the resulting method is one-stage, compact, stable
up to $\mathrm{CFL}=1$, and transports information genuinely
multi-dimensionally. The downstream consequences are remarkable: on
Cartesian grids the linearized Active Flux method is provably
stationarity- and vorticity-preserving and exhibits the correct
low-Mach limit \emph{without} any
preconditioning~\citep{barsukow2019cartesian,barsukow2021allspeed},
while the third-order theory (in 1D problems) developed by
~\cite{barsukow2021nonlinear}  closes the gap to high-order
accuracy on nonlinear systems without sacrificing compactness.

For systems of non-linear conservation laws in more than one space dimension,
however, the construction of a sufficiently accurate \emph{pointwise}
evolution operator remains the central technical difficulty, and the
literature pursues two distinct approaches. {The semi-discrete active flux approach, originally introduced by ~\cite{abgrall2023extensions} 
replaces the evolution operator with a non-standard
finite-difference reconstruction of the time derivative and integrates
the resulting ODE with a Runge-Kutta
method~\citep{abgrall2025semieuler}.} The
approach is general but reintroduces multistage diffusion and forces
small CFL numbers, sacrificing the very compactness and one-stage
purity that motivated Active Flux. The complementary line adheres to the original van~Leer-Roe ideals
and seeks an  evolution operator. The doctoral theses of
\cite{fan2017} and \cite{maeng2017} at the University of Michigan pursue this theme by
separately treating the acoustic and advective components of the
multi-dimensional Euler system; an extended Lax-Wendroff /
Cauchy-Kowalevskaya construction in the same vein is compared with
discontinuous Galerkin methods in~\citep{roe2018comparing}. The
acoustic part rests on a  result of
\cite{barsukow2022acoustic}, who derived the
exact evolution operator for the linearized acoustic equations on a
general class of (possibly non-smooth) initial data, exposing the
truly multi-dimensional structure of acoustic propagation and
providing a path for genuinely multi-dimensional Godunov and
Active Flux schemes.

{The starting point of the present work is the fully discrete formulation of Barsukow \cite{barsukow2025} in which (i) the point values
are evolved by an \emph{additive} operator split: at every node, the
acoustic subsystem is locally linearised about the frozen primitive
state and advanced by the exact evolution operator
of \cite{barsukow2022acoustic}; and (ii) the nonlinear advection
subsystem is advanced by a third-order foot-point iteration.} Cell
averages are then updated conservatively by space-time Simpson
quadrature of the unsplit Euler fluxes. The scheme is one-stage,
compact, and inherits the structure-preserving properties of its
acoustic ancestor, providing a { 
reference implementation} of evolution-operator-based multi-dimensional
Active Flux for the Euler system.  {We note that Helzel and collaborators ~\cite{chudzik2024bicharacteristics,chudzik2025fully} have recently pursued a non-split bicharacteristic
	construction,  presenting an alternative direction for the 2D Euler equations.}

The contribution of this paper is a fully discrete Active Flux scheme
that retains key elements of~\cite{barsukow2025}, but replaces
the additive composition by a \emph{transported acoustic-increment}
modification of the point-value update. The acoustic increment field
 is reconstructed across each cell as a $Q_2$
polynomial, using the increments naturally available at vertices and
edge midpoints together with a {cell-center} increment obtained by a
full-disk acoustic update; this $Q_2$ field is then evaluated at the
convective foot $P_f$, rather than added at the Eulerian node $P$ as
in~\cite{barsukow2025}. In a frozen linear  setting, it is  shown that
that this converts the additive split into a transported composition
that eliminates the leading $\mathcal{O}(\tau^2)$ piece of the
advection/acoustic cross-term defect. { In the non-linear case, this exact property cannot be proven, but is found to empirically yield higher accuracy due to an approximately multiplicative approximation.} Numerical experiments on the
isentropic vortex, planar acoustic propagation, and shear-layer
evolution show markedly improved resolution of vortex and shear
structures, with a larger stable CFL range. 
We provide Fourier analysis of the improved point schemes and also report on unsuccessful, more intrusive variants and identify the minimal-intervention modification that remains
compatible with its structure-preserving properties.

The remainder of the paper is organised as follows. Section~\ref{sec:infra} details the setup of the Active Flux infrastructure in two dimensions.
Section~\ref{sec:barsukow} reviews the fully discrete Barsukow scheme
and quantifies the cross-term defect of its additive split.
Section~\ref{sec:transported} introduces the transported-increment
method, the $Q_2$ reconstruction of the acoustic increment, and an
operator-level analysis of the modified split. {Section~\ref{sec:fourierfull} presents a Fourier analysis of the full operator , together with a Fourier analysis that
corroborates the observed behaviour at the level of the frozen acoustics.} Subsequent sections
report numerical experiments on the isentropic vortex, planar
acoustics, and shear flow.
{Additional details of the discretization and analysis is provided in the Appendix.}


\section{Active Flux Infrastructure in 2D}
\label{sec:infra}
We consider the two-dimensional compressible Euler equations
\begin{equation}
    \partial_t \U + \partial_x \F(\U)+\partial_y \G(\U)=0,
    \qquad
    \U=(\rho,\rho u,\rho v,E)^T,
\end{equation}
with ideal-gas closure
\begin{equation}
    E=\frac{p}{\gamma-1}+\frac12\rho(u^2+v^2).
\end{equation}
The corresponding primitive state is $
    \W=(\rho,u,v,p)^T,
    \quad
    c^2=\gamma p/\rho.$

The computational grid is Cartesian and periodic, with equal spacing $
    \Dx=\Dy=h.$
Cells are denoted $
    C_{ij}= [x_i,x_{i+1}]\times[y_j,y_{j+1}],$
with cell center $
    (x_{i+1/2},y_{j+1/2}).$
This work uses the following convention: conservative variables are used for cell averages, primitive variables are used for point values and reconstruction.  This is the choice made in Barsukow's multidimensional Euler Active Flux method, because it substantially simplifies the acoustic evolution operator while the average update remains conservative \citep{barsukow2025}.

\begin{figure}[t]
	\centering
	\begin{tikzpicture}[scale=1.2,>=Latex]
		\draw[step=2cm,gray!55,thin] (-2,-2) grid (4,4);
		\fill[blue!7] (0,0) rectangle (2,2);
		\draw[blue!60!black,very thick] (0,0) rectangle (2,2);
		\foreach \x in {-2,0,2,4}
		\foreach \y in {-2,0,2,4}
		\fill[black] (\x,\y) circle (2.3pt);
		\foreach \x in {-2,0,2,4}
		\foreach \y in {-1,1,3}
		\fill[teal!70!black] (\x,\y) circle (2.5pt);
		\foreach \x in {-1,1,3}
		\foreach \y in {-2,0,2,4}
		\fill[orange!85!black] (\x,\y) circle (2.5pt);
		\node[diamond,draw,fill=white,minimum size=7pt,inner sep=0pt] (cnode) at (1,1) {};
		\node[draw,fill=blue!65!black,minimum size=6pt,inner sep=0pt] (avg) at (1,1.55) {};
		\node[anchor=west] at (1.14,1.55) {$\overline{\U}_{ij}$};
		\node[anchor=west] at (1.14,1.02) {$\W^C_{ij}$};
		\node[anchor=east] (evlabel) at (-0.18,1.30) {$\W^{E_v}_{ij}$};
		\draw[-{Latex},thin] (evlabel.east) -- (0,1);
		\node[anchor=north] (ehlabel) at (1,-0.18) {$\W^{E_h}_{ij}$};
		\draw[-{Latex},thin] (ehlabel.north) -- (1,0);
		\node[anchor=north west] (vlabel) at (0.18,0.46) {$\W^V_{ij}$};
		\draw[-{Latex},thin] (vlabel.south west) -- (0,0);
		\begin{scope}[shift={(5.0,2.8)}]
			\fill[black] (0,0) circle (2.3pt); \node[right=4pt] at (0,0) {vertex family $V$};
			\fill[teal!70!black] (0,-0.48) circle (2.5pt); \node[right=4pt] at (0,-0.48) {vertical-edge family $E_v$};
			\fill[orange!85!black] (0,-0.96) circle (2.5pt); \node[right=4pt] at (0,-0.96) {horizontal-edge family $E_h$};
			\node[diamond,draw,fill=white,minimum size=7pt,inner sep=0pt] at (0,-1.44) {};
			\node[right=4pt] at (0,-1.44) {recovered center node};
			\node[draw,fill=blue!65!black,minimum size=6pt,inner sep=0pt] at (0,-1.92) {};
			\node[right=4pt] at (0,-1.92) {conservative cell average};
		\end{scope}
	\end{tikzpicture}
	\caption{Active Flux degrees of freedom and the $Q_2$ center node.  The cell
		average and three shared boundary families are independent.  The center can be
		recovered from the conservative Simpson constraint.}
	\label{fig:af-dofs}
\end{figure}

The degrees of freedom shown in Figure~\ref{fig:af-dofs} are:
\begin{itemize}
    \item conservative cell averages
$
        \avg{\U}_{ij}
        \approx
        \frac{1}{|C_{ij}|}\int_{C_{ij}}\U(x,y,t)\,\dd x\dd y;
$
    \item primitive point values at vertices,
$
        \W^V_{ij}\approx \W(x_i,y_j,t);
$
    \item primitive point values at vertical-edge midpoints,
$
        \W^{E_v}_{ij}\approx \W(x_i,y_{j+1/2},t);
$
    \item primitive point values at horizontal-edge midpoints,
$
        \W^{E_h}_{ij}\approx \W(x_{i+1/2},y_j,t).
$
\end{itemize}

Within each cell, the reconstruction is a tensor-product quadratic polynomial, i.e. a $Q_2$ or biparabolic polynomial.  Eight of its nine nodal values are already known from the point values on the boundary.  The ninth node is the cell-center value.  It is obtained by inverting the tensor-product Simpson rule in conservative variables.

Let $\U_c$ denote the unknown conservative cell-center node.  Let the four corner conservative values and four edge-midpoint conservative values be obtained by converting the primitive point values to conservative variables.  The tensor-product Simpson average is
\begin{equation}
    \avg{\U}_{ij}
    =
    \frac{1}{36}\sum_{\mathrm{corners}}\U
    +
    \frac{1}{9}\sum_{\mathrm{edge\ mids}}\U
    +
    \frac{4}{9}\U_c.
\end{equation}
Thus
\begin{equation}
    \U_c
    =
    \frac94\left[
    \avg{\U}_{ij}
    -
    \frac{1}{36}\sum_{\mathrm{corners}}\U
    -
    \frac{1}{9}\sum_{\mathrm{edge\ mids}}\U
    \right].
    \label{eq:center-node}
\end{equation}
Finally $\U_c$ is converted back to a primitive center node $\W_c$.  The primitive $Q_2$ reconstruction in cell $C_{ij}$ is then
\begin{equation}
    \W_{ij}^{\rm recon}(\xi,\eta)
    \triangleq
    \sum_{a=-1}^{1}\sum_{b=-1}^{1}
    \ell_a(\xi)\ell_b(\eta)\W_{ab}^{ij},
    \qquad
    -1\le \xi,\eta\le 1,
\end{equation}
where $\ell_{-1},\ell_0,\ell_1$ are the Lagrange polynomials at $-1,0,1$, and $\W_{ab}^{ij}$ are the nine primitive nodal values in the cell. Because each edge restriction is a one-dimensional quadratic interpolating the same three edge point values from either adjacent cell, the reconstruction is continuous across cell interfaces.  This is the characteristic Active Flux feature: the reconstruction is globally continuous, but the cell average remains a conservative degree of freedom~\citep{roe2017multidim}.

\section{Fully discrete Active Flux Scheme}
\label{sec:barsukow}

 In the fully discrete method of \cite{barsukow2025}, the point values are evolved by an approximate evolution operator built from an acoustic/advection decomposition of primitive Euler.  The average update is then performed conservatively using space-time Simpson quadrature of the unsplit Euler fluxes.

\subsection{Acoustic/advection split}

As proposed by Roe and his students~\citep{fan2017,maeng2017}, and further developed more recently~\citep{barsukow2022acoustic,barsukow2025}  the primitive Euler equations  are split into an acoustic part and an advection part. The acoustic part is
\begin{align}
    \rho_t+\rho\nabla\cdot \ubar&\triangleq 0,\label{eq:ac-rho}\\
    \ubar_t+\frac{1}{\rho}\nabla p&\triangleq 0,\label{eq:ac-u}\\
    p_t+\rho c^2\nabla\cdot\ubar&\triangleq 0.\label{eq:ac-p}
\end{align}
At each point value $P$, one freezes
$
    \rho_0\triangleq \rho^n(P),\ \
    c_0\triangleq c^n(P),
$
and introduces the acoustic scalar
\begin{equation}
    \pi\triangleq \frac{p}{\rho_0 c_0}.
\end{equation}
The locally linearized acoustic equations become
\begin{align}
    \ubar_t+c_0\nabla\pi&=0,\label{eq:lin-ac-u}\\
    \pi_t+c_0\nabla\cdot\ubar&=0.\label{eq:lin-ac-pi}
    \end{align}
This linear acoustic system is solved using the exact acoustic evolution operator.  In the Cartesian implementation, a vertex uses four quadrant wedges and an edge midpoint uses two half-plane wedges.  The exact acoustic evolution is applied to the local $Q_2$ reconstruction in each adjacent wedge, and the contributions are summed.  The paper's implementation pseudocode expresses this by precomputing polynomial-in-$r$ evolution coefficients for monomials and wedge geometries, then contracting those coefficients with the runtime polynomial coefficients \citep{barsukow2025}.

We find that it is possible to obtain explicit and exact expressions for the velocity and pressure updates in terms of the nodal coefficients, thus offering good insight for analysis. The expression for the pressure update at the vertical face node is  given in Appendix ~\ref{sec:acoustic_update}, and exploited later for Fourier analysis.

The density is recovered from the frozen acoustic invariant $
    \partial_t(\rho c_0^2-p)=0.$
If the acoustic solve returns $p_{\rm ac}$, then
\begin{equation}
    \rho_{\rm ac}=\rho_0+\frac{p_{\rm ac}-p_0}{c_0^2}.
\end{equation}

We denote the resulting acoustic point evolution by
$S_{\rm ac}(P;\tau)\W^n.$ The advection part is $
    \W_t+\ubar(\W)\cdot\nabla\W=0.$
    
Barsukow uses a third-order footpoint rule for this nonlinear advection subsystem.  For a target point $P$, the advection update is { defined as}
\begin{equation}
    S_{\rm adv}(P;\tau)\W^n
    \triangleq
    \W^n\left(P-\tau\ubar^n(P-\tau\ubar^n(P))\right).
    \label{eq:adv-footpoint}
\end{equation}
This is the multidimensional nonlinear-advection analogue of a {second order Picard  iteration.}  Barsukow proves the corresponding third-order statement for the nonlinear advection subsystem \citep{barsukow2025}.

Following this an additive point update~\citep{barsukow2025}  is performed
\begin{equation}
    \W^{n+\tau}_{\rm RB}(P)
    \triangleq
    S_{\rm ac}(P;\tau)\W^n
    +
    S_{\rm adv}(P;\tau)\W^n
    -
    \W^n(P).
    \label{eq:barsukow-additive}
\end{equation}

The update is deliberately simple.  Both subsolvers act on the same time-$n$ reconstruction.  The cell average is not updated by either split subproblem; instead, after point values have been evolved to $t^n+\Dt/2$ and $t^n+\Dt$, the conservative average update uses the unsplit Euler fluxes.

Barsukow's paper notes that the reconstruction and quadratures are third order, and the nonlinear advection subsolver is third order, but the locally linearized acoustics and additive splitting are formally second order in time \citep{barsukow2025}.

For one full time step $\Dt$, point values are computed at $t^n$, $t^{n+1/2}$, and $t^{n+1}$.  Let $
    \widehat{\F}_{i,j}^{\,m}$
denote the spatial Simpson edge flux at time level $m\in\{0,1/2,1\}$ on the vertical edge $x=x_i$.  Then the time-averaged edge flux is
\begin{equation}
    \widehat{\F}_{i,j}
    =
    \frac16\widehat{\F}_{i,j}^{\,0}
    +
    \frac46\widehat{\F}_{i,j}^{\,1/2}
    +
    \frac16\widehat{\F}_{i,j}^{\,1}.
\end{equation}
The cell average update is
\begin{equation}
    \avg{\U}^{n+1}_{ij}
    =
    \avg{\U}^{n}_{ij}
    -
    \frac{\Dt}{h}(\widehat{\F}_{i+1,j}-\widehat{\F}_{i,j})
    -
    \frac{\Dt}{h}(\widehat{\G}_{i,j+1}-\widehat{\G}_{i,j}).
    \label{eq:af-average-update}
\end{equation}
This update is conservative because it is a flux-difference update in conservative variables.

Given the history and recent developments, we refer to this scheme as the Discrete Roe-Barsukow method when comparing to other techniques. 
The semi-discrete active flux method of \cite{abgrall2025semieuler} is  presented in the Appendix and will be used for comparisons.

\subsection{Operator Error Analysis}

Linearize the primitive variables about a uniform state
\(
\WW_0\triangleq(\rho_0,u_0,v_0,p_0)^T,
\qquad
c_0^2=\gamma p_0/\rho_0.
\)
Write the perturbation as
\(
\ww=(\rho',u',v',p')^T.
\)
The frozen linearized primitive Euler equations are
\begin{align}
	\rho'_t+u_0\rho'_x+v_0\rho'_y+\rho_0(u'_x+v'_y)&=0,\label{eq:lin-rho}\\
	u'_t+u_0u'_x+v_0u'_y+\frac1{\rho_0}p'_x&=0,\label{eq:lin-u}\\
	v'_t+u_0v'_x+v_0v'_y+\frac1{\rho_0}p'_y&=0,\label{eq:lin-v}\\
	p'_t+u_0p'_x+v_0p'_y+\rho_0c_0^2(u'_x+v'_y)&=0.\label{eq:lin-p}
\end{align}
Equivalently,
\(
\partial_t\ww=(A+C)\ww,
\)
where
\(
A\triangleq -(u_0\partial_x+v_0\partial_y)I
\)
is the scalar advection generator and \(C\) is the frozen acoustic generator:
\[
C\ww\triangleq
-\begin{pmatrix}
	\rho_0(u'_x+v'_y)\\[1mm]
	p'_x/\rho_0\\[1mm]
	p'_y/\rho_0\\[1mm]
	\rho_0c_0^2(u'_x+v'_y)
\end{pmatrix}.
\]
{ For this frozen case, A and C are constant matrices.}

Let \(S_{\rm ac}(\cdot;\tau)\) and \(S_{\rm adv}(\cdot;\tau)\) denote the acoustic and
advective subsolvers.  Linearized about a uniform state, these are
\(
S_{\rm ac}(\tau)=e^{\tau C},
\qquad
S_{\rm adv}(\tau)=e^{\tau A}.
\)
Barsukow's additive point update has the linearized form
\(
E_B(\tau)=e^{\tau A}+e^{\tau C}-I.
\)
Expanding in \(\tau\),
\(
E_B(\tau)
=
I+\tau(A+C)+\frac{\tau^2}{2}(A^2+C^2)
+\frac{\tau^3}{6}(A^3+C^3)+\Ocal(\tau^4).
\)
The exact frozen unsplit evolution is
\(
E_{\rm ex}(\tau)=e^{\tau(A+C)},
\)
and its expansion is
\begin{align}
	E_{\rm ex}(\tau)
	={}&
	I+\tau(A+C)
	+\frac{\tau^2}{2}(A^2+AC+CA+C^2)\notag\\
	&+\frac{\tau^3}{6}
	(A^3+A^2C+ACA+CA^2+AC^2+CAC+C^2A+C^3)
	+\Ocal(\tau^4).
\end{align}
Hence
\(
E_B(\tau)-E_{\rm ex}(\tau)
=
-\frac{\tau^2}{2}(AC+CA)+\Ocal(\tau^3).
\)

When \(A\) and \(C\) commute, this is
\(
E_B(\tau)-E_{\rm ex}(\tau)
=
-\tau^2AC+\Ocal(\tau^3).
\)
This additive update therefore misses the symmetric
advection/acoustic cross interaction, which could have significant consequences as will be discussed.

\section{Transported acoustic-increment method}
\label{sec:transported}



 
\begin{figure}[t]
	\centering
	\begin{tikzpicture}[scale=3.5,>=stealth]
		
		\def\xmin{0}
		\def\xmax{3}
		\def\ymin{0}
		\def\ymax{3}
		
		\coordinate (P) at (2,2);
		\coordinate (P1) at (1.5,2);
		\coordinate (P2) at (1,2);
		
		\coordinate (Q) at (1.7,1.8);   
		\coordinate (D) at (2.36,2.18);   
		
		\def\R{0.45}
		
		\foreach \x in {0,1,2,3} {
			\draw[black, line width=0.35pt] (\x,0) -- (\x,3);
		}
		\foreach \y in {0,1,2,3} {
			\draw[black, line width=0.35pt] (0,\y) -- (3,\y);
		}
		
		\foreach \x in {0,1,2,3} {
			\foreach \y in {0,1,2,3} {
				\draw[black, fill=white, line width=0.4pt] (\x,\y) circle (0.035);
			}
		}
		
		\foreach \x in {0.5,1.5,2.5} {
			\foreach \y in {0,1,2,3} {
				\draw[black, fill=white, line width=0.4pt] (\x,\y) circle (0.028);
			}
		}
		
		\foreach \x in {0,1,2,3} {
			\foreach \y in {0.5,1.5,2.5} {
				\draw[black, fill=white, line width=0.4pt] (\x,\y) circle (0.028);
			}
		}
		
		\draw[black, fill=white, line width=1.0pt] (P) circle (0.05);
		\node[below right=1pt] at (P) {$P$};
		
		\fill[black] (Q) circle (0.03);
		\node[below left=1pt, align=right] at (Q) {\small $P_f$\\[3mm]\small \ \ \ \ convective foot};
		
		\fill[black] (D) circle (0.03);
		\node[above right=1pt, align=left] at (D) {\small $P+\ubar(P)\tau$\\[-1mm]\small material location at $t=\tau$};
		
		\draw[->, blue!70!black, line width=0.9pt] (Q) -- (P)
		node[midway, below=2pt, sloped] {\scriptsize $\ubar(P_f)\tau$};
		
		\draw[->, blue!70!black, line width=0.9pt] (P) -- (D)
		node[midway, above=2pt, sloped] {\scriptsize $\ubar(P)\tau$};

		\draw[red, dashed, line width=0.4pt] (P) circle (\R);
		\draw[blue, dashed, line width=0.1pt] (P1) circle (0.4);
		\draw[blue, dashed, line width=0.1pt] (1,2) circle (0.42);
		\draw[blue, dashed, line width=0.1pt] (1,1.5) circle (0.4);
		\draw[blue, dashed, line width=0.1pt] (1,1) circle (0.35);
		\draw[blue, dashed, line width=0.1pt] (1.5,1) circle (0.37);
		\draw[blue, dashed, line width=0.1pt] (2,1) circle (0.4);
		\draw[blue, dashed, line width=0.1pt] (2,1.5) circle (0.44);
		\draw[blue, dashed, line width=0.1pt] (1.5,1.5) circle (0.44);
		\draw[red,  line width=1.0pt] (Q) circle (\R);

	\end{tikzpicture}
	\caption{Schematic of the transported-increment. 
		The dashed red circle is the original acoustic disk centered at the grid node $P$. 
		That acoustic evolution is naturally associated with the material label that started at $P$, 
		which is located at $P+\ubar_0\tau$ after time $\tau$. To obtain the Eulerian value at the fixed node $P$ at time $\tau$, one instead needs acoustic information associated with the convective foot $P_f$, represented here by the red circle. This information is obtained by processing transported increments from the dashed circles.}
	\label{fig:transported_increment_schematic_2d_clean}
\end{figure}

The transported-increment modification is motivated by the insight from the above analysis. 
If the background velocity is frozen at $\ubar_0$, then the frozen linearized Euler equations 
involve the material derivative
\[
(\partial_t+\ubar_0\cdot\nabla)\WW + {C}\WW = 0,
\]
where ${C}$ denotes the frozen acoustic operator. Thus acoustics propagates relative 
to a frame moving with $\ubar_0$. Introduce the moving coordinate
\(
\bs{\xi}\triangleq \bs{x}-\ubar_0 t.
\)
Then
\[
(\partial_t+\ubar_0\cdot\nabla)\WW(\bs{x},t)=\partial_t \widetilde{\WW}(\bs{\xi},t),
\]
so in moving coordinates the frozen system reduces to a pure acoustic evolution. Thus, for
constant background conditions, the acoustic evolution computed from data centered at a point
$P$ is naturally associated with the material label that started at $P$, which is located at
$P+\ubar_0\tau$ at time $\tau$.

To recover the
Eulerian value at $P$, one instead needs information associated with the convective foot
\[
P_f \triangleq P-\ubar_0(P_f)\tau,
\qquad\text{so that}\qquad
\WW(P,\tau)=e^{\tau {C}}\WW(P_f)
\]
under the frozen constant-coefficient model. {In practice, to achieve third order accuracy, a two-stage Picard iteration is sufficient to approximate the location of the convective foot.}

\paragraph{An exact frozen acoustic update}

Note that placing the acoustic disks at the nodes, and limiting the acoustic CFL
\(c\Delta t/\Delta x \leq 0.5\), allows analytical expressions for the acoustic
update because edge-center nodes receive information from two semi-disks in the
two adjacent cells (this leads to clean expressions as in Appendix~\ref{sec:acoustic_update}), while corner nodes receive information symmetrically from
four quadrants.  This is the geometry illustrated in Figure~\ref{fig:transported_increment_schematic_2d_clean}.
It is exploited elegantly in the Appendices of ~\cite{barsukow2025}: since the disk cuts are always
centered on grid lines, the required geometric moments are only half-plane and
quadrant moments, and many of these quantities can be precomputed. {  Note that $c \Delta t/\Delta x \leq 0.5$ is not a fundamental restriction as, for instance, the limit can be extended to $c \Delta t/\Delta x \leq 1$ by considering cap moments (see below and Appendix~\ref{sec:footcentered}) for additional cells.}

When the acoustic disk is not placed at a node, as in the foot-centered variant
shown by the red circle in Figure~\ref{fig:transported_increment_schematic_2d_clean},
the situation changes substantially.  The acoustic center is now the convective
foot $P_f$  which is generally an arbitrary point inside a cell.
Even under the same acoustic CFL restriction, the disk may cut one vertical and one
horizontal grid line at off-center locations.  Thus the touched region consists
of the cell containing \(P_f\), possibly one neighboring cell in the \(x\)-direction,
possibly one neighboring cell in the \(y\)-direction, and possibly the diagonal
cell.

The author was painstakingly able to derive closed-form expressions for this off-center acoustic update by decomposing the disk integrals into three geometric building blocks: full-disk moments, single circular-cap (portion of the acoustic disk cut off by one grid line) moments, and double-cap (portion of the acoustic disk cut off by two grid lines) moments. {This type of moment calculation was introduced by ~\cite{chudzik2024bicharacteristics} in the context of a different evolution operator.}   The
usual Barsukow half-plane and quadrant moments are replaced by off-center moments
of the reconstructed polynomial data over these cap regions, together with the
corresponding radius-derivative and acoustic-history moments required by the
pressure and velocity formulas. In the limiting cases where the acoustic center
approaches a face midpoint or a corner node, these formulas reduce exactly to the
standard Barsukow semi-disk and quadrant updates. Analytical expressions  were verified to machine precision using quadrature. {Further details are provided in Appendix ~\ref{sec:footcentered}}.

Although this construction gives an exact frozen-acoustic point update centered
at the convective foot, it is not computationally attractive in its direct form.
The geometry is different for every point value and every time step, so the
moments cannot be tabulated once and reused as in the node-centered method.
Moreover, each update must dynamically identify the cut cells, translate the
local polynomial reconstruction to coordinates centered at \(P_f\), and evaluate
the single- and double-cap closed forms.  In our implementation, this direct
closed-form foot-centered update was approximately \(1000\times\) more expensive
than the original Barsukow point update, making it useful mainly as an analytical
reference rather than as a practical solver component.  Even in the linearized Euler regime, this method offered only mild improvements {at more than an order of magnitude higher} cost over a much more effective alternate which follows.

\paragraph{An efficient and accurate alternative}

We  propose a method  based on the key 
insight that, while the split scheme earlier evolves $e^{\tau {C}}\WW(P)$ is {\em not}
the correct Eulerian value at the fixed grid node $P$ at time $\tau$, it {\bf is the exact
frozen acoustic evolution associated with the material label that was at $P$ at time $0$,} and
is therefore naturally associated with the downstream point $P+\ubar_0\tau$.  

{ Define the  acoustic increment field:
\begin{equation}
	\Delta_{\rm ac}^n(\cdot;\tau)
	\triangleq S_{\rm ac}(\cdot;\tau)\W^n-\W^n.
\end{equation}

Instead of computing the acoustic increment directly at the convective  foot $P_f$ (see Figure~\ref{fig:transported_increment_schematic_2d_clean}), the present work exploits available information by a) computing acoustic increments at the nodal locations and cell centers (dashed circles in Figure~\ref{fig:transported_increment_schematic_2d_clean} ) and interpolating the acoustic increments at the convective foot $P_f$.
 
 The transported point update is
\begin{equation}
    \W^{n+\tau}_{\rm RB-TAI}(P)
    \triangleq
    S_{\rm adv}(P;\tau)\W^n
    +
    \calR_{\Delta}\left(P_f(P;\tau)\right),
    \label{eq:transported-update}
\end{equation}
where $\calR_{\Delta}$ is a $Q_2$ reconstruction of $\Delta_{ac}$.}

The acoustic increment is naturally available at vertices and edge midpoints after the acoustic update.  To reconstruct it as a $Q_2$ field, one also needs center increments.  Therefore the implementation computes an acoustic update of the cell-center primitive node $
    \W_{c,\rm ac}^{n+\tau} \triangleq S_{\rm ac}(\tau)\W_c^n,$
using the full-circle acoustic operator centered in the cell.  Then $
    \Delta_c \triangleq \W_{c,\rm ac}^{n+\tau}-\W_c^n.$
Together with the vertex and edge midpoint increments, these center increments define a cellwise $Q_2$ reconstruction $\calR_\Delta$.

Note that density is not an independent acoustic variable in the frozen Euler acoustic subsystem, and therefore its acoustic change must be tied to the pressure change through the invariant \(\rho-p/c^2\), otherwise the transported update introduces spurious thermodynamic error at the convective foot.  Hence
\begin{equation}
	\rho_{\rm ac}(P_f)
	=
	\rho^n(P_f)
	+
	\frac{p_{\rm ac}(P_f)-p^n(P_f)}{c_{P_f}^2},
\end{equation}
and therefore the transported acoustic increment used in the point update is
\begin{equation}
	R_\Delta(P_f(P;\tau))
	=
	\begin{pmatrix}
		\Delta p_/c_{P_f}^2\\
		\Delta u_{P_f}\\
		\Delta v_{P_f}\\
		\Delta p_{P_f}
	\end{pmatrix}.
\end{equation}


After applying the acoustic evolution operator for a time \(\tau\) at the vertices and edge midpoints, we obtain
\[
W^{V,{\rm ac}}_{i,j}(\tau),
\qquad
W^{E_v,{\rm ac}}_{i,j}(\tau),
\qquad
W^{E_h,{\rm ac}}_{i,j}(\tau).
\]
The corresponding acoustic increments are
\[
\Delta^V_{i,j}
=
W^{V,{\rm ac}}_{i,j}(\tau)-W^V_{i,j},
\  \ 
\Delta^{E_v}_{i,j}
=
W^{E_v,{\rm ac}}_{i,j}(\tau)-W^{E_v}_{i,j},
\ \
\Delta^{E_h}_{i,j}
=
W^{E_h,{\rm ac}}_{i,j}(\tau)-W^{E_h}_{i,j}.
\]
These eight boundary increments are not enough to define a \(Q_2\) polynomial in
a cell.  A ninth value, the cell-center increment, is also required.

Let \(C_{i,j}\) be the Cartesian cell
\(
C_{i,j}
=
[x_i,x_{i+1}]\times[y_j,y_{j+1}],
\)
with center
\(
(x_{i+1/2},y_{j+1/2}).
\)
The primitive center value \(W^C_{i,j}\) is not an independent degree of
freedom.  It is obtained from the conservative cell average by inverting the
two-dimensional Simpson rule.

\[
	U^C_{i,j}
	=
	\frac{9}{4}
	\left[
	\bar U_{i,j}
	-
	\frac{1}{36}
	\Big(
	U^V_{i,j}
	+U^V_{i+1,j}
	+U^V_{i,j+1}
	+U^V_{i+1,j+1}
	\Big)
	-
	\frac{1}{9}
	\Big(
	U^{E_h}_{i,j}
	+U^{E_h}_{i,j+1}
	+U^{E_v}_{i,j}
	+U^{E_v}_{i+1,j}
	\Big)
	\right].
\]
The primitive center value is then
\(
W^C_{i,j}=W(U^C_{i,j})=(\rho^C_{i,j},u^C_{i,j},v^C_{i,j},p^C_{i,j})^T.
\)
Together with the eight boundary point values, \(W^C_{i,j}\) defines the
cellwise primitive \(Q_2\) reconstruction \(R_W^n\).

This value is not obtained from a conservation constraint, because
\(\Delta_{\rm ac}\) is not a conservative cell-average variable.  Instead, it is
computed by applying the same locally linearized acoustic evolution operator at
the cell center. The acoustic solve at the center freezes
\(
\rho_0=\rho^C_{i,j},
\qquad
c_0=c^C_{i,j},
\)
and evolves the local primitive reconstruction \(R_W^n\) by the exact acoustic
operator centered at \((x_{i+1/2},y_{j+1/2})\):
\[
	\Delta^C_{i,j}
	=
	W^{C,{\rm ac}}_{i,j}(\tau)-W^C_{i,j}.
\]

\paragraph{Cellwise \(Q_2\) reconstruction of the acoustic increment.}

Introduce local reference coordinates on cell \(C_{i,j}\):
\[
\xi
\triangleq
\frac{2(x-x_{i+1/2})}{\Delta x},
\qquad
\eta
\triangleq
\frac{2(y-y_{j+1/2})}{\Delta y},
\qquad
(\xi,\eta)\in[-1,1]^2.
\]
Let the one-dimensional quadratic Lagrange basis at nodes
\(-1,0,+1\) be
\[
\ell_-(s)\triangleq\frac12s(s-1),
\qquad
\ell_0(s)\triangleq 1-s^2,
\qquad
\ell_+(s)\triangleq \frac12s(s+1).
\]
The cellwise reconstruction of the acoustic increment is the tensor-product
interpolant
\[
	\begin{aligned}
		R_\Delta|_{C_{i,j}}(x,y;\tau)
		&\triangleq
		\ell_-(\xi)\ell_-(\eta)\,\Delta^V_{i,j}
		+
		\ell_0(\xi)\ell_-(\eta)\,\Delta^{E_h}_{i,j}
		+
		\ell_+(\xi)\ell_-(\eta)\,\Delta^V_{i+1,j}
		\\
		&\quad
		+
		\ell_-(\xi)\ell_0(\eta)\,\Delta^{E_v}_{i,j}
		+
		\ell_0(\xi)\ell_0(\eta)\,\Delta^C_{i,j}
		+
		\ell_+(\xi)\ell_0(\eta)\,\Delta^{E_v}_{i+1,j}
		\\
		&\quad
		+
		\ell_-(\xi)\ell_+(\eta)\,\Delta^V_{i,j+1}
		+
		\ell_0(\xi)\ell_+(\eta)\,\Delta^{E_h}_{i,j+1}
		+
		\ell_+(\xi)\ell_+(\eta)\,\Delta^V_{i+1,j+1}.
	\end{aligned}
\]

The transported point update is then evaluated by locating the cell containing
the convective foot ${P_f}$,
computing its local coordinates \((\xi_{P_f},\eta_{P_f})\), and setting
\[
	W_{P_f}^{n+\tau}(P)
	=
	S_{\rm adv}(P;\tau)W^n
	+
	R_\Delta(P_f;\tau).
\]

\paragraph{Why the center increment is essential.}

The eight boundary increments determine the value of \(R_\Delta\) on the cell
boundary, but they do not determine the interior \(Q_2\) bubble
\(
\ell_0(\xi)\ell_0(\eta)
=
(1-\xi^2)(1-\eta^2).
\)
{The \(\Delta^C_{i,j}\) is computed using the same acoustic operator used at the boundary
for the following reasons

First, it makes \(R_\Delta\) a genuine \(Q_2\) interpolation of the acoustic
increment field, rather than a boundary-only extension.  For a smooth acoustic
increment field this gives the usual third-order interpolation accuracy,
\[
R_\Delta(x,y;\tau)
=
\Delta_{\rm ac}(x,y;\tau)
+
O(\Delta x^3+\Delta y^3).
\]
Under hyperbolic scaling \(\tau=O(\Delta x)=O(\Delta y)\), this interpolation
error is of the same order as the intended Active Flux spatial accuracy.

Second, it preserves constant states.  If \(W^n\) is uniform, then the exact
acoustic evolution returns the same value at vertices, edge midpoints, and cell
centers, so
\(
\Delta^V=\Delta^{E_v}=\Delta^{E_h}=\Delta^C=0.
\)
Therefore
\(
R_\Delta\equiv 0,
\)
and the transported update preserves uniform states.}

\subsection{Operator split analysis}

The transported acoustic-increment update can be written as
\begin{equation}
    W^{n+\tau}(P)
    =S_{\rm adv}(P;\tau)(W^n)
  {   +\Delta_{\rm ac}^n(P_f;\tau),}
    \label{eq:transported-physical-update}
\end{equation}
where
\begin{equation}
    \Delta_{\rm ac}^n(P_f;\tau)
    =S_{\rm ac}(P_f;\tau)W^n-W^n(P_f).
\end{equation}
After linearization about a constant state, the shift $P_f=P-\tau\ubar_0$ is exactly the frozen advection operator $e^{\tau A}$.  Thus the acoustic increment is transformed as
\begin{equation}
    \Delta_{\rm ac}(P-\tau\ubar_0;\tau)
    =e^{\tau A}(e^{\tau C}-I)w^n(P).
\end{equation}
The advective part contributes $e^{\tau A}w^n(P)$.  Hence
\begin{equation}
    E_{\rm tr}(\tau)
    =e^{\tau A}+e^{\tau A}(e^{\tau C}-I)
    =e^{\tau A}e^{\tau C}.
    \label{eq:transported-operator-justified}
\end{equation}
{This is the key intuition behind the present method: Accuracy follows from transporting the acoustic increment, not from re-solving the acoustic equation on a new state.}

{In a frozen constant-coefficient model,  $[A,C]=0$, and thus}
\begin{equation}
    E_{\rm tr}(\tau)=e^{\tau A}e^{\tau C}=e^{\tau(A+C)}.
\end{equation}
Therefore the transported update is not merely a smaller-error split in this model: it is exactly the frozen unsplit advection-acoustic point evolution.  In variable-coefficient or fully discrete settings this equality is broken by coefficient variation and interpolation error, but the leading additive cross-term defect is removed.

More generally, without assuming $[A,C]=0$,
\[
    E_{\rm tr}(\tau)-E_{\rm exact}(\tau)=\frac{\tau^2}{2}(AC-CA)+O(\tau^3)                   =\frac{\tau^2}{2}[A,C]+O(\tau^3).
\]
Thus, the transported Barsukow replaces the additive method's symmetric defect by a commutator defect.

{
\noindent {\em Remark:} It is important to recognize that the above analysis provides the underlying {\em intuition on why} the method works. When applied to non-linear problems,  the acoustic increment at each quadratic node is generated using a different local frozen state. The resulting physical increments are assembled into cellwise $Q_2$ interpolants and evaluated at a nonlinear, iteratively computed foot. Density is recovered separately from the interpolated pressure increment, and the point values are coupled to conservative cell averages through time-dependent nonlinear flux quadrature. There is therefore, in general, no single operator $C$ for which the implemented nonlinear map can be written as $e^{\tau A}e^{\tau C}$. Translation also need not commute with the cellwise reconstruction or with node-dependent freezing. The frozen identity consequently establishes neither the local truncation error nor the global convergence order of the implemented nonlinear scheme. The actual constant-state fully discrete map, including reconstruction and cell-average coupling, is analyzed through the amplification matrices in Section~\ref{sec:fourierfull}.}

\subsection{Mixed wave packet: Numerical results and analysis}

To elucidate the characteristics of the above method, we propose a numerical test problem { corresponding to the Linearized Euler equations} with a uniform subsonic background
\[
    \rho_0=1,\qquad p_0=1,\qquad
    u_0=M\cos\phi,\qquad v_0=M\sin\phi,
\]
with \(M=0.2\), \(\phi=25^\circ\), and \(\gamma=1.4\). Initialize a sum of
independent linearized primitive eigenmodes.  The entropy mode has
\(
    \rho'=\epsilon_e\cos(k_1\cdot x),\qquad
    p'=u'=v'=0.
\), and the shear mode has
\(
    u'=-\epsilon_s\widehat k_{2,y}\cos(k_2\cdot x),\qquad
    v'= \epsilon_s\widehat k_{2,x}\cos(k_2\cdot x),\qquad
    \rho'=p'=0.
\)

The acoustic perturbation is
\[
\begin{aligned}
	p' &= \epsilon \cos(k\cdot x),\ \
	\rho' &= \frac{\epsilon}{c_0^2}\cos(k\cdot x),\ \
	u' &= s\,\frac{\epsilon}{\rho_0c_0}\widehat{k}_x\cos(k\cdot x),\ \
	v' &= s\,\frac{\epsilon}{\rho_0c_0}\widehat{k}_y\cos(k\cdot x),
\end{aligned}
\]
where \(s=+1\) gives the acoustic-plus branch and \(s=-1\) gives the
acoustic-minus branch.

{
The specific packet in the experiment is $
    k_1=2\pi(1,2),\quad
    k_2=2\pi(2,1),\quad
    k_3=2\pi(3,2),\quad
    k_4=2\pi(2,3),$
with all amplitudes equal to \(10^{-6}\).  For a Fourier wave vector \(k=(k_x,k_y)\), define
\(
K \triangleq |k|,\qquad
\widehat{k} \triangleq \frac{k}{K},\qquad
\omega_0(k) \triangleq u_0k_x+v_0k_y .
\)

The mixed-packet exact linear solution is
\[
\begin{aligned}
	W_{\rm ex}(x,y,t)
	=
	W_0
	&+\epsilon_e r_e(k_1)
	\cos\!\left(k_1\cdot x-\omega_0(k_1)t\right) 
	+\epsilon_s r_s(k_2)
	\cos\!\left(k_2\cdot x-\omega_0(k_2)t\right) \\
	&+\epsilon_+ r_+(k_3)
	\cos\!\left(k_3\cdot x-\left[\omega_0(k_3)+c_0|k_3|\right]t\right)+\epsilon_- r_-(k_4)
	\cos\!\left(k_4\cdot x-\left[\omega_0(k_4)-c_0|k_4|\right]t\right).
\end{aligned}
\]

 The repeated-step error
uses
\begin{equation}
	E_{j,m}^{(n)}=
	\frac{\|\Gm_m^nq_j^{\rm ex}-g_j^nq_j^{\rm ex}\|_2}
	{\|q_j^{\rm ex}\|_2}
	\label{eq:repeated-error}
\end{equation}
which holds the physical wave number and final time \(T=1\) fixed as
\(h=\kappa\to0\).

\begin{figure}[t]
	\centering
	\includegraphics[width=0.7\linewidth]
	{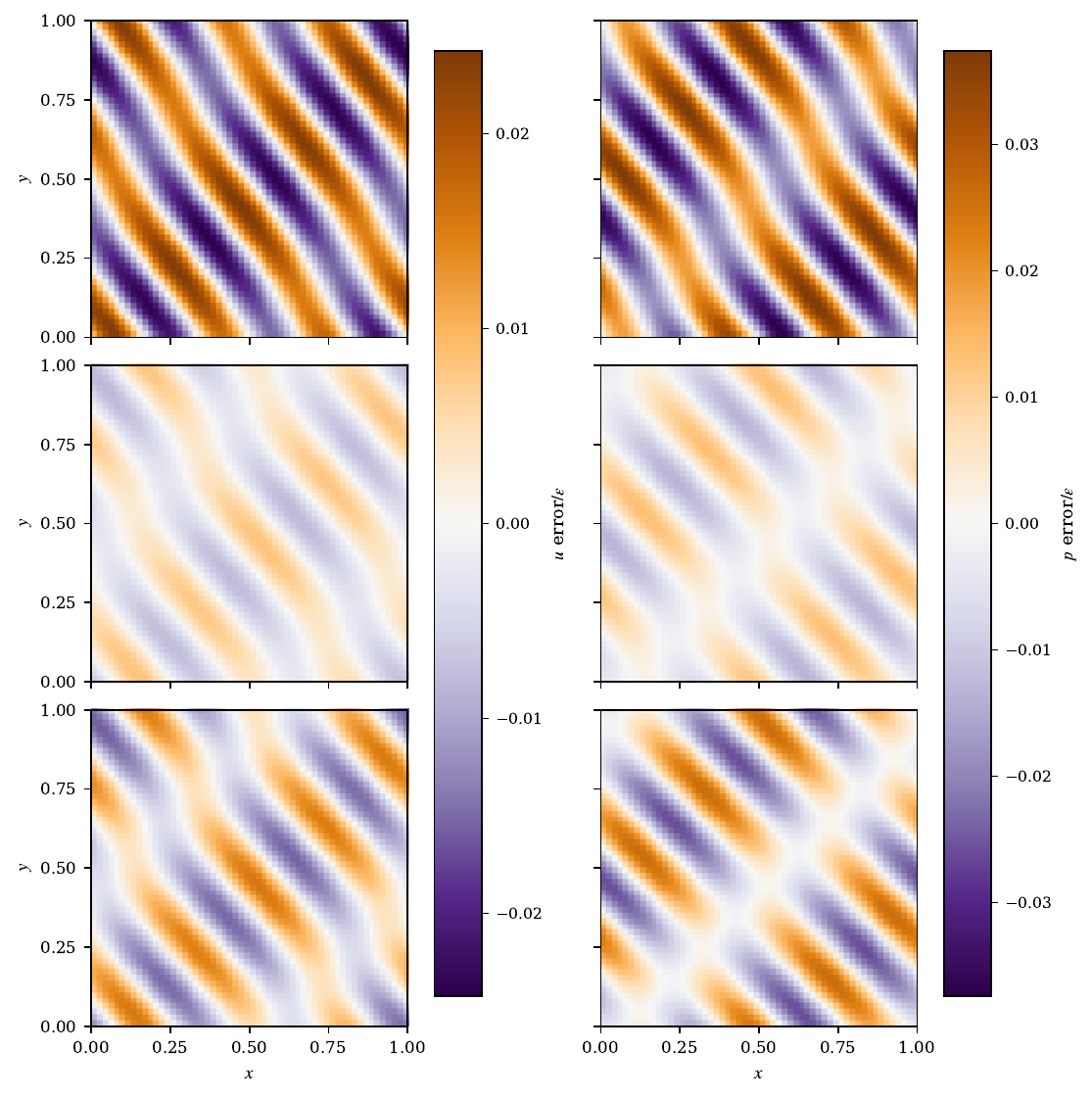}
	\caption{Mixed-packet method-minus-exact errors at \(M_0=0.2\), mean-flow
		angle \(25^\circ\), and \(\nu=0.25\), normalized by
		\(\epsilon=10^{-6}\), on a periodic \(56^2\) grid.  The comparison uses
		\(n=265\) equal steps, so \(t_n \approx 1\).  Rows from top to
		bottom are discrete RB, discrete RB-TAI, and semidiscrete AF with SSPRK3;
		columns from left to right are the \(u\)-velocity and pressure errors.  A
		single symmetric color scale is shared by all three schemes within each
		column.}
	\label{fig:m02-mixed-packet}
\end{figure}

Figure~\ref{fig:m02-mixed-packet}
compares the cell-average primitive errors after repeated updates.

\paragraph{Analysis} {To gain insight into the results,} consider a single acoustic branch Fourier mode
\(
\ww^n(P)=\epsilon r_s(k)e^{\ii k\cdot P},\quad s=\pm1.
\)
Let
\(
\alpha=\omega_0(k)\Delta t,\quad
\beta=s c_0K\Delta t .
\)
The exact frozen evolution is
\(
\ww_{\rm ex}(P,\Delta t)
=
\epsilon e^{-\ii(\alpha+\beta)}
r_s(k)e^{\ii k\cdot P}.
\)
The RB scheme yields
\[
\ww_B(P,\Delta t)
=
\epsilon
\left(
e^{-\ii\alpha}+e^{-\ii\beta}-1
\right)
r_s(k)e^{\ii k\cdot P}.
\]
\[
\noindent \textrm{Thus }	\ww_B-\ww_{\rm ex}=
	\epsilon\,\omega_0(k)\,s c_0K\,\Delta t^2\,
	r_s e^{\ii k\cdot P}
	+\OO(\Delta t^3).
\]

For any acoustic branch with \(\omega_0(k)\ne0\), \(K\ne0\), and nonzero amplitude,
this is a genuine \(\OO(\Delta t^2)\) one-step defect.  Under the hyperbolic scaling
\(
\Delta t=\nu h,
\)
the defect becomes $
\ww_B-\ww_{\rm ex}=\OO(h^2).$

We now consider the RB-TAI method. Let
\(
A=-(u_0\partial_x+v_0\partial_y)I
\)
be the frozen advection generator and \(C\) the frozen acoustic generator.  For
constant coefficients, \(A\) and \(C\) commute.  The transported increment update
with perfect evaluation of the increment is
\[
\begin{aligned}
	\ww_{\rm tr}(P)
	&=
	e^{\Delta t A}\ww^n(P)
	+
	e^{\Delta t A}\left(e^{\Delta t C}-I\right)\ww^n(P)
	=
	e^{\Delta t A}e^{\Delta t C}\ww^n(P).
\end{aligned}
\]
}
Since \([A,C]=0\),
\(
e^{\Delta t A}e^{\Delta t C}
=
e^{\Delta t(A+C)}.
\)
Thus the transported update has no \(\OO(\Delta t^2)\) defect in the frozen
constant-coefficient model.

For the acoustic Fourier branch, the same statement is just the phase identity
\(
e^{-\ii\alpha}
\left[
1+e^{-\ii\beta}-1
\right]
=
e^{-\ii(\alpha+\beta)}.
\)

Now include the transported $Q_2$ interpolation used in the point-update diagnostic.  Let \(\Rh\) be the cellwise $Q_2$
interpolant.  The transported diagnostic computes
\[
\ww_{h,{\rm f}}(P)
=
\Rh\ww^n(P_f)+\Rh\Delta_{\rm ac}(P_f),
\]
where
\(
\Delta_{\rm ac}=(e^{\Delta t C}-I)\ww^n.
\)
The exact transported value is
\(
\ww_{\rm ex}(P)=\ww^n(P_f)+\Delta_{\rm ac}(P_f).
\)
Hence
\(
\ww_{h,{\rm tr}}-\ww_{\rm ex}
=
(\Rh-I)\ww^n(P_f)+(\Rh-I)\Delta_{\rm ac}(P_f).
\)

For smooth fixed-wavenumber data, $Q_2$ interpolation gives
\(
(\Rh-I)\ww^n=\OO(h^3).
\)
Moreover,
\(
\Delta_{\rm ac}
=
(e^{\Delta t C}-I)\ww^n
=
\OO(\Delta t)
\)
for fixed \(k\).  Thus
\(
(\Rh-I)\Delta_{\rm ac}
=
\OO(\Delta t\,h^3)
=
\OO(h^4)
\)
under \(\Delta t=\OO(h)\).  Therefore
$
\ww_{h,{\rm tr}}-\ww_{\rm ex}=\OO(h^3).$
The leading visible error is the  interpolation of the transported state, not the
additive advection--acoustics split defect.

{

\section{Fourier analysis of the full operator}
\label{sec:fourierfull}
The Fourier setup and the exact acoustic map \(\A_K^\sigma\) are given in
Appendix~\ref{sec:acoustic-symbol}.  Set \(\lambda \triangleq \Delta t/h\) and
\(\nu \triangle c_0\Delta t/h\), and let \(X_K\) denote the target location.
At stage fraction \(\sigma\in\{0,1/2,1\}\), define the nondimensional
background displacement $
d_x^\sigma \triangleq \sigma\lambda u_0,\quad d_y^\sigma \triangleq \sigma\lambda v_0.$
The nonlinear two-Picard foot rule reduces exactly, after linearization about a
constant background, to reconstructed translation.  Hence
\begin{equation}
	\D_K^\sigma=\R(X_{K,x}-d_x^\sigma,X_{K,y}-d_y^\sigma)\otimes I_4.
	\label{eq:advection-map}
\end{equation}

\paragraph{Original additive RB}

Let \(\Ssel_K=e_K\otimes I_4\) select the old point value.  The original
additive RB point symbol is
\begin{equation}
	\Pp_K^\sigma=\A_K^\sigma+\D_K^\sigma-\Ssel_K.
	\label{eq:rb-point-map}
\end{equation}

\paragraph{Transported acoustic increments}

RB-TAI forms the complete nodal acoustic-increment field
\begin{equation}
	\Delta_K^\sigma=\A_K^\sigma-\Ssel_K,
	\qquad K\in\{V,E_v,E_h,C\},
	\qquad \Ssel_C=\mathcal C(z_x,z_y)\otimes I_4.
	\label{eq:tai-increments}
\end{equation}
If \(\Ih\Delta^\sigma\) denotes its cellwise-$Q_2$ interpolant (specified in
the appendix), the point map is simply
\begin{equation}
	\Pp_{K,\mathrm{TAI}}^\sigma
	=\D_K^\sigma+(\Ih\Delta^\sigma)
	(X_K-d^\sigma).
	\label{eq:tai-point-map}
\end{equation}

\paragraph{Semidiscrete Active Flux}

For the Cartesian Jacobian-split point evolution, let
\(\mathsf S_{\rm AF}(\xi,\eta)=h\widehat{\mathcal L}_{\rm AF}\) be the scaled
spatial symbol given
explicitly in Appendix~\ref{sec:sdaf-symbol}.  The method-of-lines system and
its SSPRK3 amplification matrix are
\begin{equation}
	\frac{\dd\q}{\dd t}=\frac1h\mathsf S_{\rm AF}\q,
	\qquad
	\Gm_{\rm sd}=R_3\!\left(\frac{\nu}{c_0}\mathsf S_{\rm AF}\right),
	\qquad R_3(z)=1+z+\frac{z^2}{2}+\frac{z^3}{6}.
	\label{eq:sdaf-amplification}
\end{equation}
Here \(R_3\) is evaluated as a matrix polynomial.

\subsection{The fully discrete RB amplification matrices}

For \(m\in\{\mathrm{RB},\mathrm{TAI}\}\), denote the corresponding point
map by \(\Pp_{K,m}^\sigma\).  The conservative average uses
the old, half-time, and full-time point values in space--time Simpson
quadrature.  At stage \(\sigma\), define the vertical- and horizontal-edge
state rows
\begin{align}
	E_{x,m}^\sigma=\frac{(1+z_y)\Pp_{V,m}^\sigma+4\Pp_{E_v,m}^\sigma}{6},\ \
	E_{y,m}^\sigma=\frac{(1+z_x)\Pp_{V,m}^\sigma+4\Pp_{E_h,m}^\sigma}{6},
\end{align}
and their time-Simpson averages
\begin{equation}
	\E_{x,m}=\frac{E_{x,m}^0+4E_{x,m}^{1/2}+E_{x,m}^1}{6},\qquad
	\E_{y,m}=\frac{E_{y,m}^0+4E_{y,m}^{1/2}+E_{y,m}^1}{6}.
	\label{eq:space-time-edge}
\end{equation}
Let \(B=[I_4,0,0,0]\) select the old cell average.  Linearizing the
conservative flux update and transforming its output back to primitive
perturbations gives
\begin{equation}
	G_{{\rm av},m}=B-\lambda(z_x-1)A_x\E_{x,m}
	-\lambda(z_y-1)A_y\E_{y,m}.
	\label{eq:average-row}
\end{equation}
The desired full Fourier amplification matrix is therefore
\begin{equation}
	\Gm_m(\xi,\eta;\nu,\W_0)=
	\begin{bmatrix}
		G_{{\rm av},m}\\[1mm]
		\Pp_{V,m}^1\\ \Pp_{E_v,m}^1\\ \Pp_{E_h,m}^1
	\end{bmatrix}\in\mathbb C^{16\times16}.
	\label{eq:full-symbol}
\end{equation}

An estimate of the form
$
\left|d_{h,\Delta t}(t)\right|=O(h^r),
$
is a one-step consistency statement. It does not by itself imply an $O(h^r)$ error after a fixed physical time, because the number of steps is then $O(h^{-1})$. In particular, an $O(h^3)$ one-step defect can, under ordinary accumulation, produce an $O(h^2)$ fixed-final-time error.

Temporal order refers to convergence as $\Delta t\to0$ after the spatial error has been held fixed or otherwise accounted for. Since the present method is a fully discrete construction rather than a method-of-lines discretization, temporal order must be assessed separately and is not inferred from a refinement in which $\Delta t/h$ is fixed.

The $O(h^3)$ estimate obtained for the transported point update concerns one application of the frozen constant-coefficient diagnostic. It shows that the leading $O(\Delta t^2)$ additive advection-acoustic defect is absent from that diagnostic, but it does not by itself establish third-order convergence after a fixed physical time.

The repeated-step calculation addresses a different question. It holds the physical wave number $K$ and final time $T$ fixed, sets
$
\kappa=Kh\longrightarrow0,
$
and applies the complete amplification matrix
$
N(h)=\frac{T}{\Delta t(h)}=O(h^{-1})
$
times. 

Figure ~\ref{fig:m02-accuracy} shows that only the acoustic RB branches lose an order: their one-step
and repeated-step slopes approach two, while every RB-TAI and semidiscrete AF
curve approaches three.  The matching one-step and repeated slopes show that
this distinction is not an artifact of the selected final time.

The right column of Figure~\ref{fig:m02-accuracy} is  a global test of the frozen linearized fully discrete map. The state vector contains the cell average and the vertex, vertical-edge-midpoint, and horizontal-edge-midpoint values of all four primitive variables. This calculation is not a measurement of temporal order at fixed $h$, and does not represent  a consistency analysis of the nonlinear variable-coefficient implementation.

\begin{figure}[t]
	\centering
	\includegraphics[width=0.96\linewidth]
	{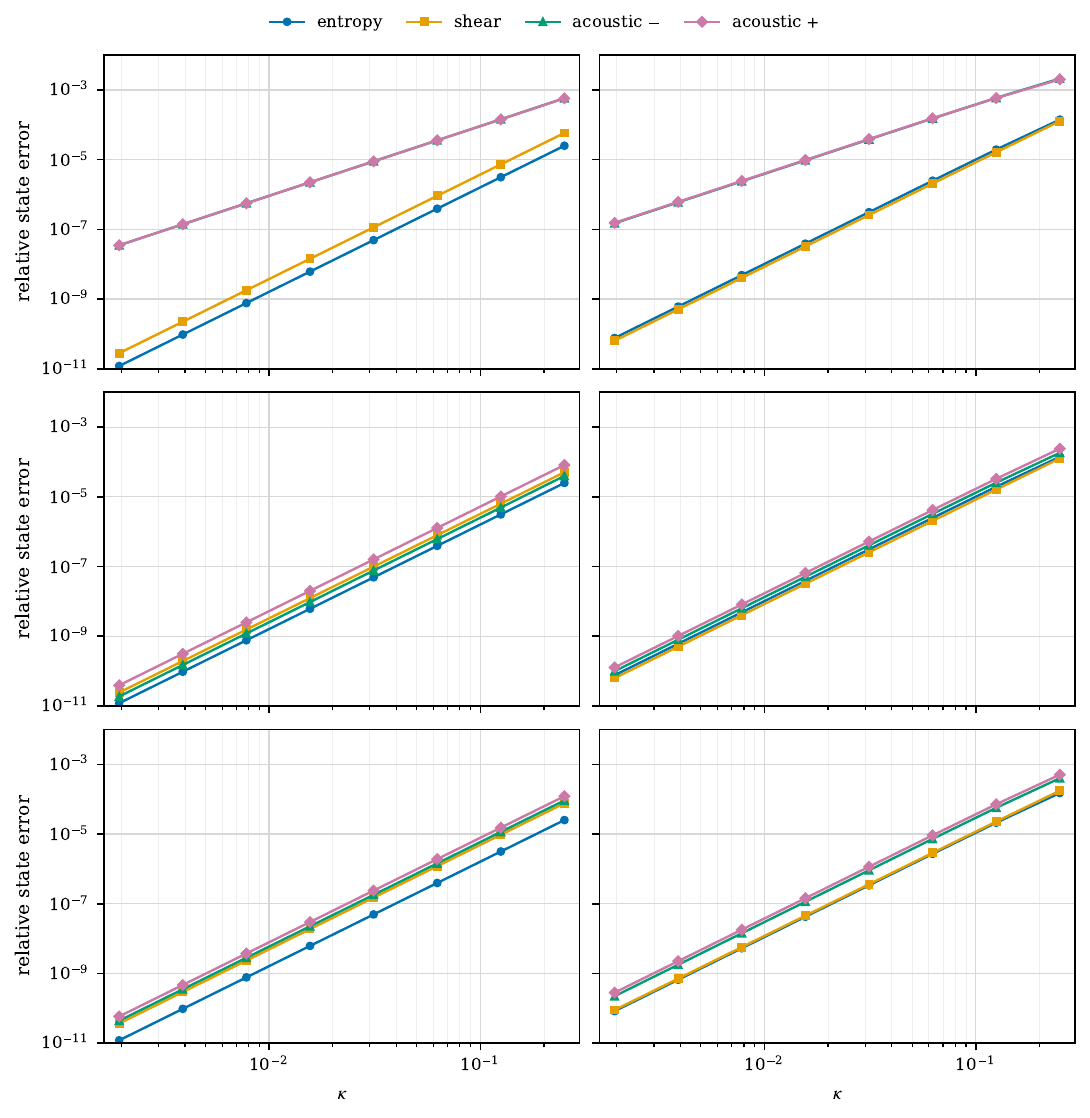}
	\caption{Relative errors for the entropy, shear, left-running acoustic, and
		right-running acoustic waves at \(\nu=0.25\).  Rows from top to bottom
		are discrete RB, discrete RB-TAI, and semidiscrete AF with SSPRK3; columns
		from left to right are the exact one-step residual 
		and the repeated-step error.  }
	\label{fig:m02-accuracy}
\end{figure}

For directional dispersion and dissipation, let \(\widetilde g_{j,m}\) be the
numerical eigenvalue matched to \(q_j^{\rm ex}\).  The four numerical eigenpairs
are assigned one-to-one by maximum normalized right-eigenvector overlap; this
resolves the exactly degenerate entropy and shear gains without duplicating a
numerical branch.  Define
\begin{equation}
	\epsilon^A_{j,m}=\bigl||\widetilde g_{j,m}|-1\bigr|,
	\qquad
	\epsilon^\phi_{j,m}=\left|\arg\!\left(
	\frac{\widetilde g_{j,m}}{g_j}\right)\right|.
	\label{eq:modal-amplitude-phase-errors}
\end{equation}

Figure~\ref{fig:m02-polar-errors} shows the absolute amplitude (left column) and phase errors (right column) at $\nu=0.25$.  The polar angle is the wave-vector direction, again showing the effectiveness of the Discrete RB-TAI method over the Discrete RB and Semi-Discrete Active Flux methods.

\begin{figure}[t]
	\centering
	\includegraphics[width=0.8\linewidth]
	{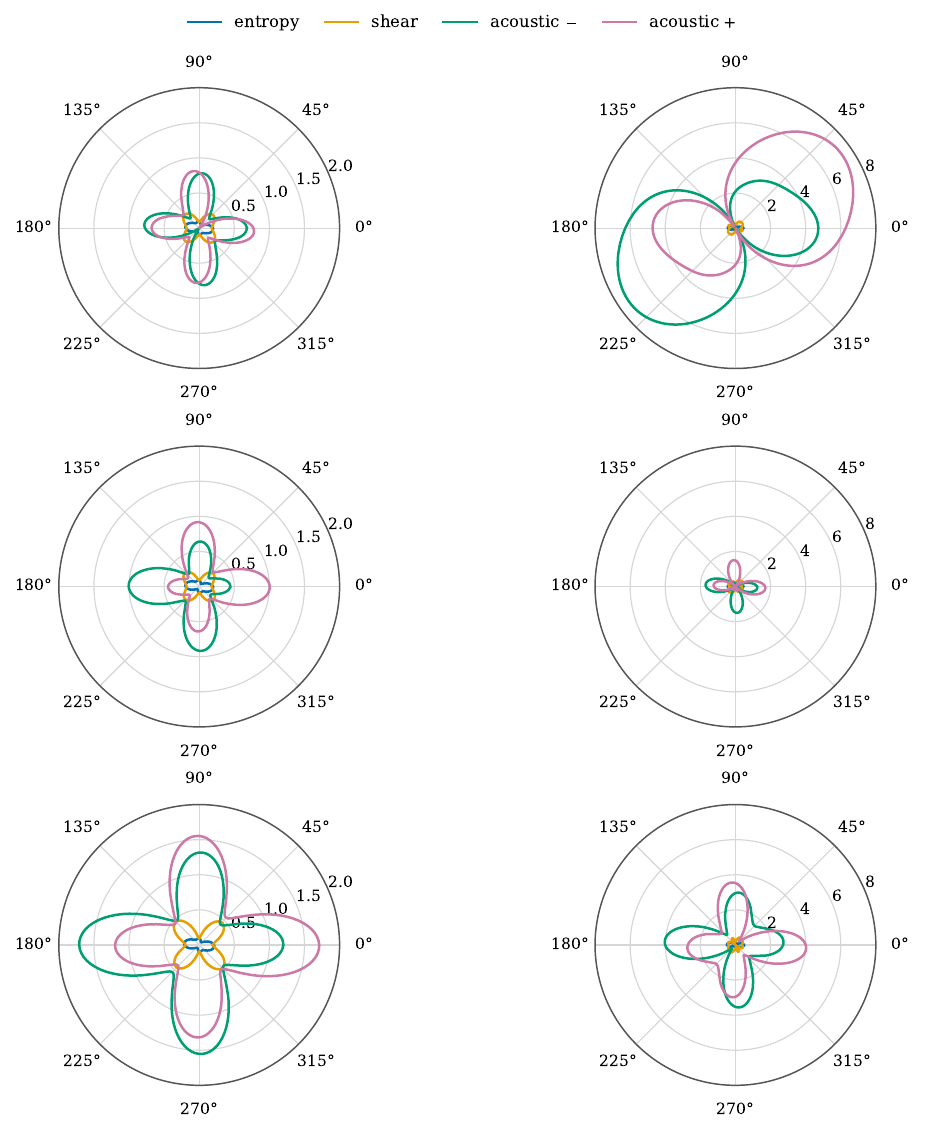}
	\caption{Directional modal errors at \(\nu=0.25\), and \(\kappa=\pi/4\), corresponding to eight
		cells per wavelength.   Rows
		from top to bottom are discrete RB, discrete RB-TAI, and semidiscrete AF with
		SSPRK3; columns from left to right are the absolute amplitude and phase errors
		in \eqref{eq:modal-amplitude-phase-errors}.  The four curves are the entropy,
		shear, left-running acoustic, and right-running acoustic waves.  Radial units
		are \(10^{-3}\) in the left column and \(10^{-4}\) radians in the right
		column, with a common radial range within each column.}
	\label{fig:m02-polar-errors}
\end{figure}

\subsection{Stability}

Maximizing over all physical and parasitic
branches gives \(\nu_{\rm crit}=0.195185\) for RB (a long-wave acoustic mode
at \(44.427^\circ\)), \(0.500000\) for RB-TAI (the wedge-geometry ceiling,
with no earlier spectral crossing), and \(0.276525\) for semidiscrete AF with
SSPRK3. 

At \(\nu=0.25\), RB is therefore globally unstable, whereas RB-TAI and
semidiscrete AF are spectrally stable.  Sampled powers through 1000 steps give
maximum norms 3.87 and 3.85 for the latter two methods, respectively, with no
continuing growth.  Thus the first row of Fig.~\ref{fig:m02-accuracy} is only a
directional accuracy diagnostic at a globally inadmissible CFL.
It is also noted that the above linear stability analysis is only valid for this particular
base flow configuration

\subsection{One-step point acoustic error}
For a more resolved view of the one-step acoustic error, let
\(e_{E_v,p}^T\) select pressure at the vertical-edge midpoint and define the
wave vector by
\((\xi,\eta)=\kappa(\cos\vartheta,\sin\vartheta)\). 
The exact frozen-acoustic pressure gain for the same branch is $
g_{p,{\rm ex}}^{s}(\xi,\eta,\nu)=\exp(-i\,s\,\nu\kappa).$

Set
\begin{equation}
	g_{p,m}^+(\kappa,\vartheta;\nu)=
	\frac{e_{E_v,p}^T\Gm_mq_+^{\rm ex}}
	{e_{E_v,p}^Tq_+^{\rm ex}}.
	\label{eq:pressure-response}
\end{equation}

Figure~\ref{fig:m02-pressure-gain-loci} shows that RB-TAI is closest to the exact pressure response.  The large transverse RB loop grows with both \(\kappa\) and
$\nu$, while semidiscrete AF follows RB-TAI at the smaller $\nu$ values but
separates more clearly at \(\nu=0.25\).   At \(\kappa=1\), the directional error peaks mildly near
\(\nu=0.3\) and then contracts toward \(\nu=0.5\), so proximity to the
wedge-geometry ceiling improves this one-step acoustic response.

From Figure~\ref{fig:m02-tai-pressure-gain-high-cfl} we see that at
\(\kappa=1\), the maximum directional values of
\(|g_{p,{\rm TAI}}^+-g_{p,{\rm ex}}^+|\) are respectively
\(8.48\times10^{-3}\), \(8.75\times10^{-3}\), \(7.02\times10^{-3}\),
and \(5.15\times10^{-3}\); thus the $\nu$ dependence is mildly nonmonotone
and accuracy improves as \(\nu\) approaches the exact-wedge geometry
ceiling \(\nu=1/2\).

\begin{figure}[t]
	\centering
	\includegraphics[width=0.98\linewidth]
	{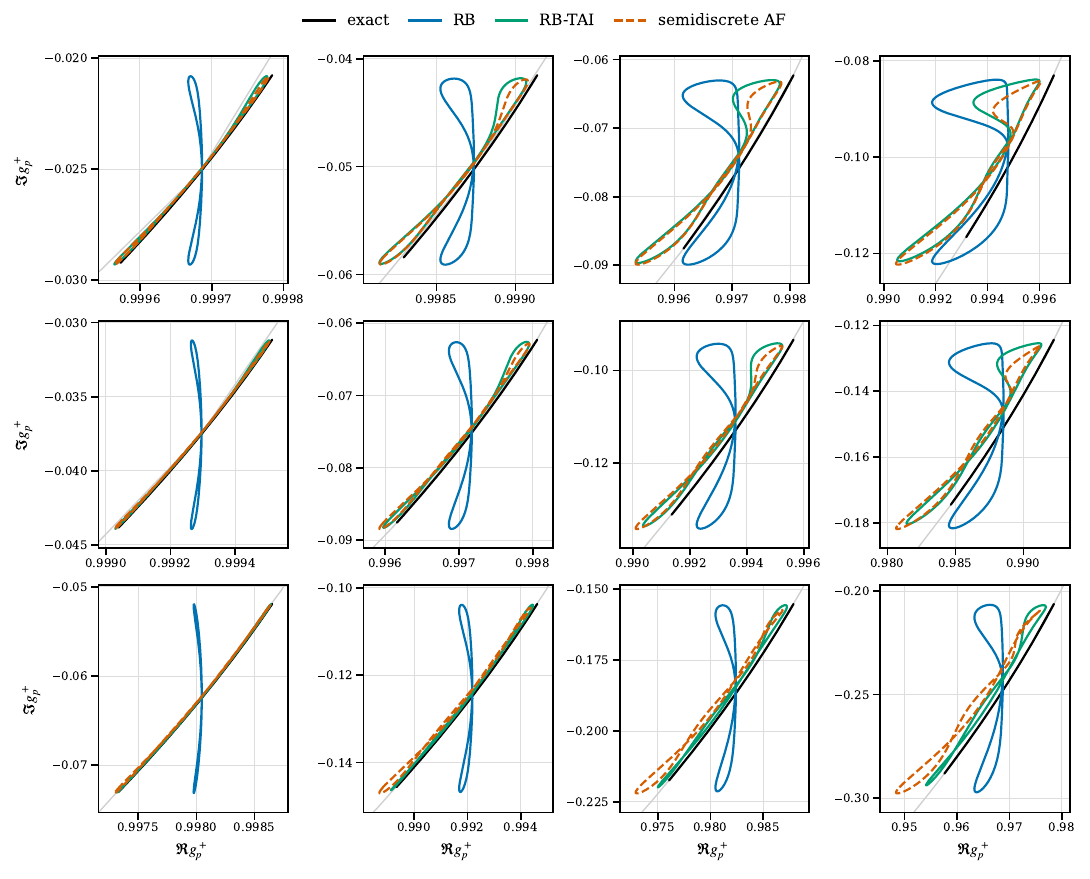}
	\caption{Complex-plane loci of the vertical-edge
		pressure response \eqref{eq:pressure-response}  as the wave-vector angle \(\vartheta\) traverses
		\([0,2\pi]\).  Rows from top to bottom use acoustic
		CFL numbers \(\nu=0.10\), \(0.15\), and \(0.25\); columns from left to right
		use \(\kappa=0.25\), \(0.50\), \(0.75\), and \(1.00\).  Curves are exact
		(black), discrete RB (blue), discrete RB-TAI (green), and semidiscrete AF
		with SSPRK3 (dashed orange); the pale curve is the unit circle.  Each panel
		is zoomed separately in its real and imaginary coordinates, so curve
		separation, rather than apparent locus aspect ratio, measures the error.
		The \(\nu=0.25\) RB row is a one-step diagnostic beyond RB's global
		stability limit.}
	\label{fig:m02-pressure-gain-loci}
\end{figure}

\begin{figure}[p]
	\centering
	\includegraphics[width=0.96\linewidth]
	{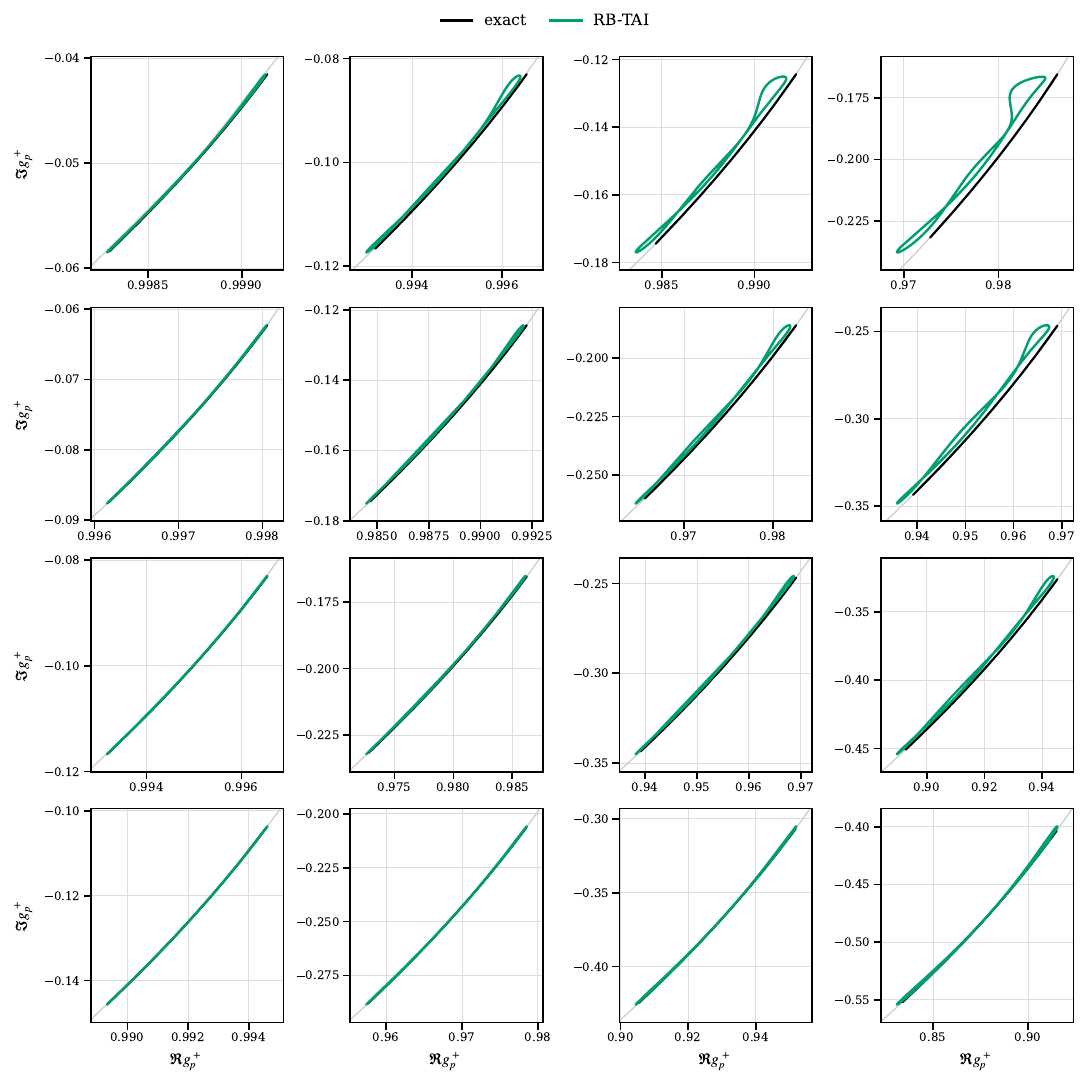}
	\caption{Continuation of Figure~\ref{fig:m02-pressure-gain-loci} for RB-TAI
		only.  Rows from top to bottom use \(\nu=0.20\), \(0.30\), \(0.40\), and
		\(0.50\); columns from left to right use \(\kappa=0.25\), \(0.50\),
		\(0.75\), and \(1.00\).  The wave-vector angle again traverses
		\([0,2\pi]\), and each complex-plane panel is zoomed independently.  Curves are exact (black) and RB-TAI (green), with
		the unit circle shown in pale gray.}
	\label{fig:m02-tai-pressure-gain-high-cfl}
\end{figure}
}

{

\section{Algorithm and Numerical Experiments}
The previous example showed good performance in the linearized Euler limit. Now, we present the full algorithm for a non-linear case, and perform a series of numerical experiments where the non-linearity is important. All tests solve the two-dimensional ideal-gas Euler equations with
\(\gamma=1.4\).   Cell averages
are initialized by quadrature, while active point values are sampled directly
from the primitive initial data. In all the experiments the CFL number is defined as $\frac{\Delta t}{\Delta x}(|\mathbf{u}|+c)_{max}.$ No limiters are used in the active flux simulations that follow. For the entropy plots, the mathematical entropy $-\int_A \rho \log{\frac{p}{\rho^\gamma}} d A$ is used. {  All of the results presented in this section correspond to an acoustic CFL $c \Delta t/\Delta x \leq 0.5$.}

\subsection{Algorithm}
Let $R_h^n W$ denote the globally continuous, cellwise primitive $Q_2$ reconstruction obtained from the time step $n$ boundary point values and the center value recovered from the conservative Simpson constraint.  

For each $
\tau\in\left\{\frac{\Delta t}{2},\Delta t\right\},
$ perform the following operations independently:

\begin{enumerate}
	\item At every vertex, vertical-edge midpoint, horizontal-edge midpoint, and recovered cell center $X_K$, set
	
	$$
	W_K^n \triangleq R_h^n W(X_K),
	\qquad
	c_K^2 \triangleq \frac{\gamma p_K^n}{\rho_K^n},
	\qquad
	Z_K \triangleq \rho_K^n c_K.
	$$
	
	Freeze the acoustic coefficients at $(\rho_K^n,c_K)$ and apply the node-centered exact acoustic operator for the duration $\tau$ to the time-$n$ reconstruction. A vertex uses four quadrant contributions, an edge midpoint uses the two adjacent half-disk contributions, and a cell center uses the full-disk contribution. Denote the returned physical variables by
	$
	u_{K,{\rm ac}}(\tau),
	v_{K,{\rm ac}}(\tau),
	p_{K,{\rm ac}}(\tau),
	$
	and form the physical acoustic increments
	$
	\Delta u_K(\tau)\triangleq u_{K,{\rm ac}}(\tau)-u_K^n,
	$
	$
	\Delta v_K(\tau) \triangleq v_{K,{\rm ac}}(\tau)-v_K^n,
	\Delta p_K(\tau) \triangleq p_{K,{\rm ac}}(\tau)-p_K^n.
	$
	
	The interpolated quantities are these physical increments, not the locally normalized acoustic variable $\pi \triangleq p/Z_K$.
	
	\item In each cell, construct the three tensor-product quadratic interpolants
	$
	I_h\Delta u(\cdot;\tau),
	I_h\Delta v(\cdot;\tau),
	I_h\Delta p(\cdot;\tau)
	$
	from the eight boundary increments and the corresponding center increment.
	
	\item For every persistent target point $P$, compute the two-Picard convective foot using only the velocity component of $R_h^nW$:
	$$
	P_f^{(1)}
	=
	P-\tau\boldsymbol{u}_h^n(P);
	P_f
	=
	P-\tau\boldsymbol{u}_h^n(P_f^{(1)}),
	$$
	
	where
	$
	\boldsymbol{u}_h^n(x)
	\triangleq
	\left(
	[R_h^nW(x)]_u,
	[R_h^nW(x)]_v
	\right)^T.
	$
	
	 Locate the cell containing $P_f$ and evaluate
	$
	W_f \triangleq R_h^nW(P_f)
	=
	(\rho_f,u_f,v_f,p_f)^T.
	$
	
	Thus both velocity samples in the foot iteration and the advected state come from the same time-$n$ primitive reconstruction. Neither half-time values nor acoustically corrected velocities are used when computing the full-time foot.
	
	\item Define the sound speed entering the density correction by
	$
	c_{P_f}^2
	\triangleq
	\frac{\gamma p_f}{\rho_f}.
	$ Evaluate the three acoustic-increment interpolants in the cell containing $P_f$:
	$$
	\Delta u_f =(I_h\Delta u)(P_f;\tau),
	\qquad
	\Delta v_f=(I_h\Delta v)(P_f;\tau),
	\qquad
	\Delta p_f=(I_h\Delta p)(P_f;\tau).
	$$
	
	The updated point value is
	
	$$
	W_P^{n+\tau}
	=
	\begin{pmatrix}
		\rho_f+\Delta p_f/c_{P_f}^2\\
		u_f+\Delta u_f\\
		v_f+\Delta v_f\\
		p_f+\Delta p_f
	\end{pmatrix}.
	$$
	
	Thus density is not interpolated as an independent acoustic increment. It is recovered at the foot from the interpolated pressure increment using the frozen acoustic invariant.
\end{enumerate}

The half-time and full-time point values are both computed directly from $R_h^nW$; the full-time update is not a second half-step initialized from the half-time state. After both point stages have been obtained, the unsplit Euler flux is evaluated at the old, half-time, and full-time boundary point values. For example, the spatial Simpson flux on a vertical edge is

$$
\widehat{F}_{i,j}^{\,\sigma}
=
\frac{1}{6}
\left[
F(W_{i,j}^{V,n+\sigma})
+
4F(W_{i,j}^{E_v,n+\sigma})
+
F(W_{i,j+1}^{V,n+\sigma})
\right],
\qquad
\sigma\in\left\{0,\frac{1}{2},1\right\}.
$$

The time-integrated edge flux is

$$
\widehat{F}_{i,j}
=
\frac{1}{6}
\left(
\widehat{F}_{i,j}^{\,0}
+
4\widehat{F}_{i,j}^{\,1/2}
+
\widehat{F}_{i,j}^{\,1}
\right),
$$

with an analogous definition of $\widehat{G}_{i,j}$ on a horizontal edge. The conservative average is then updated by

$$
\overline{U}_{ij}^{\,n+1}
=
\overline{U}_{ij}^{\,n}
-
\frac{\Delta t}{\Delta x}
\left(
\widehat{F}_{i+1,j}-\widehat{F}_{i,j}
\right)
-
\frac{\Delta t}{\Delta y}
\left(
\widehat{G}_{i,j+1}-\widehat{G}_{i,j}
\right).
$$

There is no split acoustic or advective update of the cell average and no half-time cell average. Before the next time step, the new center state and primitive $Q_2$ reconstruction are recovered from $\overline{U}^{,n+1}$ and the full-time boundary point values.

The node-centered acoustic increments may be computed before foot tracing in an implementation, but both operations depend only on $R_h^nW$. The final operation is the advected time-$n$ state plus the acoustic increment evaluated at the same foot; it is not an acoustic solve initialized by the advected state.
}

\begin{figure}
	\centering
	\includegraphics[width=0.49\textwidth]{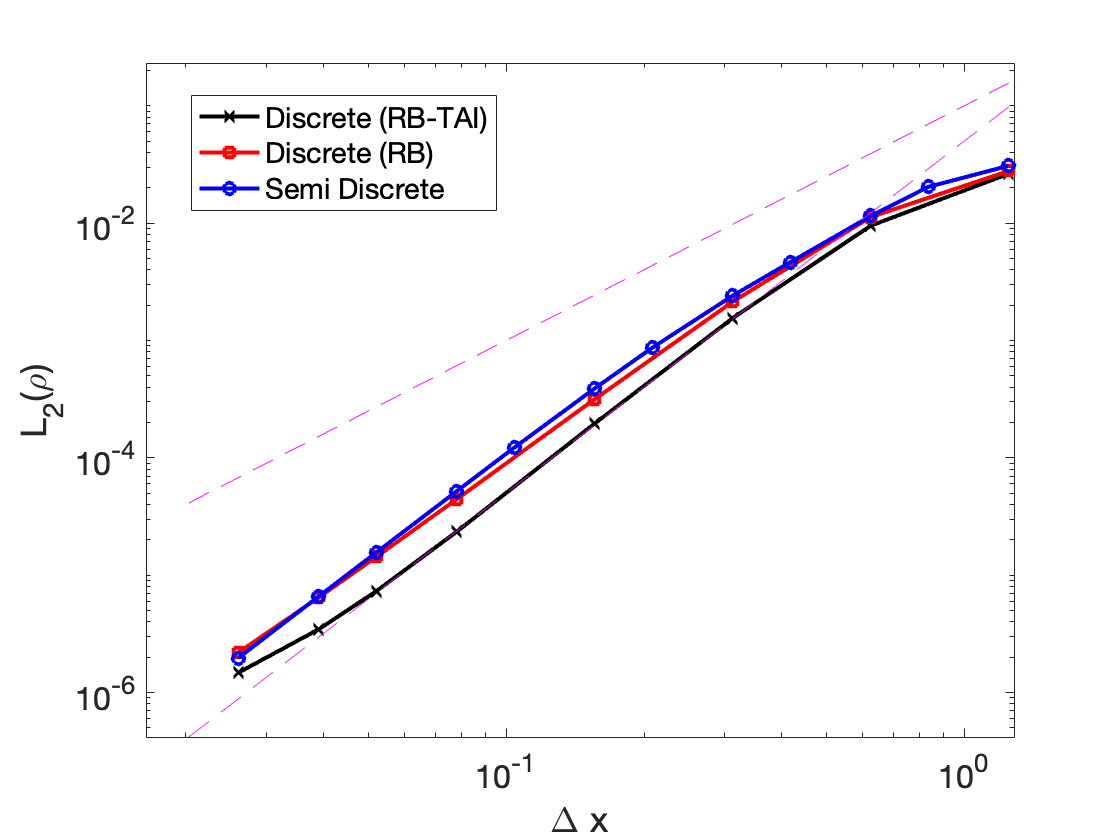}
	\includegraphics[width=0.49\textwidth]{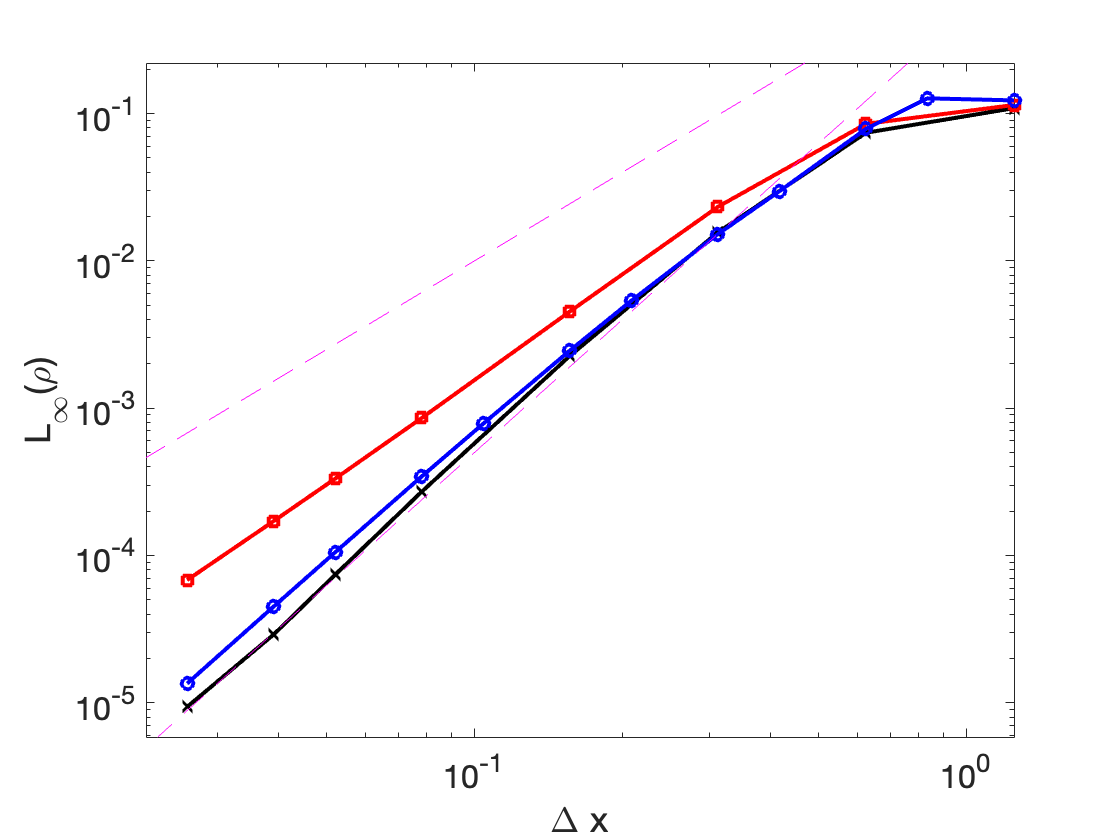}
	\caption{Isentropic vortex convection. Left: Average density error (Left: $L_2$ and Right: $L_\infty$) after 1 period with grids ranging from $8^2$ to  $384^2$. Dashed lines show 2nd and 3rd order trends.}
	\label{fig:vortex_conv}
\end{figure}

\begin{figure}
	\centering
	\includegraphics[width=0.49\textwidth]{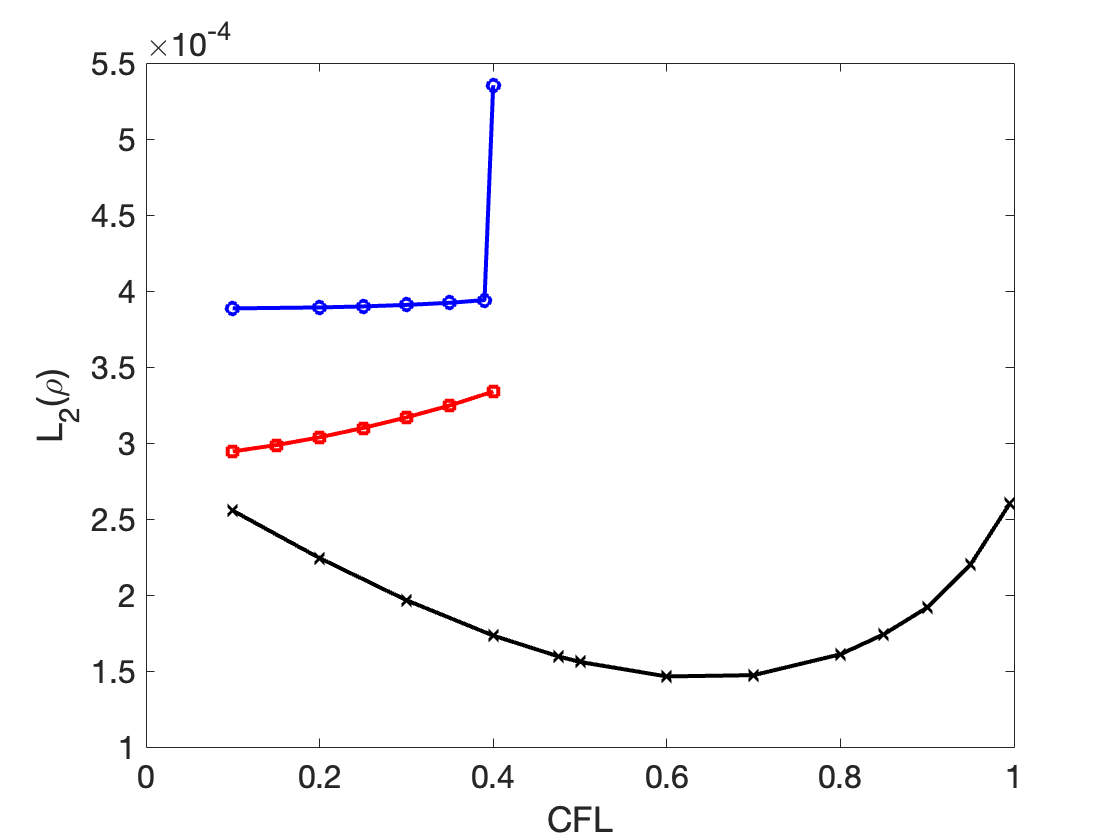}
\includegraphics[width=0.49\textwidth]{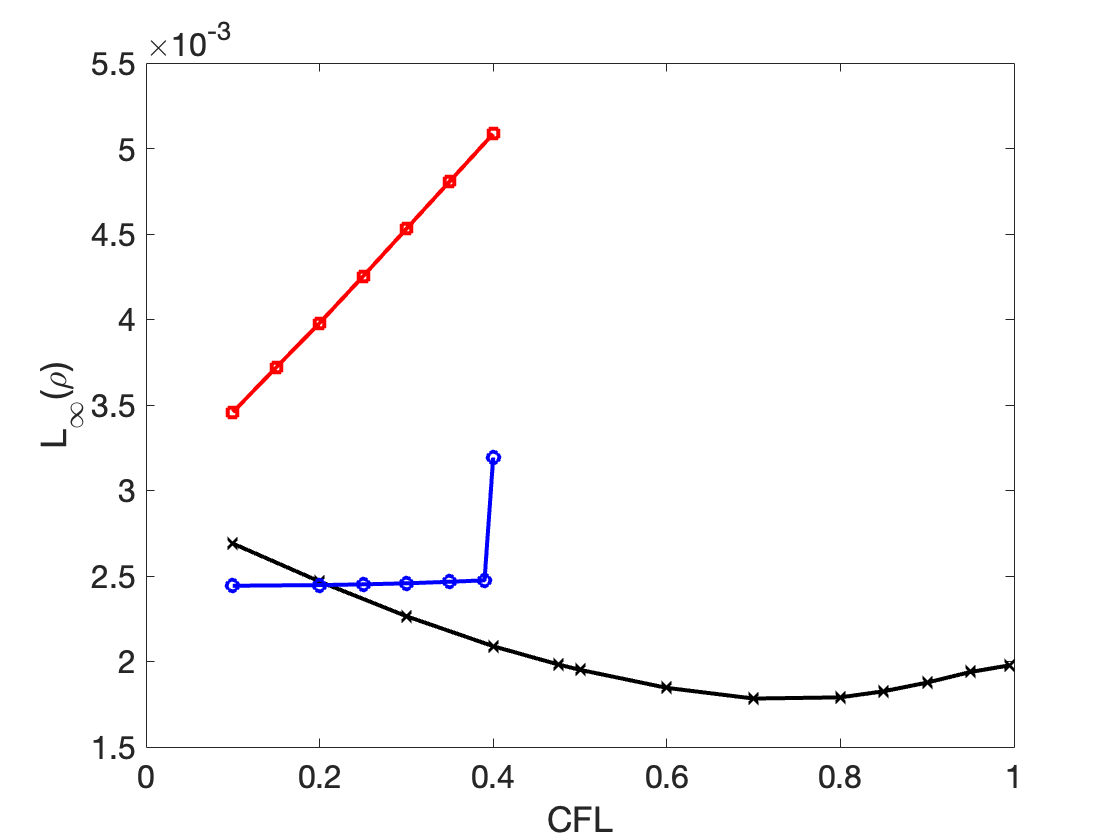}
	\caption{Isentropic vortex convection: Average density error ($L_2$ and $L_\infty$) on $64^2$ grids as a function of CFL.}
	\label{fig:vortex_conv1}
\end{figure}

\subsection{Isentropic vortex.}
The domain is \([0,10]\times[-5,5]\), with periodic boundary conditions.  The
vortex is centered initially at \((x_0,y_0)=(5,0)\) and is convected by the
background velocity \((u_\infty,v_\infty)=(1,0)\).  With \(\beta=5\),
\[
x_c(t)=x_0+u_\infty t,\qquad y_c(t)=y_0+v_\infty t,
\]
and
\[
\delta u=-{\beta\over 2\pi}\exp\left({1-r^2\over 2}\right)(y-y_c(t)),
\qquad
\delta v={\beta\over 2\pi}\exp\left({1-r^2\over 2}\right)(x-x_c(t)).
\]
The temperature-like variable is
\(
T=1-{(\gamma-1)\beta^2\over 8\gamma\pi^2}\exp(1-r^2),
\)
and the primitive state is
\(
\rho=T^{1/(\gamma-1)},\qquad
u=u_\infty+\delta u,\qquad
v=v_\infty+\delta v,\qquad
p=\rho^\gamma .
\)
A final time of 10 units is taken, which is enough for one periodic pass. The exact solution is thus the same as the initial conditon.  {The present RB-TAI scheme, the orginal RB scheme (Section~\ref{sec:barsukow}) and the semi-discrete active flux technique (Section~\ref{sec:semidiscrete}) are compared with each other on the same grids.}

{The grid-refinement study on the left of Figure~\ref{fig:vortex_conv} at a CFL number of $0.3$
shows that the RB-TAI scheme follows an approximately third order trend in $L_2$ error at practically useful resolutions ($16^2$-$128^2$) while deviating towards lower order accuracy in the asymptotic limit. At every resolution the discrete RB-TAI scheme
is more accurate than the discrete RB scheme. The exact frozen foot update (described in Section~\ref{sec:transported}) was actually found to be slightly less accurate than RB-TAI in this problem. This difference is more pronounced in the $L_\infty$ error plot.

The CFL sweep in Figure~\ref{fig:vortex_conv1} shows that the semi-discrete scheme loses stability slightly
above $\mathrm{CFL}  > 0.38$, the discrete RB scheme has a marginally higher stability region. For this problem, the present approach remains stable all the way to
$\mathrm{CFL}=1$, at which timestep the acoustic CFL is 0.39.} The structural reason is visible already at the
level of linearized analysis. In the frozen-velocity setting, where
$\bm u_0$ is constant, $S_{\mathrm{adv}}(\tau)\bm W^n(P)=\bm W^n(P_f)$, so the transported point update collapses to
\begin{equation}
	\W_{\mathrm{TAI}}^{n+\tau}(P)=
	e^{\tau C}\,e^{\tau A}\,\bm W^n(P).
	\label{eq:lie-trotter-id}
\end{equation}
That is, the transported point update is the
multiplicative composition of the two subsolvers,
whereas the additive RB update is $e^{\tau A}+e^{\tau C}-I$. Both
$e^{\tau A}$ and $e^{\tau C}$ are unitary on the linearised Euler
system; their product is therefore unitary and more stable than the sum of unitaries. We are careful to remark that the actual operator will not be unitary : the impact of interpolation,  nonlinear freezing, primitive-variable point evolution, and conservative cell-average coupling will have to be accounted for. Generally, in every problem that we have evaluated, the present approach is found to be more stable than the discrete RB scheme, yet not as dramatically so as in this problem.

A more striking feature of the right panel of
Figure~\ref{fig:vortex_conv} is the U-shape of the RB-TAI error curve. {This non-monotone behaviour
was reproduced on $16^2$, $32^2$ and $64^2$ grids, with the minima drifting to lower CFL numbers with increasing resolution, suggesting an}
intrinsic competition between two distinct error mechanisms in the
RB-TAI update, rather than a grid-specific artefact.
At high CFL, the dominant
contribution is the per-step truncation of the splitting
(beyond local linearisation, i.e.\ the residual non-commutator terms at higher order) together with the foot-point displacement
error in the nonlinear advective subsolver. 
At low CFL, the method takes many small semi-Lagrangian point updates. The accumulated interpolation/reconstruction error of the transported point values can then become visible. The acoustic-increment interpolation error is smaller in the frozen smooth analysis because the increment itself is $O(\tau)$, but nonlinear freezing and repeated primitive-variable reconstruction may produce additional low-CFL accumulation.


The flat error curves of Discrete RB and the semi-discrete scheme do not
exhibit this U-shape because their dominant error might be a result of the
\emph{additive split} cross-term defect of magnitude
$\tfrac12\tau^2(AC+CA)$ \emph{per step}, which accumulates over
$N\sim 1/\tau$ steps to a contribution of order $\tau$ in the final
time. Once the transport modification removes this dominant defect,
the residual reconstruction-accumulation and per-step truncation
errors become comparable, and the U-shape 
emerges clearly. We regard this as a useful diagnostic: the curvature
of the RB-TAI error vs. CFL is, in itself, evidence that the leading
operator-level defect of the additive split has been suppressed.

\clearpage

\begin{figure}
	\centering
	\includegraphics[width=1\textwidth]{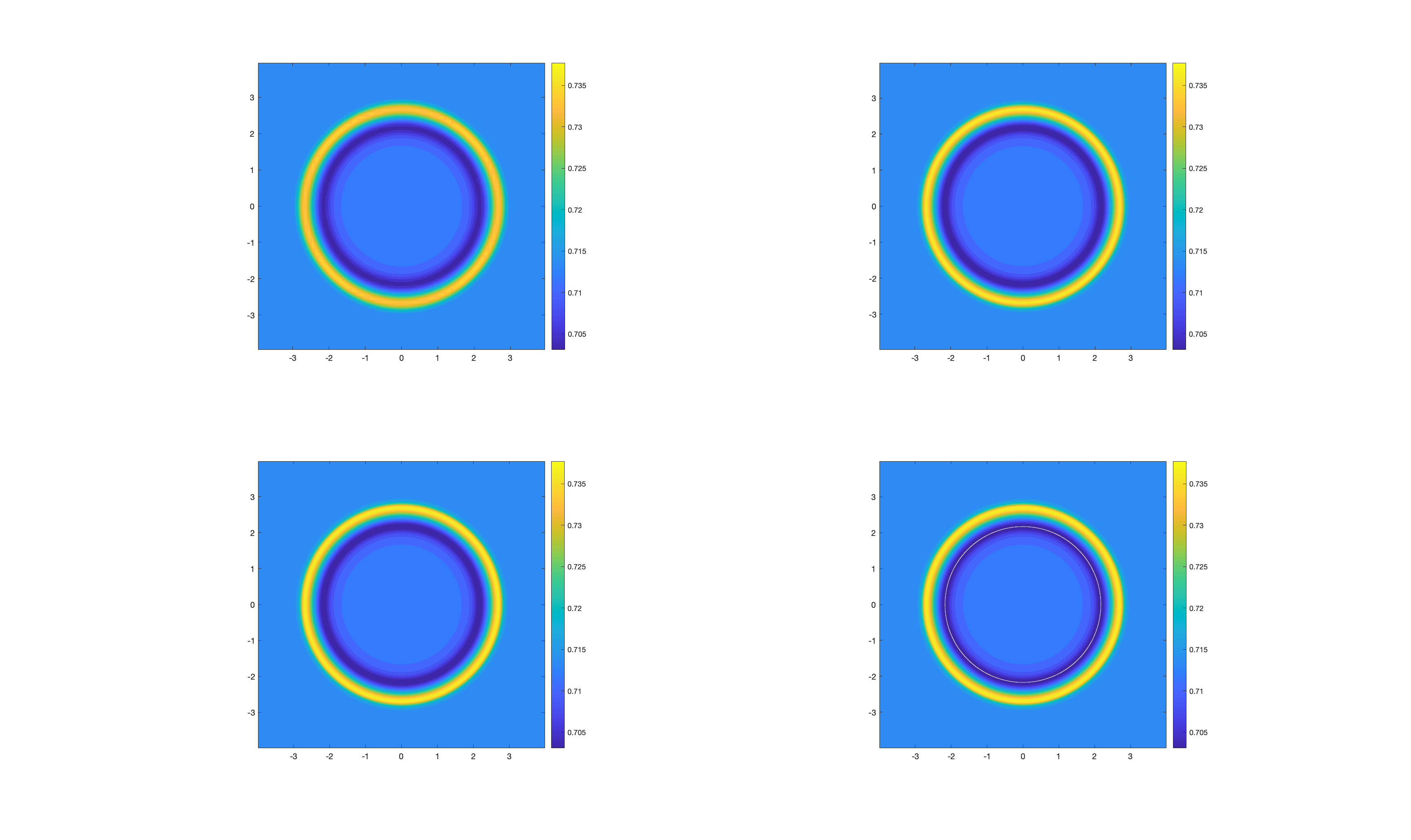}
	\includegraphics[width=1\textwidth]{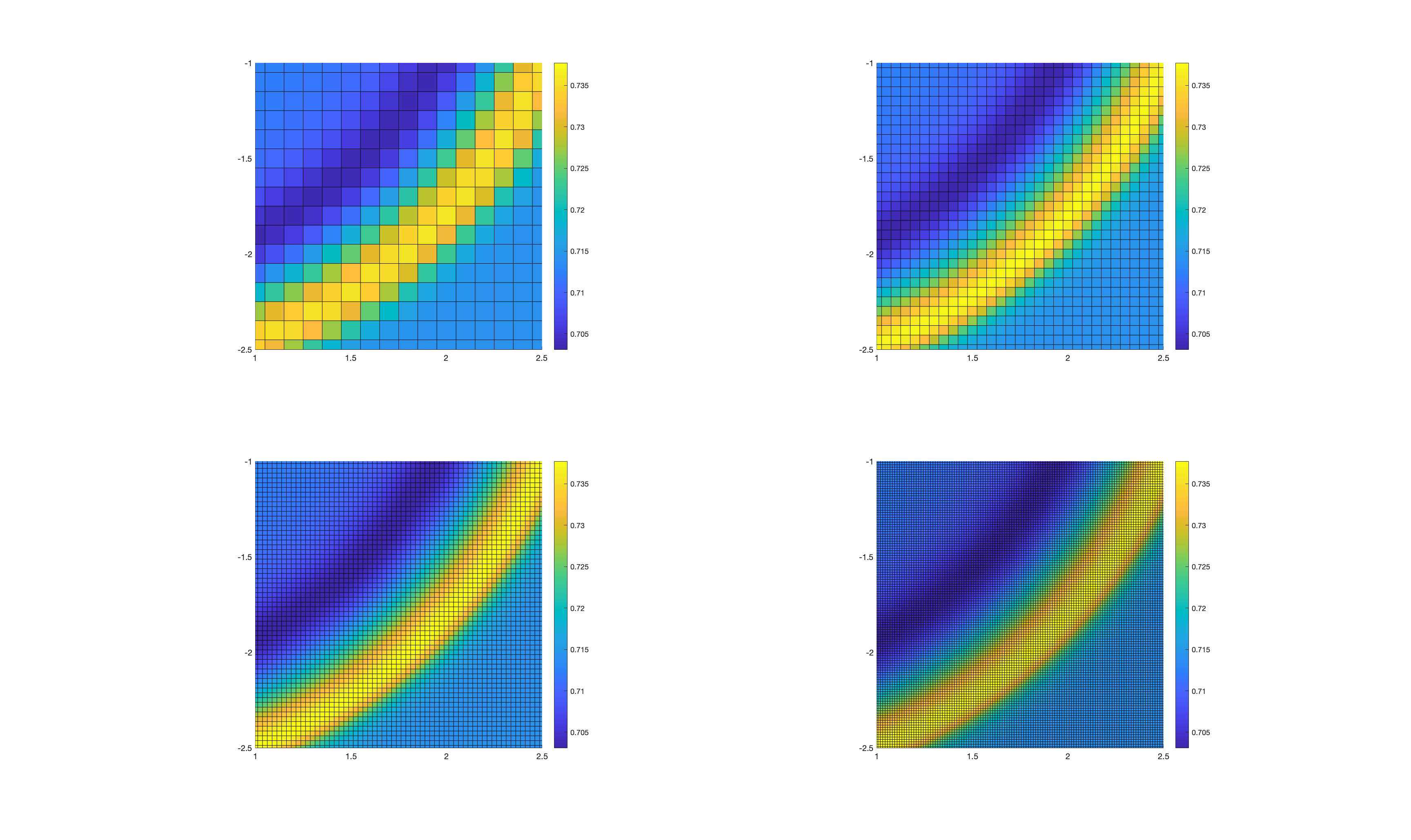}
	\caption{Evolution of Gaussian Pulse at t=2.5 at $\Delta x=0.1, 0.05,0.025,0.0125$ with RB-TAI. Cell average values shown.}
	\label{fig:gaussian_pulse}
\end{figure}

\begin{figure}
	\centering
	\includegraphics[width=0.49\textwidth]{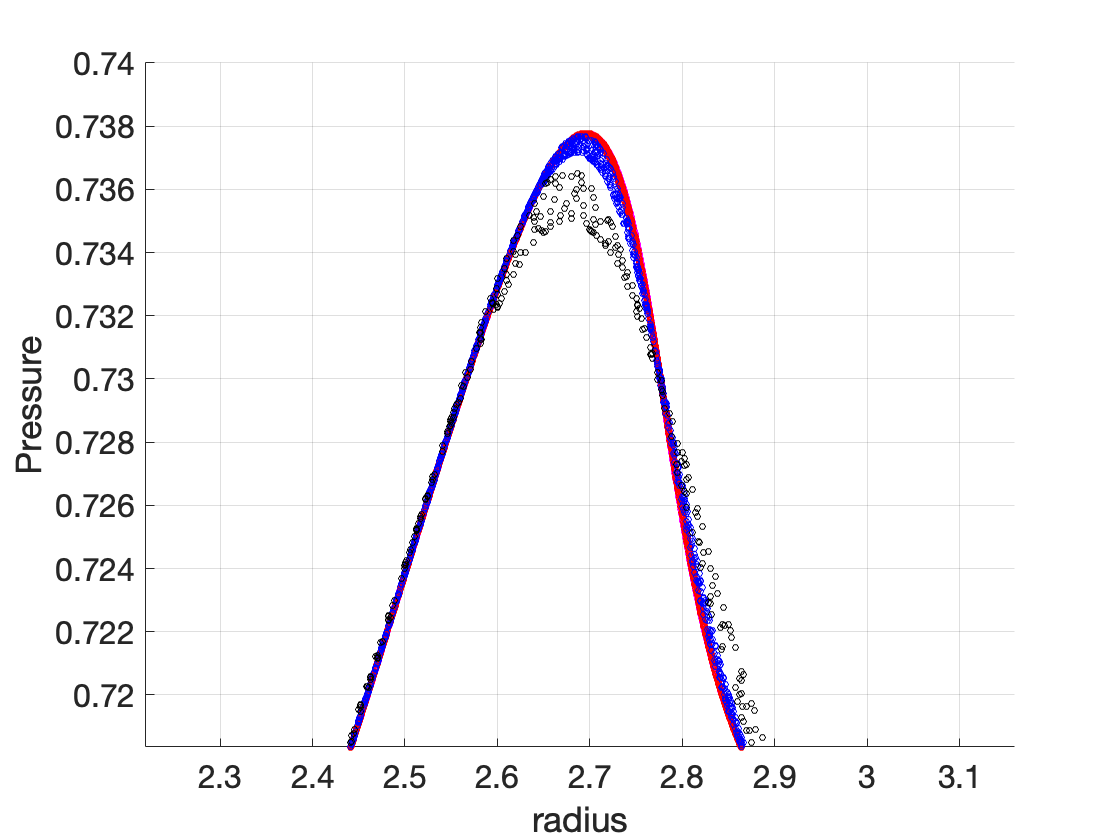}
	\includegraphics[width=0.49\textwidth]{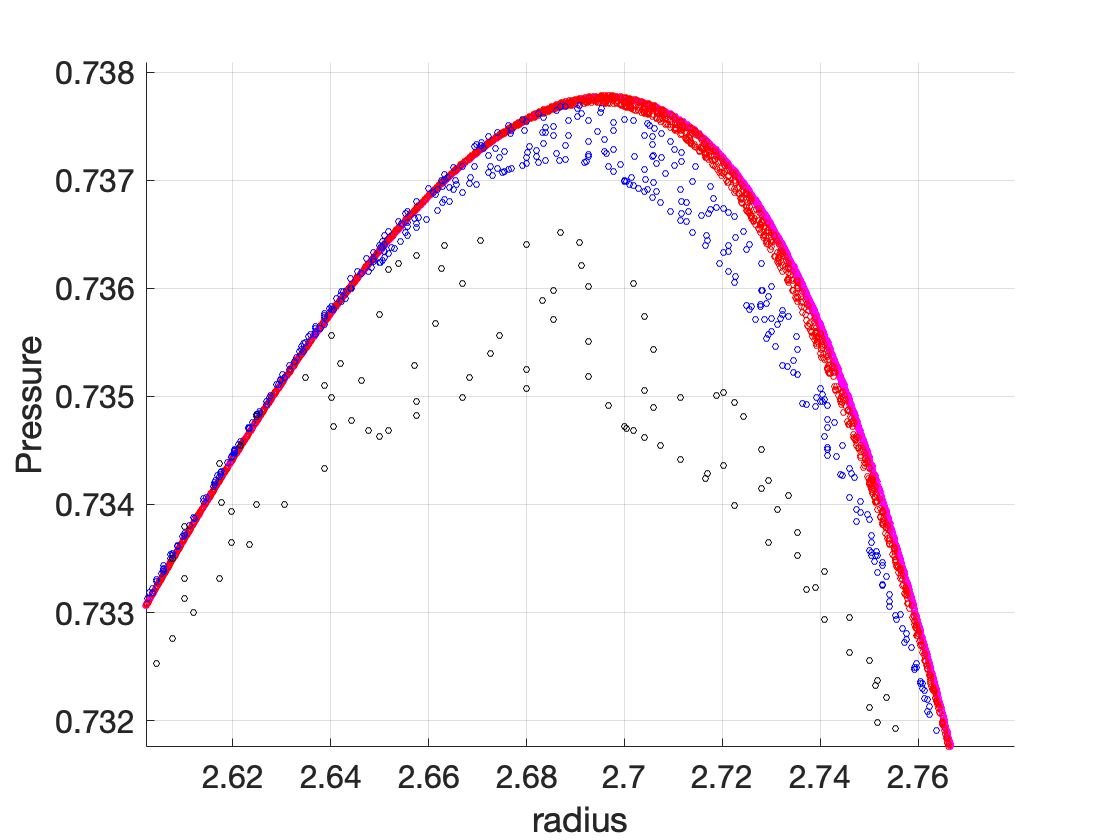}
	\caption{Evolution of Gaussian Pulse at t=2.5. Black: $\Delta x=0.1$; Blue: $\Delta x=0.05$; Red: $\Delta x=0.025$; Magenta: $\Delta x=0.0125$. Nodal values shown.}
	\label{fig:gaussian_pulse_radial}
\end{figure}

\subsection{Acoustic Pulse Evolution}
This example is adapted from the thesis work of \cite{fan2017}, where the non-linear acoustics equations were considered. Here, we modify the problem (as the original conditions lead to a shock) in the context of the Euler equations.
The domain is \([-4,4]\times[-4,4]\).  The initial condition is a Gaussian pulse.
\[
\rho=1+A\exp(-\alpha r^2),\qquad
u=0,\qquad v=0,\qquad p={1\over \gamma}\rho^\gamma .
\]
For the nonlinear Euler acoustic pulse runs, we used \(\alpha=20\) and  \(A=0.25\).  The
profiles shown in the comparison plots are taken at \(t=2.5\). { For the results presented below, time steps corresponding to $CFL_{RB}=0.2, CFL_{SD}=0.2, CFL_{RB-TAI} = 0.475$ were used. The highest resolution RB-TAI solution on a $640 \times 640$ grid is taken as the reference.}

\begin{figure}
	\centering
	\includegraphics[width=0.47\textwidth]{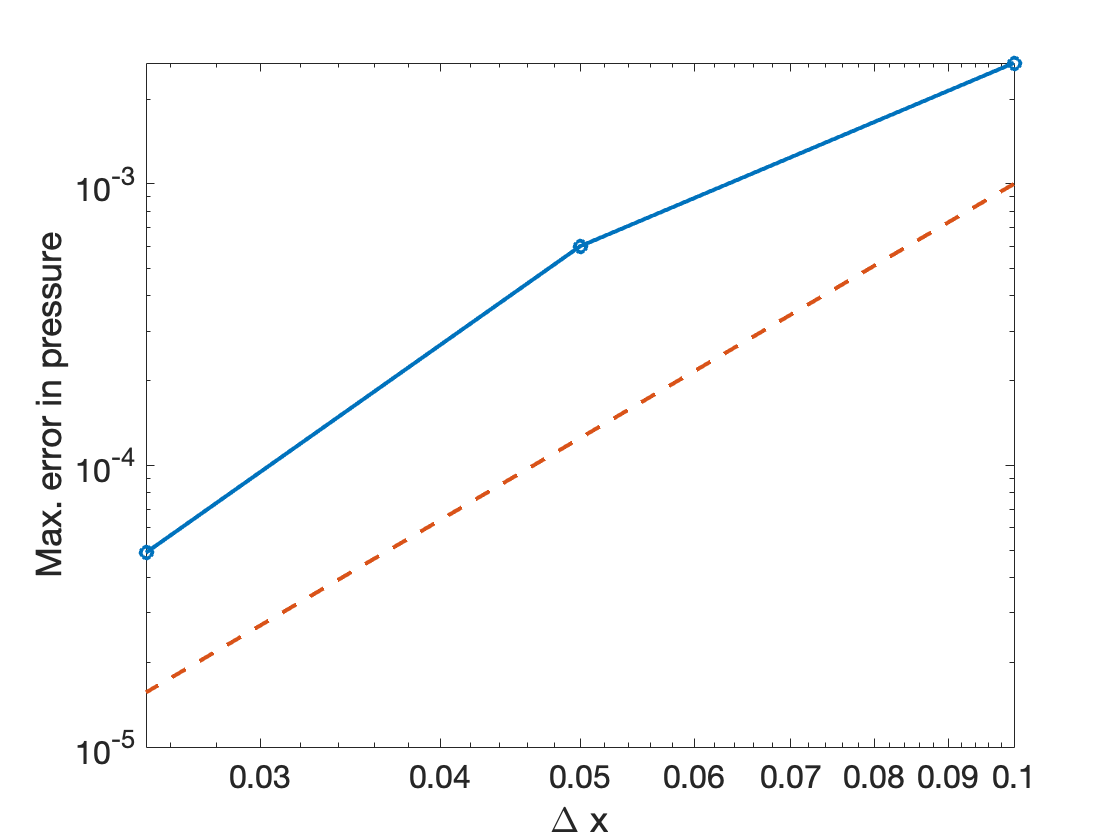}
	\includegraphics[width=0.47\textwidth]{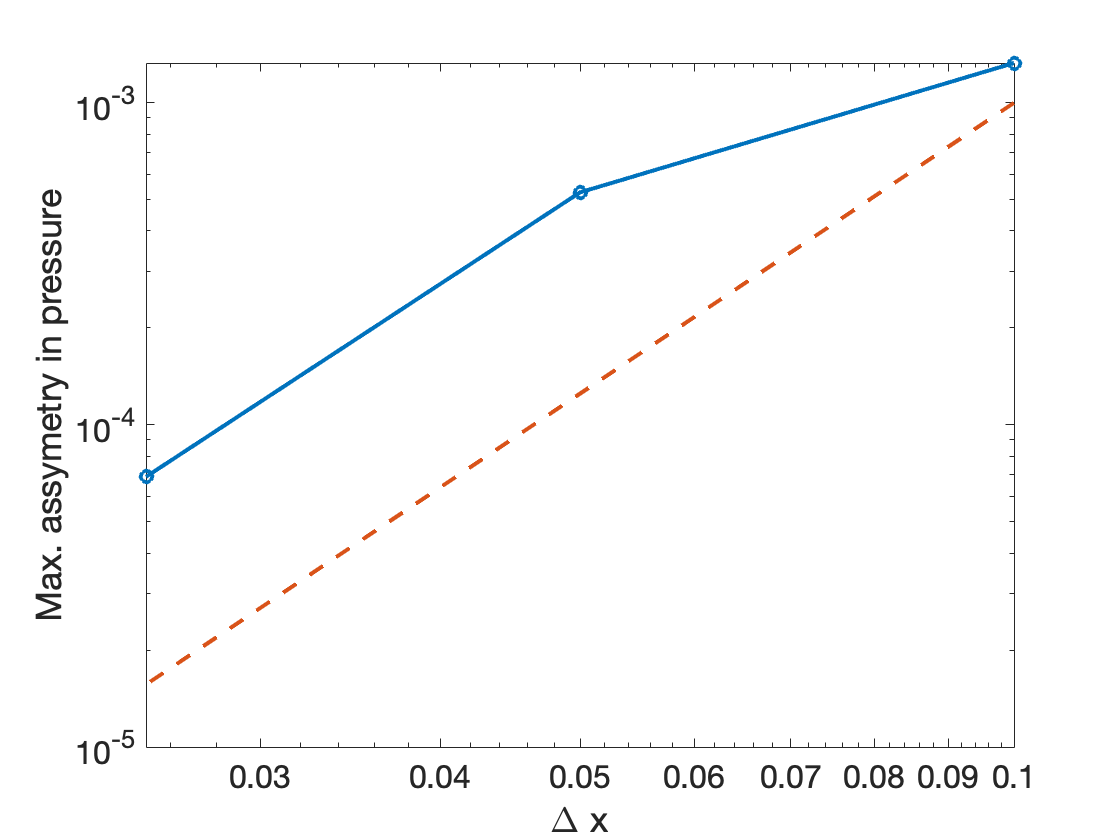}
	\caption{Maximum error and  asymmetry in pressure for  Gaussian Pulse at t=2.5.{ Dashed lines show 3rd order trend.}}
	\label{fig:gaussian_pulse_error}
\end{figure}

\begin{figure}
	\centering
	\includegraphics[width=0.5\textwidth]{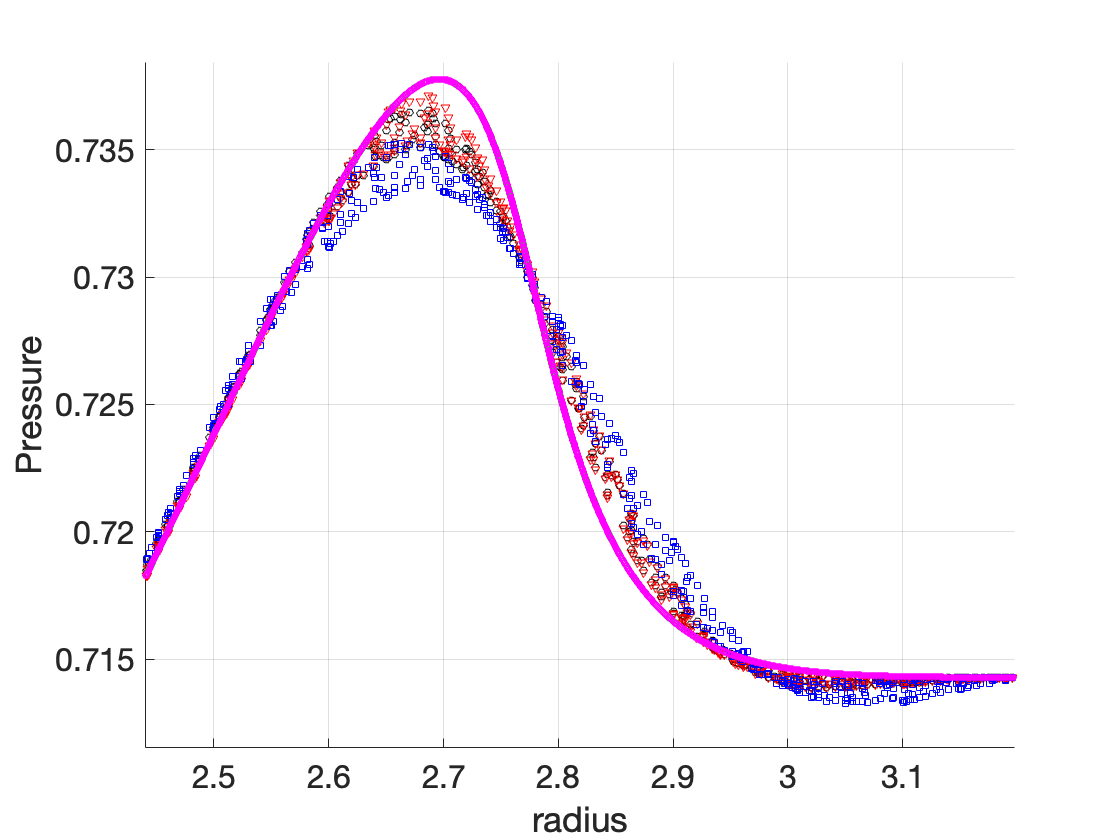}
	\caption{Zoomed in view of Gaussian pulse evolution at t=2.5, $\Delta x= 0.1$, i.e., ($40\times40$ grid).  Blue: Semi Discrete; Black: Discrete (RB-TAI); Red: Discrete (RB); Magenta: Reference. }
	\label{fig:gaussian_pulse_compare}
\end{figure}

The radial Gaussian pulse is, by construction, an acoustic-only test
in its early time evolution: the initial velocity vanishes
identically, so the advective subsolver acts only through the {
self-induced velocities} created by the pulse itself. The test
therefore primarily probes the acoustic component of each scheme.
Figure~\ref{fig:gaussian_pulse} shows the converged radial wave
structure for RB-TAI at different resolutions. 
The discrete pulse shape and amplitude approach the 
reference monotonically, and the wave fronts {display only a small asymmetry}
at every resolution: a direct visual confirmation that the
truly multidimensional acoustic operator of \cite{barsukow2022acoustic}
embedded in the scheme produces no detectable directional bias on
the underlying Cartesian grid. 

Figure~\ref{fig:gaussian_pulse_radial} plots the  nodal solutions on a radial grid. The coarse grid solutions display a small
amplitude under-prediction and a slight asymmetry across the
peak, both of which are systematically reduced under refinement;
the four curves collapse onto the high resolution reference profile as
$\Dx\to 0$. It is notable that even at the extremely coarse resolution ($\Delta x= 0.1$), {the maximum assymmetry in pressure is $< 0.15 \%$} and spread over 2 cells!

Figure~\ref{fig:gaussian_pulse_error} reports the maximum
nodal pressure error and the maximum radial asymmetry as a
function of $\Dx$. Both appear to converge close to third order. We regard the
asymmetry, in particular, as a direct probe of the multidimensional
fidelity of the acoustic operator: a dimensionally-split or
one-dimensional Riemann-solver-based scheme would generate $O(1)$
asymmetry on a Cartesian grid through axis-aligned upwinding
bias~\citep{barsukow2022acoustic}, whereas the asymmetries seen here
are pure discretisation residue and converge cleanly with the grid.

Figure~\ref{fig:gaussian_pulse_compare} contrasts the three Active
Flux variants on the coarsest grid ($\Dx=0.1$). The differences
between RB and RB-TAI are small here, as expected: in the absence of
significant background advection and at low Mach number, the convective foot $P_f$ remains very
close to the Eulerian node $P$, so the transport modification is
almost inert. The semi-discrete solution shows a noticeably broader
and more diffusive pulse, reflecting  additional numerical damping. This is consistent with the interpretation of the
isentropic vortex test: the principal advantage of the transported
discrete scheme over the semi-discrete one lies in the elimination of
multistage diffusion and in the multiplicative composition of
information transport.

{Figure~\ref{fig:gaussian_pulse_compare} also shows that the undershoot near radius $\approx$ 3.05 is less pronounced for the RB-TAI scheme. While this is promising, the fact that the recovery of the cell-center state is not a convex combination of the cell average and boundary point values has to be given consideration in challenging problems that demand the design of positivity preservation strategies.}

\clearpage

\begin{figure}
	\centering
	\includegraphics[width=0.6\textwidth]{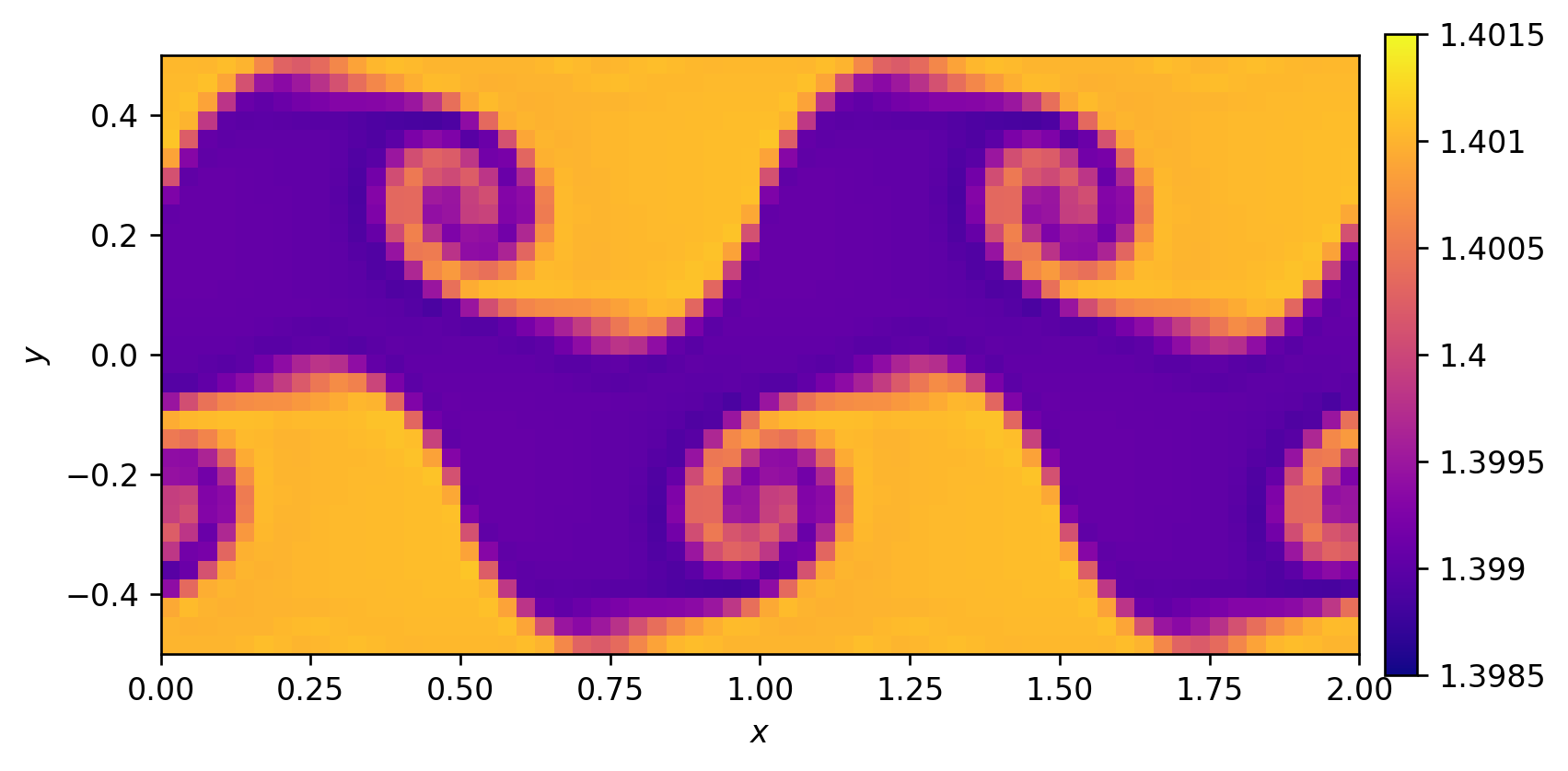}
	\includegraphics[width=0.6\textwidth]{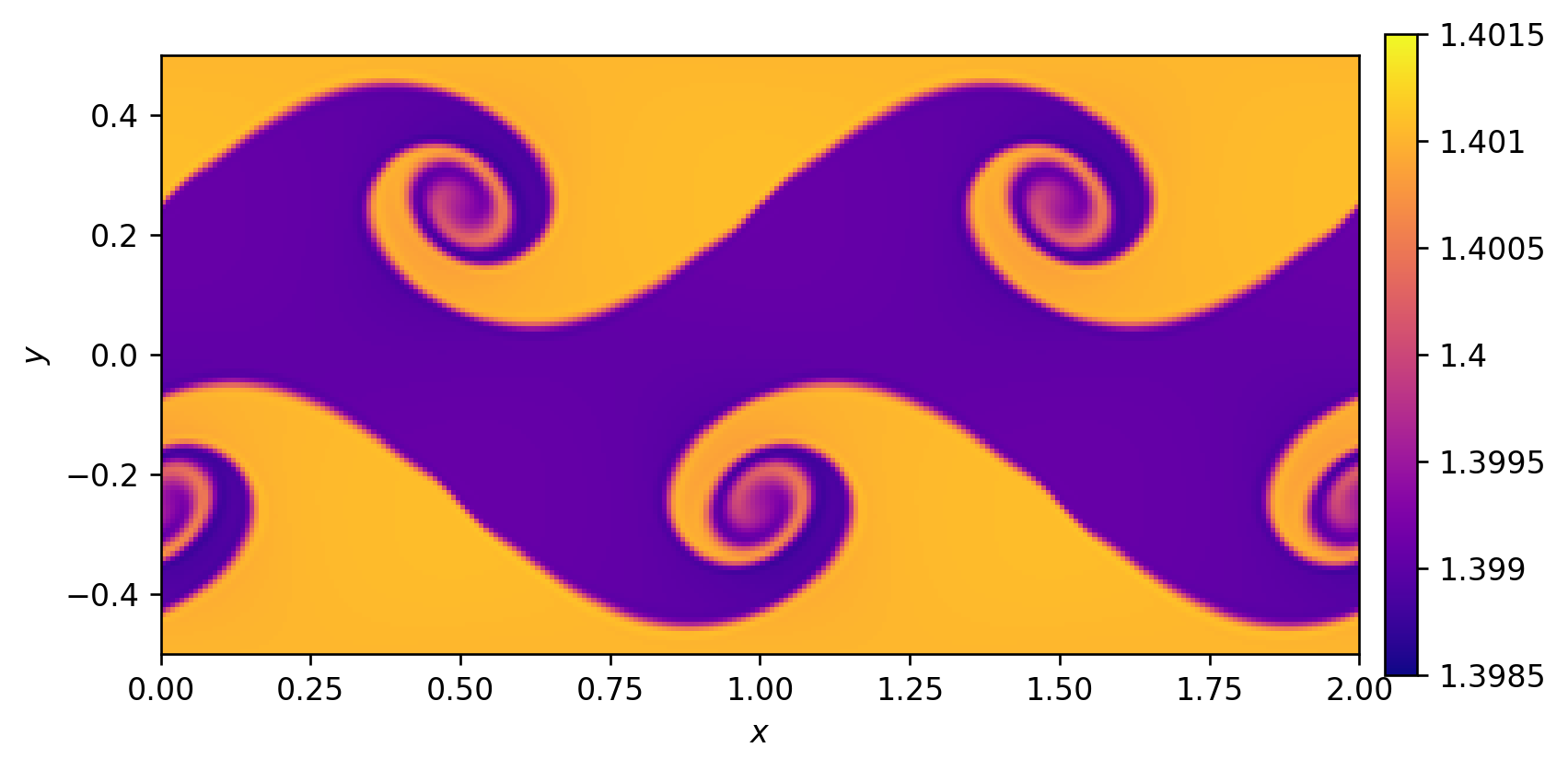}
	\caption{Density contours for shear problem on $64\times 32$ and $256\times 128$ grids using Discrete (RB-TAI) }
	\label{fig:sheardens}
\end{figure}

\begin{figure}
	\centering
	\includegraphics[width=0.49\textwidth]{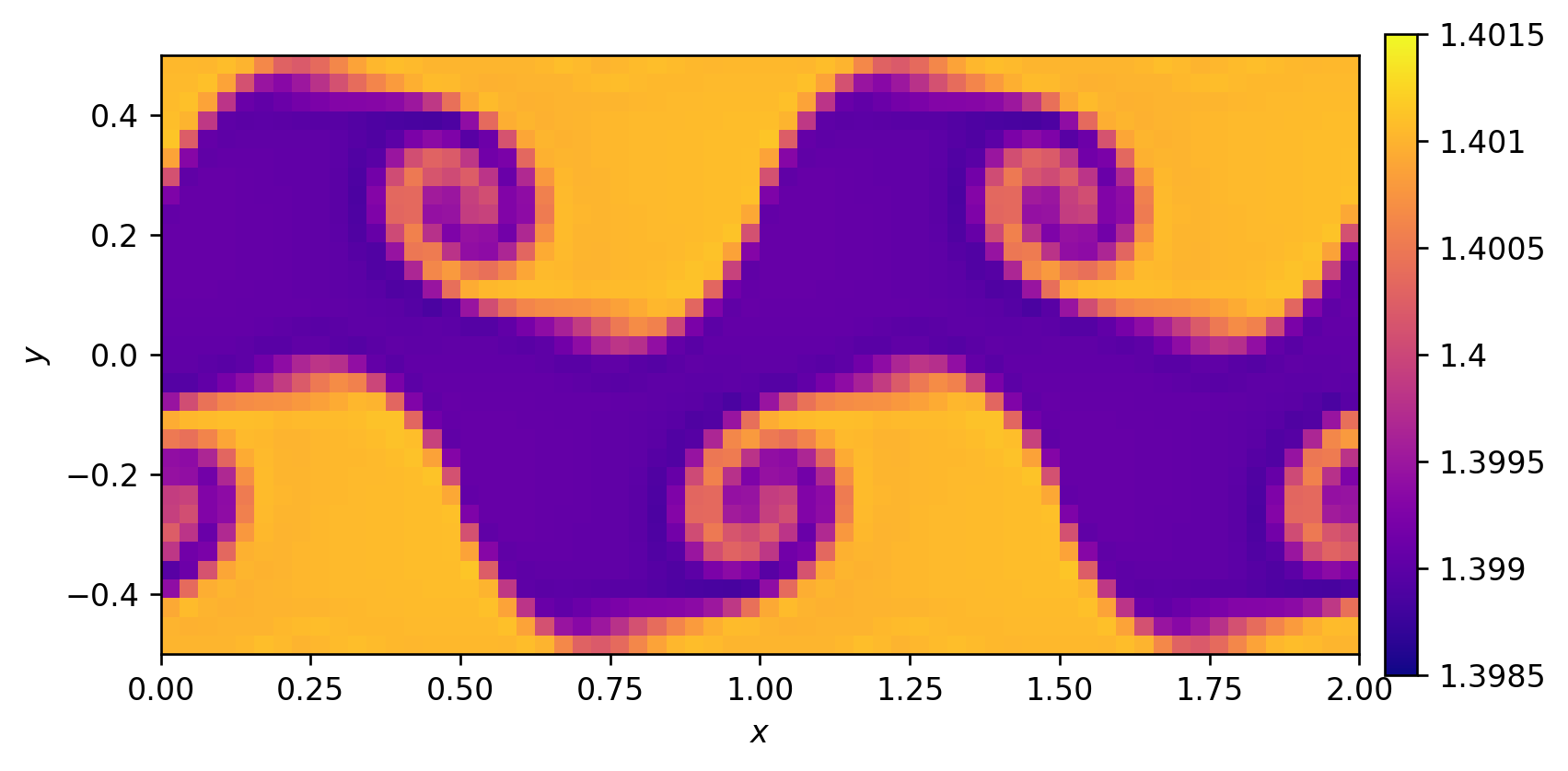}
	\includegraphics[width=0.49\textwidth]{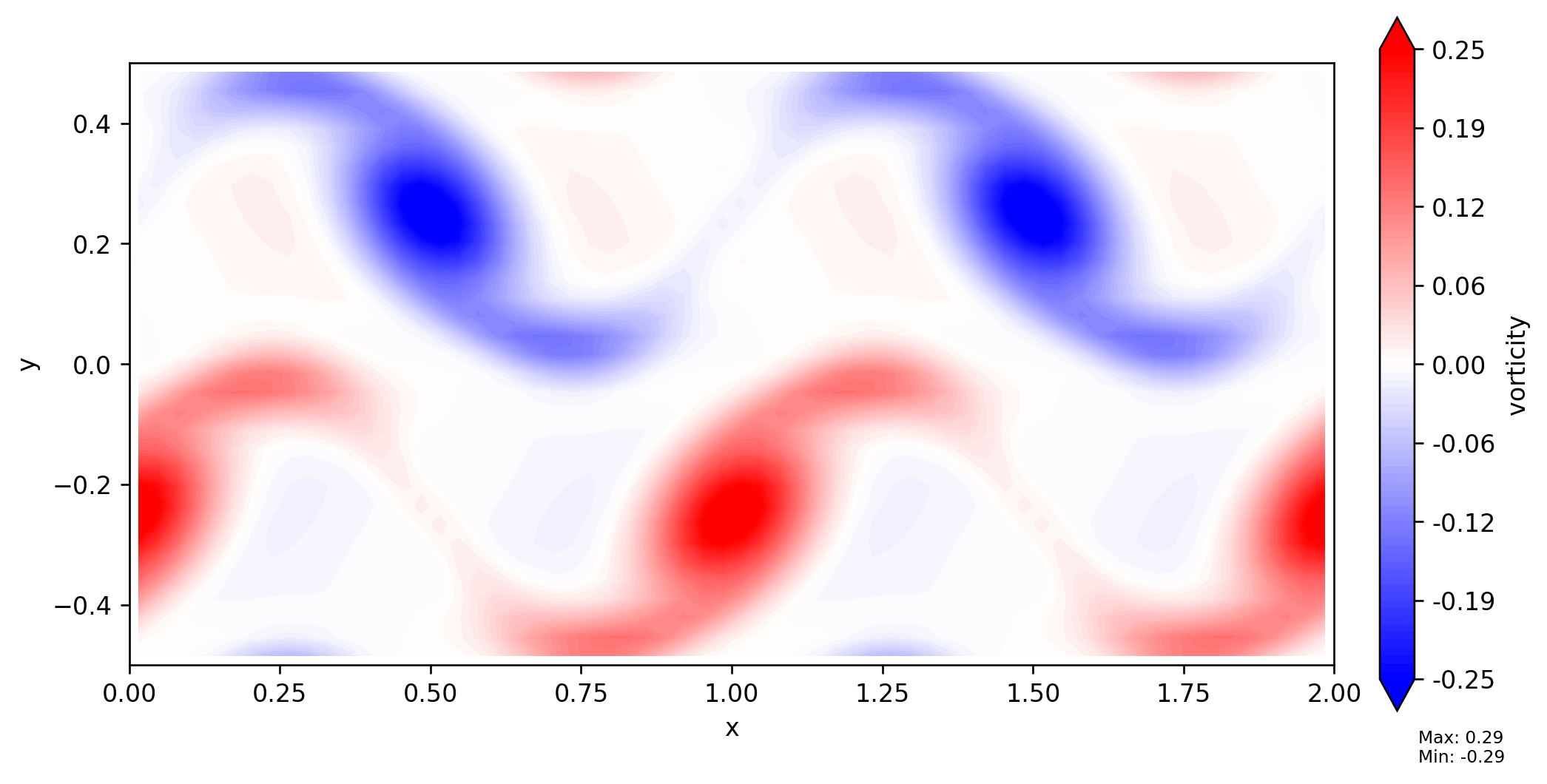}
	\includegraphics[width=0.49\textwidth]{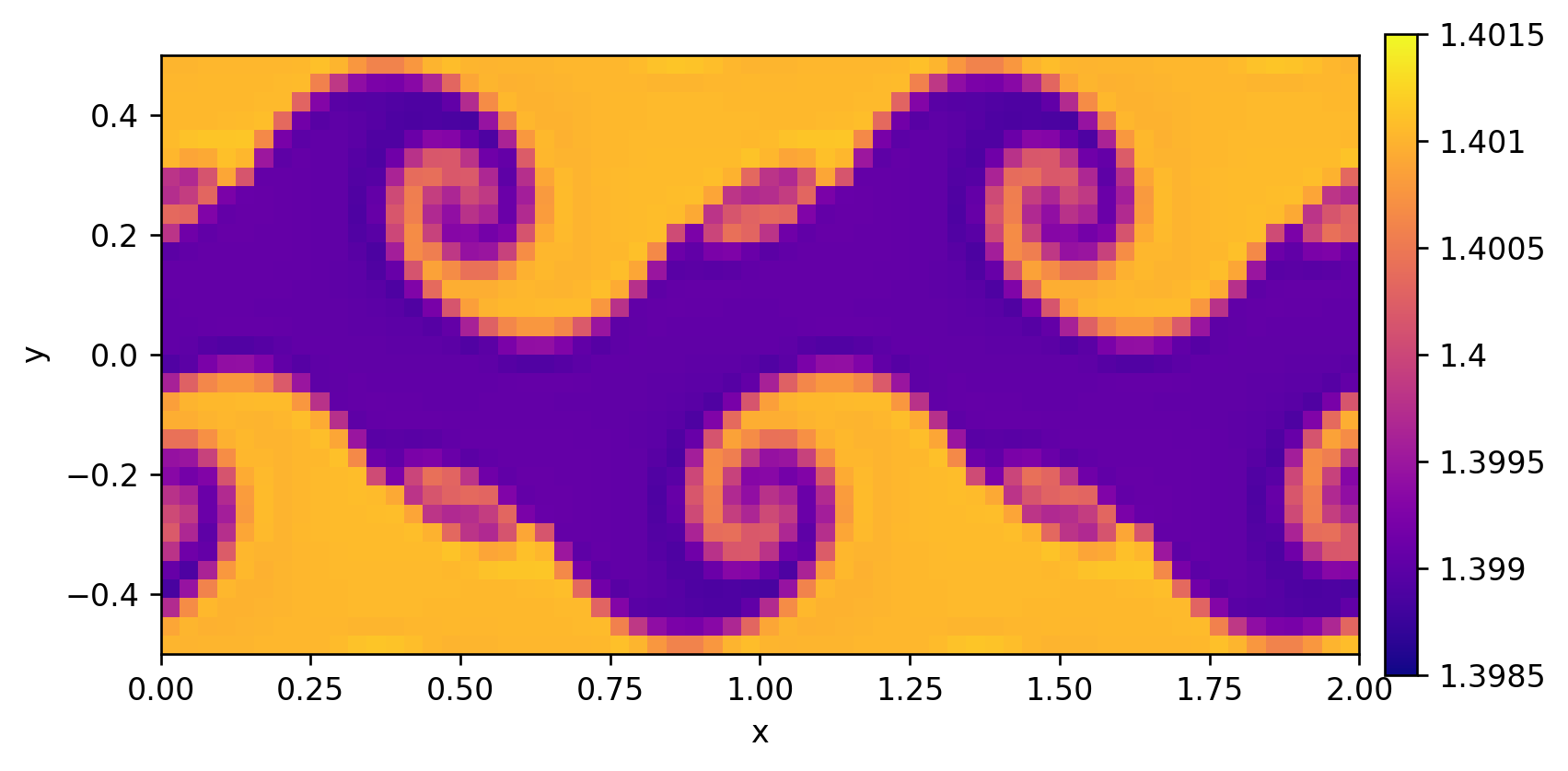}
	\includegraphics[width=0.49\textwidth]{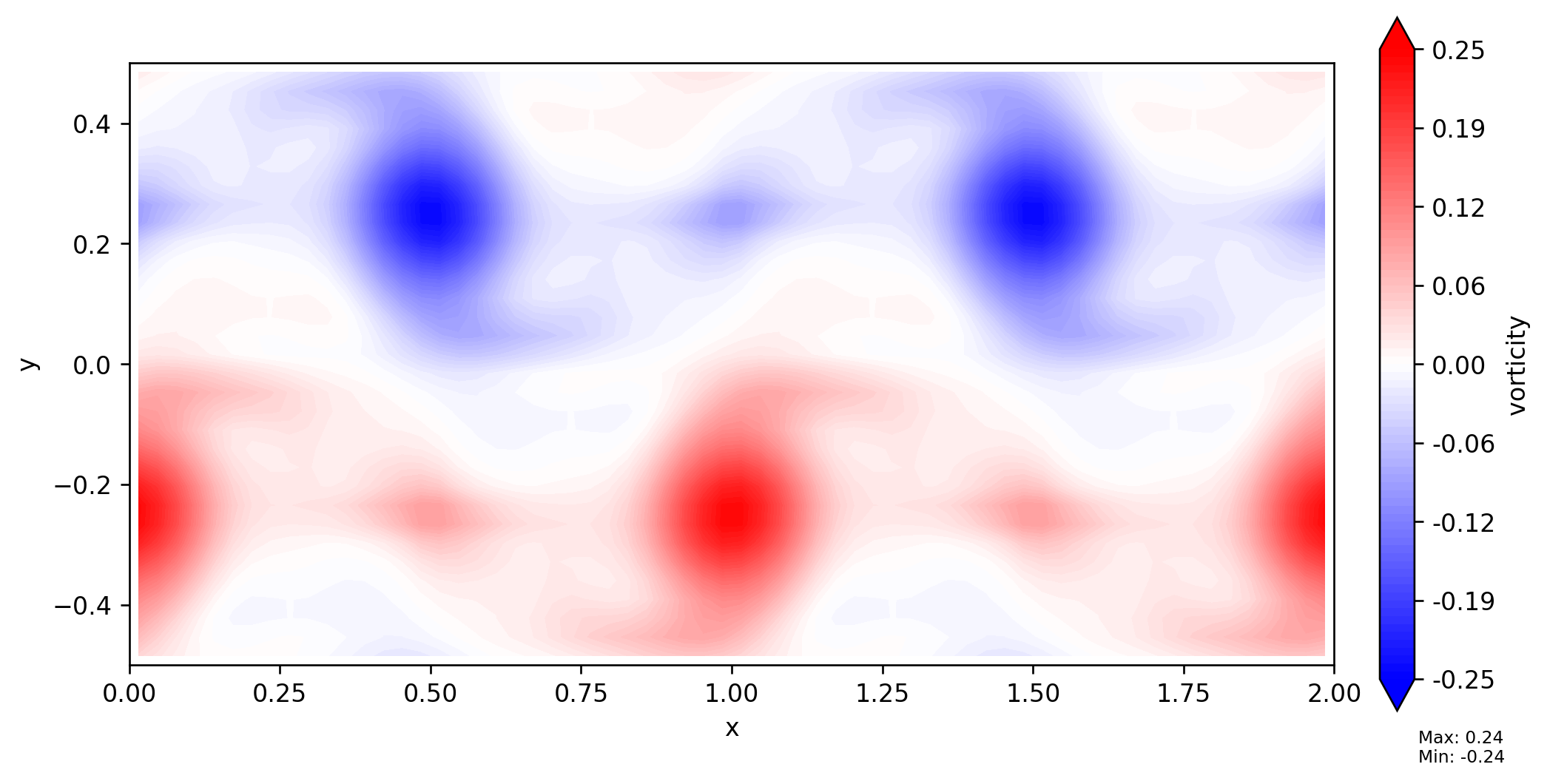}
	\caption{Shear Problem on $64\times 32$ grids. Top: Discrete (RB-TAI), Bottom: DG (P1) }
	\label{fig:shearcomp}
\end{figure}

\subsection{Low Mach Number Shear }

This test case~\cite{leidi2024performance, abgrall2025semieuler} involves the perturbation of a shear layer.
The domain is \([0,2]\times[-0.5,0.5]\), with periodic boundary conditions. 
Define
\[
a(y)=
\begin{cases}
	{1\over 2}\left[1+\sin(16\pi(y+0.25))\right],
	& -9/32 \le y < -7/32,\\
	1, & -7/32 \le y < 7/32,\\
	{1\over 2}\left[1-\sin(16\pi(y-0.25))\right],
	& 7/32 \le y < 9/32,\\
	0, & \text{otherwise}.
\end{cases}
\]
The primitive initial data are
\[
\rho=\gamma+R(1-2a(y)),\qquad
u=M(1-2a(y)),\qquad
v=\delta M\sin(2\pi x),\qquad
p=1 .
\]
Here
\(
R=10^{-3},\qquad M=10^{-2},\qquad \delta=0.1 .
\)
The standard shear-layer visualizations are shown at \(t=80\) in Figure~\ref{fig:sheardens}. The RB-TAI solution on a very coarse ($64 \times 32$) mesh appears to track the correct features: primary vortices with a clean wraparound
between them, when compared to a  refined mesh ($256 \times 128$) solution ; no secondary structures are present. Further confirmation of this reference solution can be seen in Appendix ~\ref{sec:Addplot} on higher resolution meshes with a Continuous Galerkin scheme. For this problem, as the Mach number is very small, the difference between the RB and RB-TAI methods was minimal, except for the fact that RB-TAI was stable up to a CFL number of 0.75, whereas the RB scheme was stable until a CFL number of 0.5. The Semi Discrete method was stable up to a CFL number of 0.3. { We note that, with the baseline implementation of $S_{ac}$ at the edge midpoints, it is not recommended to operate beyond an acoustic CFL of 0.5. For the results presented below, time steps corresponding to $CFL_{RB}=0.2, CFL_{SD}=0.2, CFL_{RB-TAI} = 0.475$ were used.}

Figure~\ref{fig:shearcomp} reports coarse mesh
solutions of the same problem produced by a P1 DG method with the same number of degrees of freedom. The DG solution  develops a spurious secondary
vortex pair between the primary rollers, accompanied by ringing
along the shear interface. These features are possibly artefacts
of one-dimensional Riemann-solver-based upwinding and low Mach numbers,
where the numerical pressure-velocity coupling can pollute the solution anisotropically~\citep{guillard1999,dellacherie2010}.

{The fact that the secondary structures disappear under mesh refinement in the DG and CG methods suggests that the secondary vortices are numerical artifacts, and that the multidimensional-transport character of the Active Flux methods are beneficial.}
 
More details of the DG and CG method are provided in Appendix \ref{sec:Addplot}. We are careful to remark that this is {\em not} meant to be a direct comparison of the accuracy of Active flux vs DG and CG methods. That exercise requires more careful baselining (e.g.~\citep{roe2018comparing,barsukow2026equivalence,barsukow2025semi}), and the present problem is not an ideal test case to evaluate DG methods. Rather the present example reinforces  low-Mach number characteristics of Active Flux methods \citep{barsukow2019cartesian,barsukow2021allspeed} without any adhoc fixes. {This test also suggests that the absence of secondary structure in the coarse-grid results reflects the benefit of  multidimensional information transport built into the method, rather than a coincidence
of resolution.}

The vorticity time history (Figure~\ref{fig:shearvort} in Appendix ~\ref{sec:Addplot}) shows that the symmetric, alternating-sign shear sheets are preserved as coherent structures throughout the run to $t=80$, with the maximum vorticity drifting
only modestly between $t=20$  and $t=80$. For a method run at $\mathcal{M}=10^{-2}$ on
a barely-resolved mesh over $\sim 8\times 10^3$ time steps, this is
the {structural behaviour} expected of a scheme that
propagates  genuine multidimensional acoustic information without
preconditioning.

\begin{figure}
	\centering
	\includegraphics[width=0.49\textwidth]{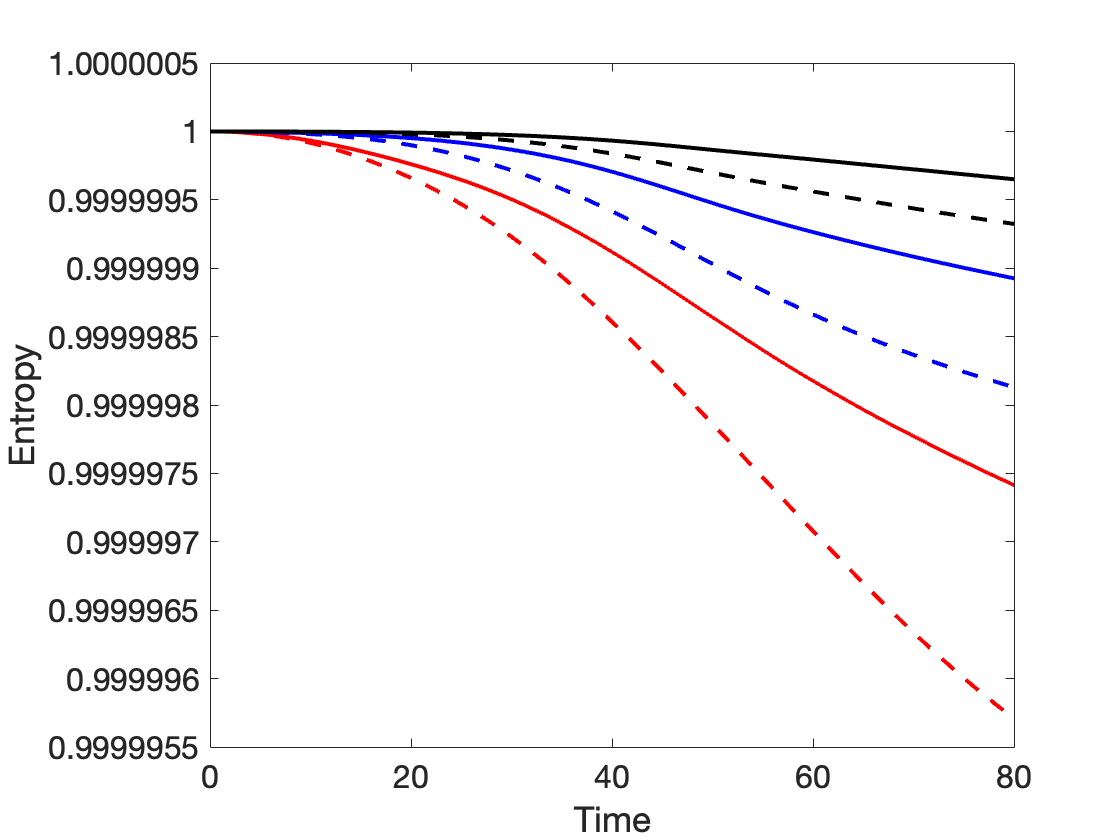}
	\includegraphics[width=0.49\textwidth]{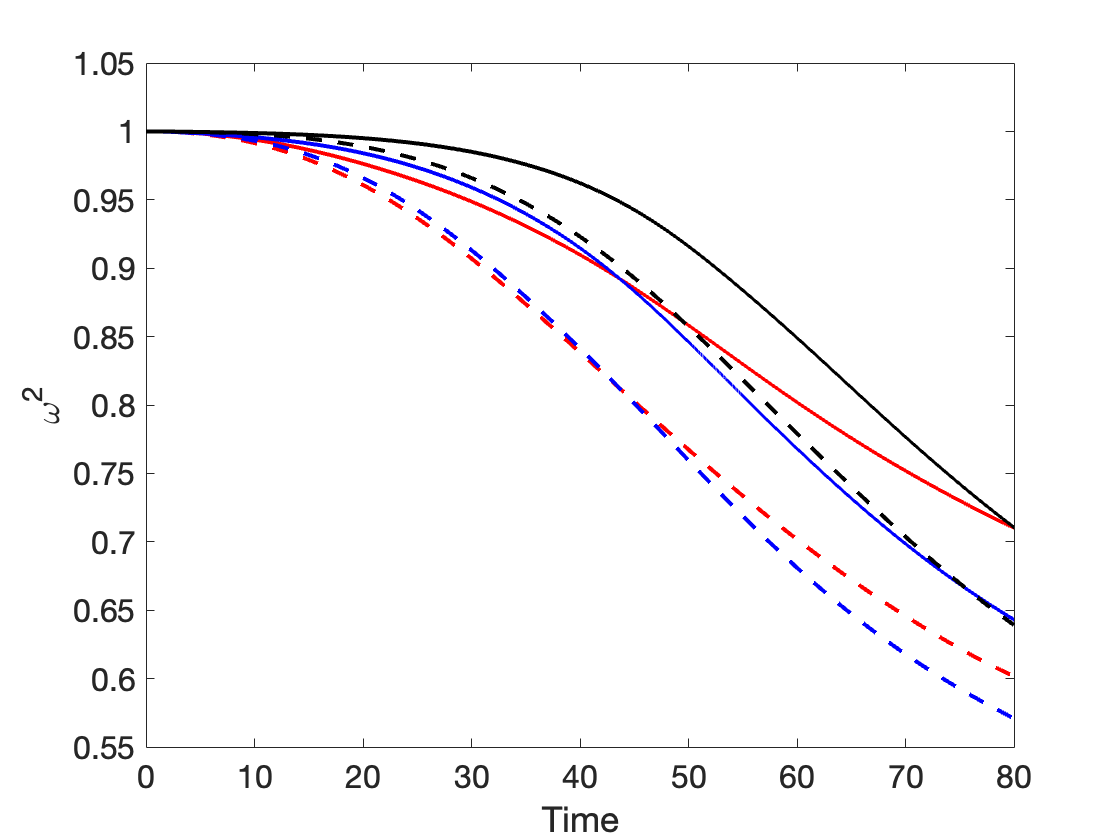}
		\caption{Normalized integral quantities. Red: $64 \times 32$; Blue: $128 \times 64$; Black: $256 \times 128$. Solid lines: Discrete (RB-TAI); Dashed lines: Semi Discrete.}
	\label{fig:shearstats}
\end{figure}

Figure~\ref{fig:shearstats} shows {a  low level of} entropy production even on the coarsest mesh, and also highlights the fact that the discrete active flux technique is less dissipative than the Semi Discrete method. Further, it was confirmed that for all three active flux variants,{ the integral of vorticity was at machine precision $1 \times 10^{-18}$ throughout the duration of the simulation, which amounts to 32,000 time steps in the finest mesh.} The integral of the square of the vorticity, however, does show some decline (Figure~\ref{fig:shearstats}), with the discrete versions of the method showing smaller declines.

\subsection{Robustness test: Compressible, highly unstable Kelvin Helmholtz configuration}
The next example involves a compressible and highly unstable shear problem which has been used to stress test high order methods in very recent publications~\citep{glaubitz2025generalized,ranocha2025robustness,bercik2026stable}. 
The computational domain is \([-1,1]\times[-1,1]\) with periodic boundary
conditions.  The initial condition is 
\[
B(y)=\tanh(15y+7.5)-\tanh(15y-7.5),
\]
\[
\rho(x,y,0)=\frac12+\frac34 B(y),\qquad
u(x,y,0)=\frac12\bigl(B(y)-1\bigr),
\]
\[
v(x,y,0)=\frac{1}{10}\sin(2\pi x),\qquad
p(x,y,0)=1 .
\]
The goal is to evolve the solution to \(t=15\). Most of the numerical methods in the above publications~\citep{glaubitz2025generalized,ranocha2025robustness,bercik2026stable} fail by $t < 5$ units. All three active flux variants  spanning under-resolved meshes ($64^2, 128^2, 256^2$) demonstrated robustness.

 The Semi Discrete method required a CFL $< 0.2$, whereas the discrete methods were stable for CFL $<0.5$, again with the RB-TAI variant being marginally more stable. { Time steps corresponding to $CFL_{RB}=0.25, CFL_{SD}=0.15, CFL_{RB-TAI} = 0.475$ were used.}
  In contrast to the prior low Mach number test case, there were noticeable differences between  RB and RB-TAI, as the latter showed noticeably sharper features as can be seen from Figure~\ref{fig:shearnew64}. Further evidence of this is shown in Appendix~\ref{sec:Addplot}. 

\begin{figure}
	\centering
	\includegraphics[width=0.49\textwidth]{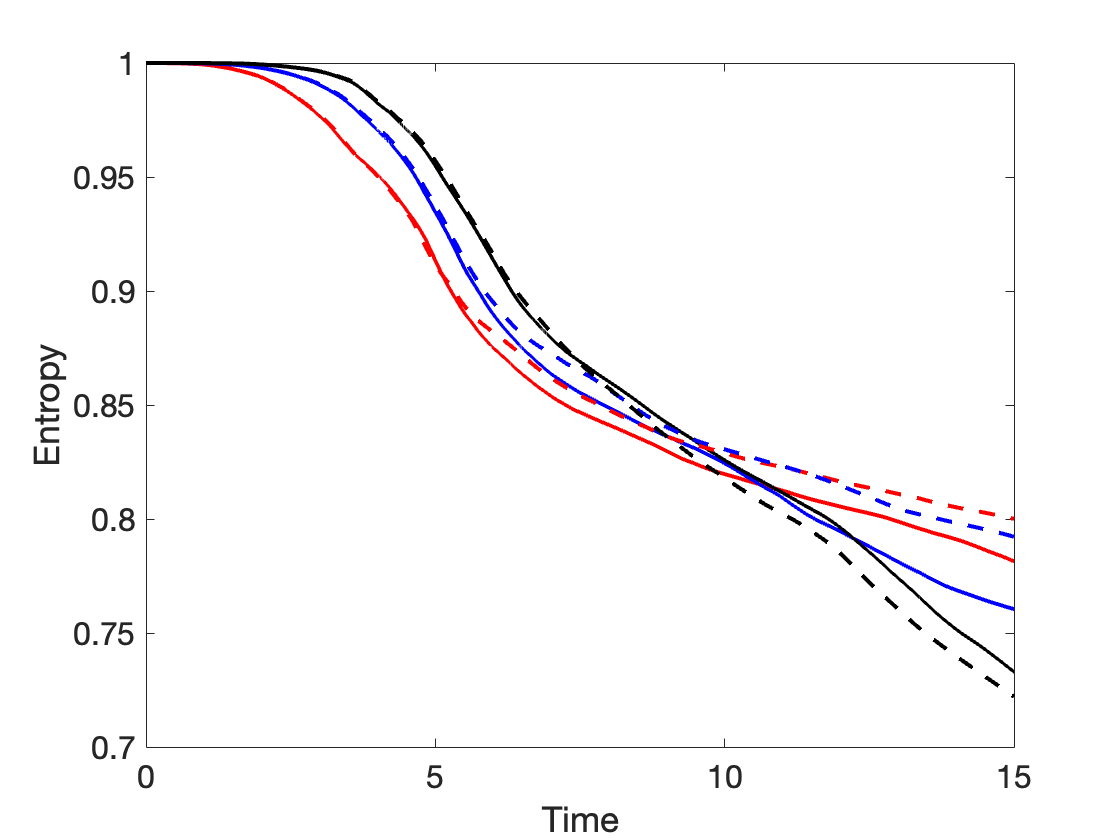}
	\caption{Normalized Entropy. Red: $64^2$; Blue: $128^2$; Black: $256^2$. Solid lines: Discrete (RB-TAI); Dashed lines: Semi Discrete.}
	\label{fig:shearstatsnew}
\end{figure}

{It is notable that without limiters or positivity preservation, only very robust methods, such as the entropy stable discontinuous Galerkin method  described in ~\cite{chan2022entropy} are able to run successfully beyond t = 5 units for this problem.}
Figure~\ref{fig:shearstatsnew} shows that the entropy production is much more pronounced compared to the nearly incompressible case, giving more credence to the comprehensiveness of the Active Flux method.

\begin{figure}
	\centering
		\includegraphics[width=0.32\textwidth]{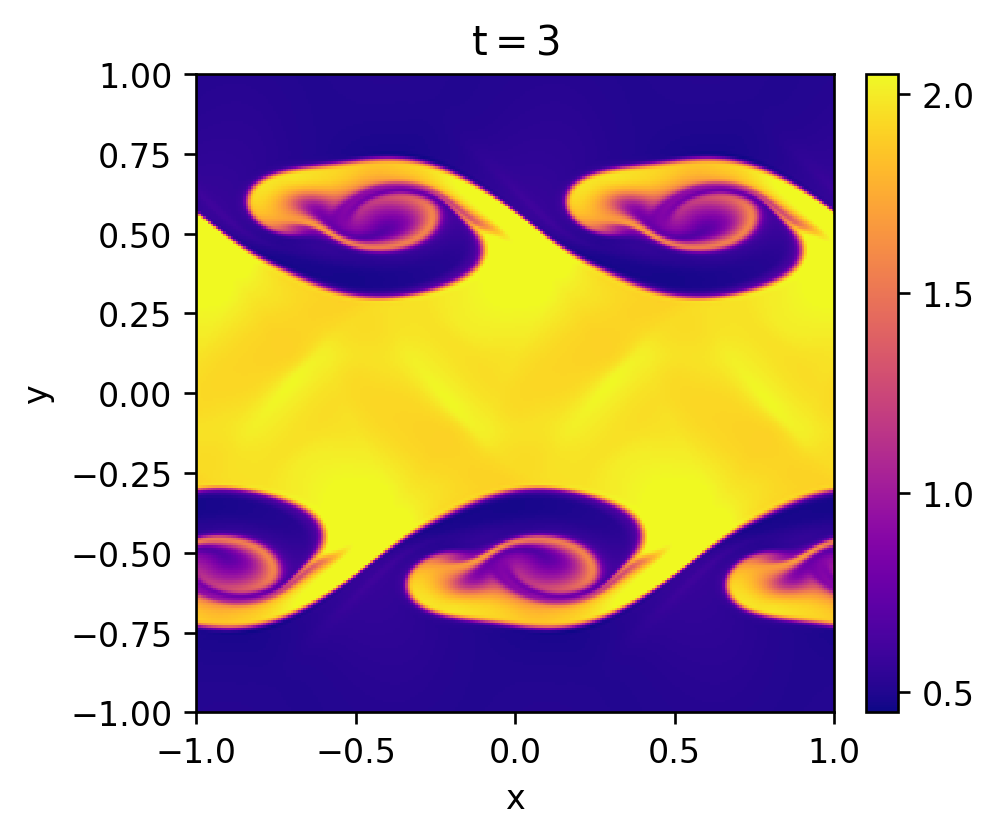}
				\includegraphics[width=0.32\textwidth]{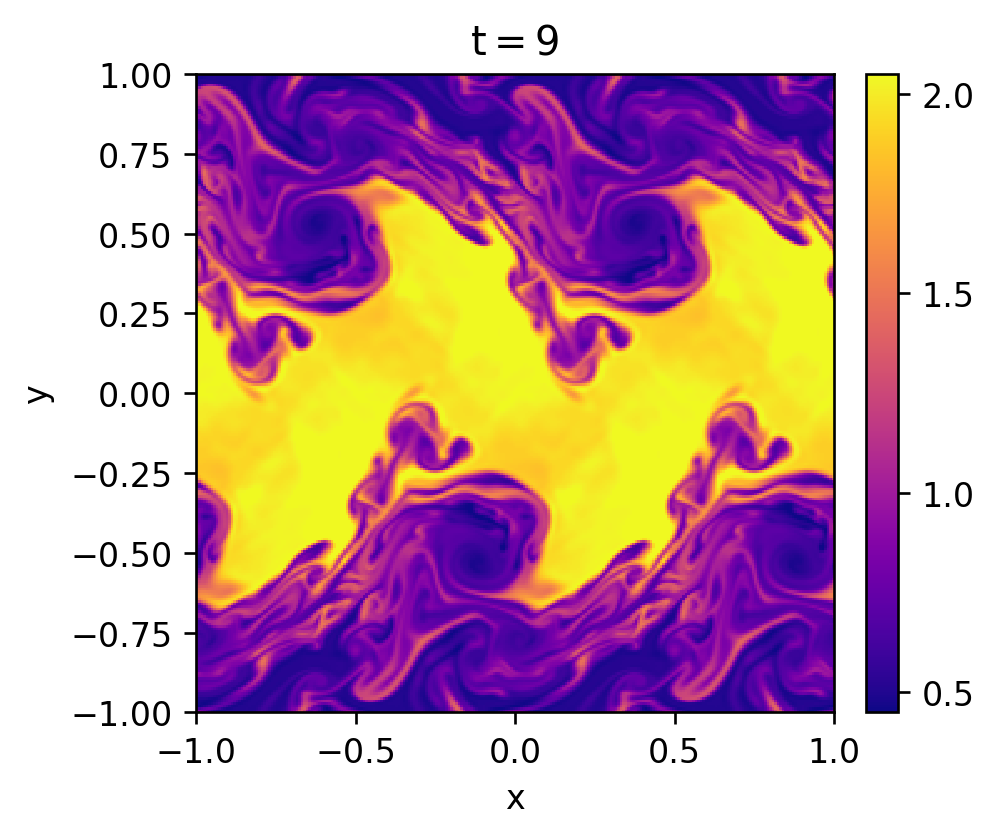}
					\includegraphics[width=0.32\textwidth]{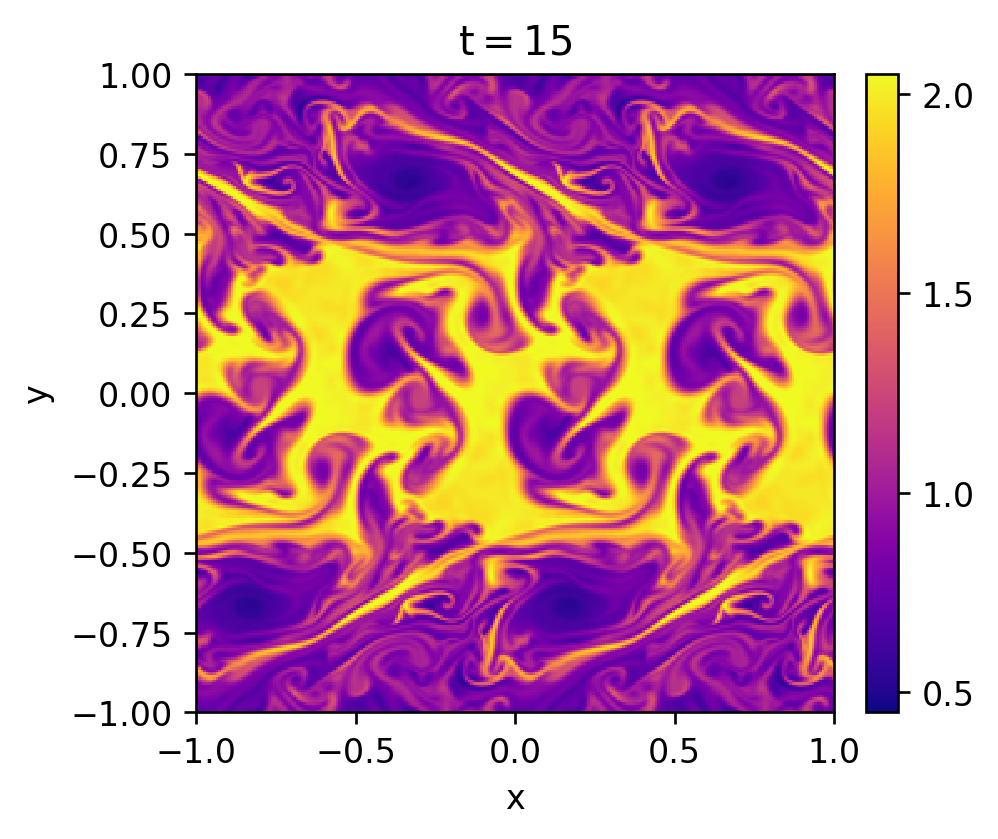}
		\caption{Under-resolved KH problem on $256^2$ mesh, Discrete RB-TAI}
	\label{fig:kh256}
	\end{figure}



\begin{figure}
	\centering
	\includegraphics[width=0.49\textwidth]{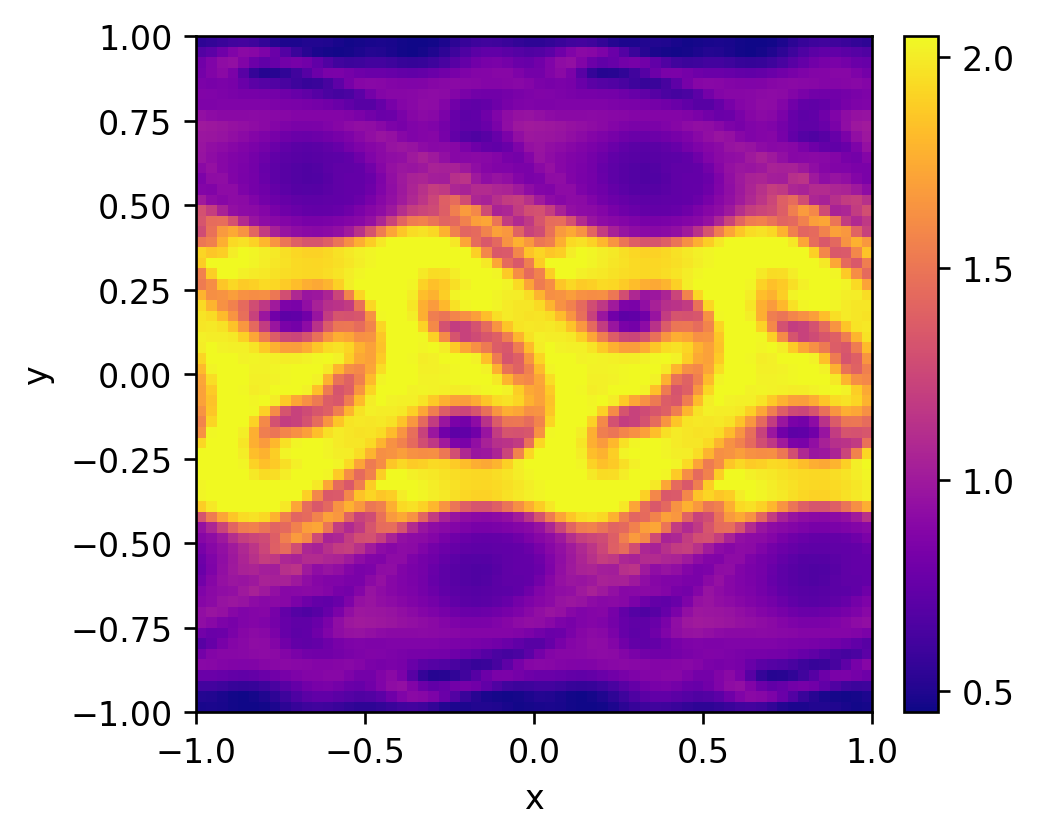}
	\includegraphics[width=0.49\textwidth]{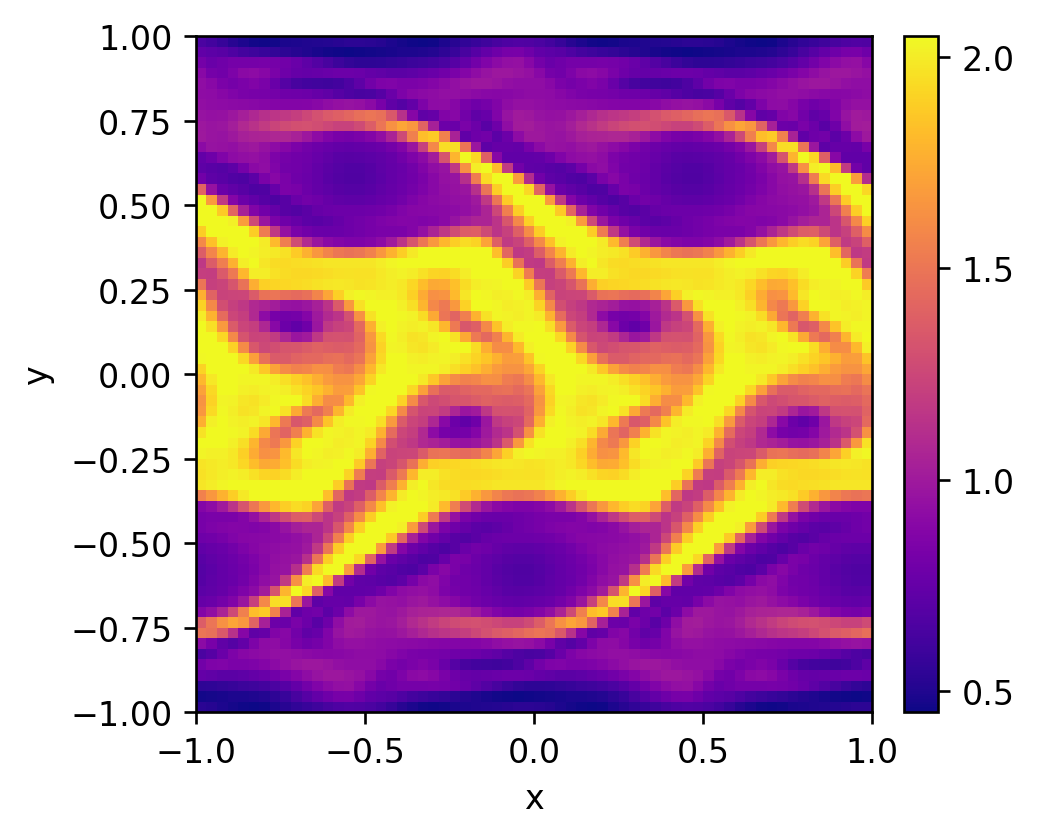}
		\caption{Density contours for under-resolved KH problem on $64^2$ mesh. Left: Discrete RB, Right: Discrete RB-TAI}
	\label{fig:shearnew64}
\end{figure}


\section{Conclusions and Perspectives}
This work introduced a transported acoustic-increment version of the fully discrete multidimensional Active Flux method for the compressible Euler equations. In the additive Roe–Barsukow~\citep{barsukow2025} point update, the acoustic evolution computed at a grid point is added back at that same Eulerian point, although in the presence of background advection that acoustic information is naturally associated with a material label that has moved. Reconstructing the acoustic increment and evaluating it at the convective foot restores this missing transport mechanism without abandoning the compact Active Flux degrees of freedom or the conservative cell-average update. 

The transported increment exploits the fact that the locally linearized solve still contains valuable information, and thus adds only the reconstruction and evaluation of the acoustic increment at a convective foot, and uses the same conservative space-time Simpson update for cell averages. The cell-center increment is essential in this construction: without it, the {$Q_2$ } increment field is not determined in the cell interior, and the transported correction would be a boundary extension rather than a genuine Active Flux reconstruction.

{Intuitively,} the additive update has the form $e^{\tau A}+e^{\tau C}-I$ and misses the symmetric second-order advection–acoustic cross term. The transported-increment update instead gives $e^{\tau A}(e^{\tau C}-I)+e^{\tau A}=e^{\tau A}e^{\tau C}$. Thus, for constant frozen coefficients, where the advective and acoustic generators commute, the transported update reproduces the exact unsplit frozen evolution. {In more general linear settings,  the transported form replaces the additive symmetric defect by a commutator defect, which is the expected residual when the advective and acoustic generators do not commute.}

The numerical results support this operator-level interpretation across a broad range of tests. 
\begin{enumerate}
\item In the mixed  wave-packet { linearized Euler} diagnostic, the additive fully discrete point update exhibits the expected second-order one-step defect, while the transported-increment update recovers third-order convergence and approaches the accuracy of a direct convective-foot acoustic update. 

\item In the isentropic vortex problem, the transported method {follows third-order trends in a practical resolution limit  for the non-linear Euler equations} while reducing the error constant relative to both the additive fully discrete method and the semi-discrete Active Flux reference. In the tested vortex configuration, it also remains stable up to total CFL number one, giving an empirical high-CFL behavior consistent with the multiplicative structure of the frozen point update.

\item The acoustic-pulse experiment shows that the transport modification does not degrade the underlying multidimensional acoustic propagation. The radial pulse remains symmetric on Cartesian grids, and the measured pressure error and radial-asymmetry error decrease at approximately third order. This is important because the proposed method does not replace the exact locally linearized acoustic operator; it only transports its increment. The symmetry behavior therefore indicates that the transported correction remains compatible with the multidimensional acoustic structure on which Active Flux is built.

\item The low-Mach shear-layer experiment probes a different regime. At $M=10^{-2}$, the transported method preserves the large-scale coherent shear structures on coarse meshes and avoids the secondary vortex structures observed in DG/CG comparisons. {The integral of vorticity remains at machine precision throughout the computation, while the entropy change is minuscule,  thus consistent with a nearly incompressible, smooth shear evolution.} 

\item The compressible Kelvin–Helmholtz experiment provides a complementary robustness test. In this highly unstable under-resolved configuration, the Active Flux variants evolve to the target time  on meshes as coarse as $64^2$. The discrete methods admit larger time steps than the semi-discrete reference, and the transported variant gives visibly sharper structures in the under-resolved regime. The entropy history shows stronger dissipation than in the low-Mach shear case, as expected for a compressible flow developing fine-scale structure. Thus the entropy behavior is consistent across the two shear tests: nearly conservative in the smooth low-Mach case, and more dissipative in the compressible unstable case.
\end{enumerate}

Several related alternatives were considered. A direct convective-foot acoustic disk integration was derived. This construction is useful as an analytical reference and recovers the standard Barsukow half-plane and quadrant formulas in limiting cases. {However, in our prototype, the direct off-center update was more than an order of magnitude more expensive than the original point update, while offering only slight accuracy benefits in linear problems and slightly lower accuracy in non-linear problems.} Midpoint refreezing, direction-dependent acoustic refreezing, and Strang-type compositions are also natural possibilities, but they move the method away from the single-stage, compact, exact-acoustic Active Flux structure.  The transported-increment method {is therefore best viewed as a minimal correction that captures the missing moving-frame placement of acoustic information while leaving the core Roe–Barsukow infrastructure intact.}

Several questions remain open. {Provable third order accuracy for multi-dimensional non-linear problems is a formidable challenge.} The present work  does not address discontinuous-flow stabilization. All Active Flux computations reported here were performed without limiters or positivity corrections. The compressible Kelvin-Helmholtz calculation is a useful robustness test, but it is not a substitute for shock-capturing evidence.  Extending the transported-increment method to flows with shocks, contacts, and near-vacuum states will require nonlinear limiting{, and positivity preservation strategies}.

\section*{Acknowledgments}
This work was supported by Los Alamos National Laboratory under the grant `Advanced Algorithms for Multiphysics Analysis'. The author is deeply indebted to the late Professor Phil Roe of the University of Michigan for several discussions on related topics, and ultimately inspiring the ideas that led to this work. This paper is dedicated to his memory. 

This work has also benefitted from discussions with Prof. Praveen Chandrasekhar (TIFR, Bangalore), Prof. Christian Klingenberg (Wuerzburg), Ms. Moon Hazarika (U. Michigan), Josh Dolence \& Chad Meyer (Los Alamos). Bjorn Kierulf, Eric Johnsen, and Praveen Chandrasekhar provided DG/CG solution comparisons. 

OpenAI's Codex was used for coding and numerical postprocessing, and all outputs were verified multiple times by the author.  Scientific content, reasoning, and conclusions are the author's own.

\section*{Appendix}

\appendix

\section{Semi-discrete Active Flux}
\label{sec:semidiscrete}

The semi-discrete Active Flux method~\citep{abgrall2025semieuler} is the method-of-lines analogue.  One keeps the same degrees of freedom and reconstruction, but evolves the degrees of freedom by an ODE system instead of using a fully discrete evolution operator.

Let
$
	Z_h(t)=\{\avg{\U}_{ij}(t),\W^V_{ij}(t),\W^{E_v}_{ij}(t),\W^{E_h}_{ij}(t)\}_{i,j}
$
collect all degrees of freedom.  The reconstruction operator is denoted by
$
	\W_h(\cdot,t)=\calR_h Z_h(t).
$

A semi-discrete point update evaluates the right-hand side using derivatives of the reconstruction, or using an equivalent finite-difference formula based on the Active Flux degrees of freedom.  In compact notation,
\begin{equation}
	\frac{\dd}{\dd t}\W_P(t)=\calL^{P}_{\rm SD}(Z_h(t)).
\end{equation}
For example, using the reconstructed derivative directly gives
\begin{equation}
	\calL^{P}_{\rm SD}(Z_h)=
	-\begin{pmatrix}
		u\rho_x+v\rho_y+\rho(u_x+v_y)\\[2pt]
		u u_x+v u_y+p_x/\rho\\[2pt]
		u v_x+v v_y+p_y/\rho\\[2pt]
		u p_x+v p_y+\gamma p(u_x+v_y)
	\end{pmatrix}_{P,\;\W_h=\calR_h Z_h}.
\end{equation}
At vertices and edge midpoints the derivative of a globally continuous $Q_2$ reconstruction is not generally single-valued; a semi-discrete implementation must choose a consistent derivative average or a finite-difference residual.  This is one reason why the fully discrete evolution-operator form is attractive: the acoustic evolution operator naturally uses the adjacent wedges/half-planes around a point. {The reference implementation in this work faithfully follows the finite difference Jacobian splitting introduced in ~\cite{abgrall2025semieuler}.}

The conservative cell average satisfies
\begin{equation}
	\frac{\dd}{\dd t}\avg{\U}_{ij}
	=
	-\frac{1}{h}\left(\widehat{\F}_{i+1,j}-\widehat{\F}_{i,j}\right)
	-\frac{1}{h}\left(\widehat{\G}_{i,j+1}-\widehat{\G}_{i,j}\right),
\end{equation}
where the instantaneous edge fluxes are approximated by Simpson quadrature along each edge.  The semi-discrete ODE therefore has the form
\begin{equation}
		\frac{\dd}{\dd t}Z_h(t)=\calL_{\rm SD}(Z_h(t)).
\end{equation}
{The SSPRK3 Runge--Kutta method is used for time integration.}

\section{Explicit Form of Acoustic Update at Vertical edge node}
\label{sec:acoustic_update}
In the active flux solver implementation, the acoustic point value updates are written in the form of monomial coefficients as given by \cite{barsukow2025}. However, it is possible to write these updates in a more explicit form that also provides more insight into the underlying algorithm, and makes it more amenable for analysis. The update of the point node at the center of the vertical cell face is given as an example in this appendix.

\subsection{Right half contribution}
Let \(U^R,V^R,P^R\) be the right-cell nodal arrays for \(u,v,p\).  Set
\[
\nu=\frac{\alpha}{h}\triangleq \frac{c_*\Delta t}{h}.
\]
The right half-plane acoustic contribution is the direct nodal stencil
\[
	\begin{aligned}
		U_R\triangleq{}&
		\langle A^R_{uu},U^R\rangle
		+\langle A^R_{uv},V^R\rangle
		+\frac1{Z_*}\langle A^R_{up},P^R\rangle,\\
		V_R\triangleq{}&
		\langle A^R_{uv},U^R\rangle
		+\langle A^R_{vv},V^R\rangle
		+\frac1{Z_*}\langle A^R_{vp},P^R\rangle,\\
		P_R\triangleq{}&
		Z_*\langle A^R_{up},U^R\rangle
		+Z_*\langle A^R_{vp},V^R\rangle
		+\langle A^R_{pp},P^R\rangle .
\end{aligned}
\]
The right-cell nodal stencils are
\begingroup
\scriptsize
\setlength{\arraycolsep}{3pt}
\[
A^R_{uu}=
\begin{bmatrix}
	\frac{\nu^{3}(2\nu-3)}{6} &
	-\frac{4\nu^{4}-6\nu^{3}-6\nu^{2}+9\nu-3}{6} &
	\frac{\nu^{3}(2\nu-3)}{6} \\
	-\frac{2\nu^{3}(\nu-1)}{3} &
	\frac{2\nu(\nu-1)(2\nu^{2}-3)}{3} &
	-\frac{2\nu^{3}(\nu-1)}{3} \\
	\frac{\nu^{3}(2\nu-1)}{6} &
	-\frac{\nu(2\nu-1)(2\nu^{2}-3)}{6} &
	\frac{\nu^{3}(2\nu-1)}{6}
\end{bmatrix},
\]
\[
A^R_{uv}=
\begin{bmatrix}
	-\frac{\nu(4\nu^{2}-9\nu+6)}{12} & 0 &
	\frac{\nu(4\nu^{2}-9\nu+6)}{12} \\
	\frac{\nu^{2}(2\nu-3)}{3} & 0 &
	-\frac{\nu^{2}(2\nu-3)}{3} \\
	-\frac{\nu^{2}(4\nu-3)}{12} & 0 &
	\frac{\nu^{2}(4\nu-3)}{12}
\end{bmatrix},
\]
\[
A^R_{up}=
\begin{bmatrix}
	-\frac{\nu^{2}(\nu-1)^{2}}{2} &
	\frac{(\nu-1)(2\nu^{3}-2\nu^{2}-2\nu+1)}{2} &
	-\frac{\nu^{2}(\nu-1)^{2}}{2} \\
	\frac{\nu^{3}(3\nu-4)}{3} &
	-\frac{2\nu(3\nu^{3}-4\nu^{2}-3\nu+3)}{3} &
	\frac{\nu^{3}(3\nu-4)}{3} \\
	-\frac{\nu^{3}(3\nu-2)}{6} &
	\frac{\nu(6\nu^{3}-4\nu^{2}-6\nu+3)}{6} &
	-\frac{\nu^{3}(3\nu-2)}{6}
\end{bmatrix}.
\]
\[
A^R_{vv}=
\begin{bmatrix}
	\frac{\nu^{2}(\nu^{2}-3\nu+3)}{3} &
	-\frac{4\nu^{4}-12\nu^{3}+12\nu^{2}-3}{6} &
	\frac{\nu^{2}(\nu^{2}-3\nu+3)}{3}\\
	-\frac{2\nu^{3}(\nu-2)}{3} &
	\frac{4\nu^{3}(\nu-2)}{3} &
	-\frac{2\nu^{3}(\nu-2)}{3}\\
	\frac{\nu^{3}(\nu-1)}{3} &
	-\frac{2\nu^{3}(\nu-1)}{3} &
	\frac{\nu^{3}(\nu-1)}{3}
\end{bmatrix},
\]
\[
A^R_{vp}=
\begin{bmatrix}
	\frac{\nu(4\nu^{2}-9\nu+6)}{12} & 0 &
	-\frac{\nu(4\nu^{2}-9\nu+6)}{12}\\
	-\frac{\nu^{2}(2\nu-3)}{3} & 0 &
	\frac{\nu^{2}(2\nu-3)}{3}\\
	\frac{\nu^{2}(4\nu-3)}{12} & 0 &
	-\frac{\nu^{2}(4\nu-3)}{12}
\end{bmatrix},
\]
\[
A^R_{pp}=
\begin{bmatrix}
	\frac{\nu^{2}(4\nu^{2}-9\nu+6)}{6} &
	-\frac{8\nu^{4}-18\nu^{3}+6\nu^{2}+9\nu-3}{6} &
	\frac{\nu^{2}(4\nu^{2}-9\nu+6)}{6}\\
	-\frac{2\nu^{3}(2\nu-3)}{3} &
	\frac{2\nu(4\nu^{3}-6\nu^{2}-3\nu+3)}{3} &
	-\frac{2\nu^{3}(2\nu-3)}{3}\\
	\frac{\nu^{3}(4\nu-3)}{6} &
	-\frac{\nu(8\nu^{3}-6\nu^{2}-6\nu+3)}{6} &
	\frac{\nu^{3}(4\nu-3)}{6}
\end{bmatrix}.
\]
\endgroup

\subsection{Left half contribution}

Let \(U^L,V^L,P^L\) be the left-cell nodal arrays.  The left half-plane
contribution is
\[
	\begin{aligned}
		U_L={}&
		\langle A^L_{uu},U^L\rangle
		+\langle A^L_{uv},V^L\rangle
		+\frac1{Z_*}\langle A^L_{up},P^L\rangle,\\
		V_L={}&
		\langle A^L_{uv},U^L\rangle
		+\langle A^L_{vv},V^L\rangle
		+\frac1{Z_*}\langle A^L_{vp},P^L\rangle,\\
		P_L={}&
		Z_*\langle A^L_{up},U^L\rangle
		+Z_*\langle A^L_{vp},V^L\rangle
		+\langle A^L_{pp},P^L\rangle .
\end{aligned}
\]
The left-cell nodal stencils are
\begingroup
\scriptsize
\setlength{\arraycolsep}{3pt}
\[
A^L_{uu}=
\begin{bmatrix}
	\frac{\nu^{3}(2\nu-1)}{6} &
	-\frac{\nu(2\nu-1)(2\nu^{2}-3)}{6} &
	\frac{\nu^{3}(2\nu-1)}{6} \\
	-\frac{2\nu^{3}(\nu-1)}{3} &
	\frac{2\nu(\nu-1)(2\nu^{2}-3)}{3} &
	-\frac{2\nu^{3}(\nu-1)}{3} \\
	\frac{\nu^{3}(2\nu-3)}{6} &
	-\frac{4\nu^{4}-6\nu^{3}-6\nu^{2}+9\nu-3}{6} &
	\frac{\nu^{3}(2\nu-3)}{6}
\end{bmatrix},
\]
\[
A^L_{uv}=
\begin{bmatrix}
	\frac{\nu^{2}(4\nu-3)}{12} & 0 &
	-\frac{\nu^{2}(4\nu-3)}{12} \\
	-\frac{\nu^{2}(2\nu-3)}{3} & 0 &
	\frac{\nu^{2}(2\nu-3)}{3} \\
	\frac{\nu(4\nu^{2}-9\nu+6)}{12} & 0 &
	-\frac{\nu(4\nu^{2}-9\nu+6)}{12}
\end{bmatrix},
\]
\[
A^L_{up}=
\begin{bmatrix}
	\frac{\nu^{3}(3\nu-2)}{6} &
	-\frac{\nu(6\nu^{3}-4\nu^{2}-6\nu+3)}{6} &
	\frac{\nu^{3}(3\nu-2)}{6} \\
	-\frac{\nu^{3}(3\nu-4)}{3} &
	\frac{2\nu(3\nu^{3}-4\nu^{2}-3\nu+3)}{3} &
	-\frac{\nu^{3}(3\nu-4)}{3} \\
	\frac{\nu^{2}(\nu-1)^{2}}{2} &
	-\frac{(\nu-1)(2\nu^{3}-2\nu^{2}-2\nu+1)}{2} &
	\frac{\nu^{2}(\nu-1)^{2}}{2}
\end{bmatrix}.
\]
\[
A^L_{vv}=
\begin{bmatrix}
	\frac{\nu^{3}(\nu-1)}{3} &
	-\frac{2\nu^{3}(\nu-1)}{3} &
	\frac{\nu^{3}(\nu-1)}{3}\\
	-\frac{2\nu^{3}(\nu-2)}{3} &
	\frac{4\nu^{3}(\nu-2)}{3} &
	-\frac{2\nu^{3}(\nu-2)}{3}\\
	\frac{\nu^{2}(\nu^{2}-3\nu+3)}{3} &
	-\frac{4\nu^{4}-12\nu^{3}+12\nu^{2}-3}{6} &
	\frac{\nu^{2}(\nu^{2}-3\nu+3)}{3}
\end{bmatrix},
\]
\[
A^L_{vp}=
\begin{bmatrix}
	\frac{\nu^{2}(4\nu-3)}{12} & 0 &
	-\frac{\nu^{2}(4\nu-3)}{12}\\
	-\frac{\nu^{2}(2\nu-3)}{3} & 0 &
	\frac{\nu^{2}(2\nu-3)}{3}\\
	\frac{\nu(4\nu^{2}-9\nu+6)}{12} & 0 &
	-\frac{\nu(4\nu^{2}-9\nu+6)}{12}
\end{bmatrix},
\]
\[
A^L_{pp}=
\begin{bmatrix}
	\frac{\nu^{3}(4\nu-3)}{6} &
	-\frac{\nu(8\nu^{3}-6\nu^{2}-6\nu+3)}{6} &
	\frac{\nu^{3}(4\nu-3)}{6}\\
	-\frac{2\nu^{3}(2\nu-3)}{3} &
	\frac{2\nu(4\nu^{3}-6\nu^{2}-3\nu+3)}{3} &
	-\frac{2\nu^{3}(2\nu-3)}{3}\\
	\frac{\nu^{2}(4\nu^{2}-9\nu+6)}{6} &
	-\frac{8\nu^{4}-18\nu^{3}+6\nu^{2}+9\nu-3}{6} &
	\frac{\nu^{2}(4\nu^{2}-9\nu+6)}{6}
\end{bmatrix}.
\]
\endgroup

The acoustic point value at the vertical face node is
\begin{equation}
		u_{\rm ac}^{n+1}(P)=U_R+U_L,\qquad
		v_{\rm ac}^{n+1}(P)=V_R+V_L,\qquad
		p_{\rm ac}^{n+1}(P)=P_R+P_L .
\end{equation}
The density follows from the frozen acoustic invariant
\(\rho-p/c_*^2=\text{constant}\) at the point:
\begin{equation}
		\rho_{\rm ac}^{n+1}(P)
		=
		\rho_*+\frac{p_{\rm ac}^{n+1}(P)-p_*}{c_*^2}.
\end{equation}

{
\section{Foot-centered implementation for total CFL below 0.5}
\label{sec:footcentered}

In contrast to the transported increment approach which interpolates nodal increments at the material foot F, the foot-centered implementation computes the increment directly at F. The baseline implementation of this algorithm is given below for a total CFL $\leq 0.5$:

For a boundary point \(P\) and a substep \(\tau\in\{\Dt/2,\Dt\}\), compute the
location of the material foot (e.g. using two Picard steps ) and freeze the acoustic coefficients there:
\begin{equation}
	\W_*=\W_h^n(F),\quad
	c_*=\sqrt{\gamma p_*/\rho_*},\quad Z_*=\rho_*c_*,\quad R=c_*\tau .
	\label{eq:fc-foot-freeze}
\end{equation}
With \(\delta=F-P\), the compact footprint condition on an equal Cartesian
mesh is
\begin{equation}
	|\delta_x|+R<\frac h2,\qquad |\delta_y|+R<\frac h2.
	\label{eq:fc-half-footprint}
\end{equation}
It follows, for example, from \(\Dt s_{\max}/h<1/2\).  Thus, as shown in Figure~\ref{fig:fc-moment-primitives}, a vertex disk
crosses at most its two grid lines, a vertical- or horizontal-edge disk only
its normal grid line, and a center disk no grid line. 

Center and scale the disk by \((x,y)=F+R(\xi,\eta)\), let
\(B_1=\{\xi^2+\eta^2<1\}\), and use the projected-sphere weight
\(\Psi_{\rm sph}(\xi,\eta)=2(1-\xi^2-\eta^2)^{-1/2}\).  If
\(\widehat\Omega_r(R)\subset B_1\) is the part belonging to one adjacent
piecewise-$Q_2$ cell, define the area, radius-derivative, and history moments by
\begin{equation}
			\mathcal M^r_{ij}(R)=\int_{\widehat\Omega_r(R)}
			\xi^i\eta^j\Psi_{\rm sph}\,\dd\xi\dd\eta,\ \
			\mathcal M^{r,D}_{ij,m}(R)=
			\frac{\dd}{\dd R}\!\left(R^m\mathcal M^r_{ij}(R)\right),\ \
			\mathcal M^{r,H}_{ij,q}(R)=
			\int_{R_b^r}^{R}s^{q-1}\mathcal M^r_{ij}(s)\,\dd s .
	\label{eq:fc-moment-families}
\end{equation}
Here \(R_b^r\) is the fragment birth radius; only indices generated by
\(0\le a,b\le2\) are required.

Only three geometric primitives occur.  For \(\lambda=d_x/R\) and
\(\mu=d_y/R\), their area moments are
\begin{equation}
	\begin{aligned}
		\mathcal M^F_{ij}&=\int_{B_1}\xi^i\eta^j\Psi_{\rm sph}\,
		\dd\xi\dd\eta,\\[-1mm]
		\mathcal M^x_{ij}(\lambda)&=
		\int_{B_1\cap\{\xi\ge\lambda\}}
		\xi^i\eta^j\Psi_{\rm sph}\,\dd\xi\dd\eta,\qquad
		\mathcal M^y_{ij}(\mu)=\mathcal M^x_{ji}(\mu),\\[-1mm]
		\mathcal M^{xy}_{ij}(\lambda,\mu)&=
		\int_{B_1\cap\{\xi\ge\lambda,\,\eta\ge\mu\}}
		\xi^i\eta^j\Psi_{\rm sph}\,\dd\xi\dd\eta .
	\end{aligned}
	\label{eq:fc-cap-moments}
\end{equation}

\begin{figure}[!ht]
	\centering
	\vspace{-1mm}
	\begin{tikzpicture}[
		x=1.40cm,y=1.40cm,
		disk/.style={draw=blue!60!black,line width=0.7pt},
		grid/.style={draw=black!45,dashed,line width=0.45pt},
		axis/.style={draw=black!70,-{Latex[length=1.3mm]},line width=0.45pt},
		transport/.style={draw=black!55,-{Latex[length=1.4mm]},line width=0.45pt},
		target/.style={fill=black},
		foot/.style={fill=red!70!black},
		every node/.style={font=\scriptsize}]
		\begin{scope}[shift={(0,0)}]
			\coordinate (PA) at (0,0.60);
			\coordinate (FA) at (-0.48,0.48);
			\fill[blue!13] (FA) circle (0.40);
			\draw[grid] (-1,-1) grid[step=1] (1,1);
			\draw[disk] (FA) circle (0.40);
			\draw[transport] (PA)--(FA);
			\draw[axis] (FA)--++(0.28,0) node[right] {$\xi$};
			\draw[axis] (FA)--++(0,0.28) node[above] {$\eta$};
			\fill[target] (PA) circle (1.15pt) node[above right] {$P$};
			\fill[foot] (FA) circle (1.25pt) node[below left] {$F$};
			\node at (0,-1.25) {full disk $\mathcal M^F$};
		\end{scope}
		\begin{scope}[shift={(2.65,0)}]
			\coordinate (PB) at (0,0.60);
			\coordinate (FB) at (-0.18,0.48);
			\begin{scope}
				\clip (FB) circle (0.40);
				\fill[orange!32] (0,0) rectangle (1,1);
			\end{scope}
			\draw[grid] (-1,-1) grid[step=1] (1,1);
			\draw[disk] (FB) circle (0.40);
			\draw[transport] (PB)--(FB);
			\draw[axis] (FB)--++(0.28,0) node[right] {$\xi$};
			\draw[axis] (FB)--++(0,0.28) node[above] {$\eta$};
			\fill[target] (PB) circle (1.15pt) node[above right] {$P$};
			\fill[foot] (FB) circle (1.25pt) node[below left] {$F$};
			\node[below right] at (0,0.09) {$\lambda$};
			\node at (0,-1.25) {$x$-cap $\mathcal M^x(\lambda)$};
		\end{scope}
		\begin{scope}[shift={(5.30,0)}]
			\coordinate (PC) at (0,0);
			\coordinate (FC) at (-0.18,-0.14);
			\begin{scope}
				\clip (FC) circle (0.40);
				\fill[green!28] (0,0) rectangle (1,1);
			\end{scope}
			\draw[grid] (-1,-1) grid[step=1] (1,1);
			\draw[disk] (FC) circle (0.40);
			\draw[transport] (PC)--(FC);
			\draw[axis] (FC)--++(0.28,0) node[right] {$\xi$};
			\draw[axis] (FC)--++(0,0.28) node[above] {$\eta$};
			\fill[target] (PC) circle (1.15pt) node[above right] {$P$};
			\fill[foot] (FC) circle (1.25pt) node[below left] {$F$};
			\node[below right] at (0,-0.35) {$\lambda$};
			\node[above left] at (-0.34,0) {$\mu$};
			\node at (0,-1.25) {double cap $\mathcal M^{xy}(\lambda,\mu)$};
		\end{scope}
	\end{tikzpicture}
	\caption{Moment primitives in their local grid geometry.  Dashed lines are grid
		lines, with the same cell size in every panel; \(P\) is the target, \(F\) the
		disk center, and the short arrows at \(F\) define \((\xi,\eta)\).  Shading
		identifies the full disk, cap, and double cap.}
	\label{fig:fc-moment-primitives}
	\vspace{-1mm}
\end{figure}

Opposite orientations multiply the cap moments by \(s_x^i\), \(s_y^j\), or
\(s_x^is_y^j\).  The double cap is empty if
\(d_x^2+d_y^2\ge R^2\).  Closed primitives supply every member of the
\(\mathcal M\) family.

Let \(\mathcal M_F,\mathcal M_x,\mathcal M_y,\mathcal M_{xy}\) denote the
complete area, derivative, and history moment arrays.  Inclusion--exclusion
assigns the possible vertex-cell fragments by
\begin{equation}
	\mathcal M_{00}=\mathcal M_F-\mathcal M_x-\mathcal M_y+\mathcal M_{xy},\quad
	\mathcal M_{10}=\mathcal M_x-\mathcal M_{xy},\quad
	\mathcal M_{01}=\mathcal M_y-\mathcal M_{xy},\quad
	\mathcal M_{11}=\mathcal M_{xy}.
	\label{eq:fc-role-rows}
\end{equation}
For an edge, retain only the corresponding one-cap split; for a center, use
\(\mathcal M_F\).  Express each contributing $Q_2$ polynomial in monomials about
\(F\), contract its \((p/Z_*,u,v)\) coefficients with the assigned moment
array, and sum the fragments.  This gives the exact frozen-acoustic
spherical-mean update without angular or radial quadrature.  Set
\(p^A=Z_*\pi^A\) and use
\(\rho^A=\rho_h(F)+(p^A-p_h(F))/c_*^2\) for the linear symbol, or the
entropy-consistent nonlinear closure
\(\rho^A=\rho_h(F)[p^A/p_h(F)]^{1/\gamma}\).  Applying this point map at
\(\tau=\Dt/2,\Dt\) completes the standard space--time Simpson average update
\eqref{eq:space-time-edge}--\eqref{eq:average-row}.

\section{von Neumann analysis quantities}

\label{sec:acoustic-symbol}

Cells are \(C_{ij}=[x_i,x_{i+1}]\times[y_j,y_{j+1}]\), with equal spacing
\(h\).  After linearizing the conservative-to-primitive conversion about
\(\W_0\), it is convenient to express every family in primitive perturbations.
This is a blockwise similarity transformation and therefore leaves all
eigenvalues unchanged.  Use the Fourier ansatz
\begin{equation}
	\begin{bmatrix}
		\overline{\w}_{ij}\\ \w^V_{ij}\\ \w^{E_v}_{ij}\\ \w^{E_h}_{ij}
	\end{bmatrix}
	=z_x^iz_y^j\q,
	\qquad z_x=\e^{\ii\xi},\quad z_y=\e^{\ii\eta},\quad
	\q\in\mathbb C^{16}.
	\label{eq:fourier-state}
\end{equation}

Let a scalar family amplitude be ordered as
\(\widehat\phi=(\overline\phi,\phi^V,\phi^{E_v},\phi^{E_h})^T\).
Linearization of the conservative Simpson inversion gives exactly
\begin{equation}
	\widehat\phi^C=\mathcal C(z_x,z_y)\widehat\phi,
	\qquad
	\mathcal C=\left[
	\frac94,
	-\frac{(1+z_x)(1+z_y)}{16},
	-\frac{1+z_x}{4},
	-\frac{1+z_y}{4}
	\right].
	\label{eq:center-row}
\end{equation}
This row is the full-state counterpart of the center-node Fourier factor used
in Section~6 of the manuscript.

Define the following nine \(1\times4\) rows, which map scalar family
amplitudes to the $Q_2$ nodes of cell \(C_{ij}\):
\begin{equation}
	\mathsf N(z_x,z_y)=
	\begin{bmatrix}
		e_V & e_{E_h} & z_xe_V\\
		e_{E_v} & \mathcal C & z_xe_{E_v}\\
		z_ye_V & z_ye_{E_h} & z_xz_ye_V
	\end{bmatrix},
	\label{eq:nodal-row-matrix}
\end{equation}
where \(e_V=(0,1,0,0)\), \(e_{E_v}=(0,0,1,0)\), and
\(e_{E_h}=(0,0,0,1)\).  The entry \(N_{mn}\) is itself a row.  For cell
\(C_{i+r,j+s}\), every row in \ref{eq:nodal-row-matrix} is multiplied by
\(z_x^rz_y^s\).

For later use, the scalar reconstruction row at a general nondimensional point
\((x,y)\), measured in cells from \((x_i,y_j)\), is
\begin{equation}
	\R(x,y)=z_x^rz_y^s
	\sum_{m=0}^2\sum_{n=0}^2
	\ell_m(2[x-r-\tfrac12])\ell_n(2[y-s-\tfrac12])N_{mn},
	\quad r=\lfloor x\rfloor,\quad s=\lfloor y\rfloor,
	\label{eq:reconstruction-row}
\end{equation}
where \(\ell_0(a)=a(a-1)/2\), \(\ell_1(a)=1-a^2\), and
\(\ell_2(a)=a(a+1)/2\).  Equation~\eqref{eq:reconstruction-row} includes the
phase change when a foot lies in a neighboring cell.

The exact acoustic update can be converted directly into a Fourier matrix
without numerical quadrature.  This section states that conversion in a form
that covers vertices and both edge-midpoint families.

Introduce \(\pi=p/Z_0\), so that frozen acoustics is
\begin{equation}
	\bm u_t+c_0\nabla\pi=0,\qquad
	\pi_t+c_0\nabla\!\cdot\bm u=0.
\end{equation}
For a wedge \(\Theta=[\theta_0,\theta_1]\), define
\begin{equation}
	\omega_\Theta=\frac{\theta_1-\theta_0}{2\pi},\qquad
	\mu_{mn}^\Theta=\int_\Theta\cos^m\theta\sin^n\theta\dd\theta,
	\qquad
	\zeta_j=\frac{\sqrt\pi\,\Gamma((j+1)/2)}{\Gamma(j/2+1)}.
\end{equation}
For a monomial \(x^ay^b\), put \(q=a+b\),
\begin{equation}
	\chi_q=\begin{cases}1,&q=0,\\1+1/q,&q>0,\end{cases}
	\qquad g_q=q+1+3\chi_q,
\end{equation}
and abbreviate
\begin{align}
	M_{ab}&=\mu_{ab}^\Theta\zeta_{q+1},&
	X_{ab}&=\mu_{a+2,b}^\Theta\zeta_{q+3},&
	Y_{ab}&=\mu_{a+1,b+1}^\Theta\zeta_{q+3},\\
	Z_{ab}&=\mu_{a,b+2}^\Theta\zeta_{q+3},&
	P^x_{ab}&=(q+2)\mu_{a+1,b}^\Theta\zeta_{q+2},&
	P^y_{ab}&=(q+2)\mu_{a,b+1}^\Theta\zeta_{q+2}.
\end{align}
With output and input both ordered as \((\pi,u,v)\), all nine component
formulas collapse to the symmetric block
\begin{equation}
		\mathsf C_{ab}^{\Theta}
		=\frac1{4\pi}
		\begin{bmatrix}
			(q+1)M_{ab}&-P^x_{ab}&-P^y_{ab}\\
			-P^x_{ab}&g_qX_{ab}-\chi_qM_{ab}&g_qY_{ab}\\
			-P^y_{ab}&g_qY_{ab}&g_qZ_{ab}-\chi_qM_{ab}
		\end{bmatrix}
		+\frac{2\omega_\Theta}{3}\mathbf1_{\{q=0\}}
		\operatorname{diag}(0,1,1).
	\label{eq:compact-monomial-block}
\end{equation}
The contribution at radius \(r=c_0\tau/h\) is
\(r^{a+b}\mathsf C_{ab}^{\Theta}\).   Since \(a,b\le2\), every nodal acoustic entry is a degree-four
polynomial in \(r\), not a fourth-order Taylor approximation.

\subsection{Semidiscrete method}
\label{sec:sdaf-symbol}

Let \(A_d=R_d\Lambda_dR_d^{-1}\) be the primitive Euler Jacobian in direction
\(d\in\{x,y\}\), where
\begin{equation}
	\begin{split}
		R_x&=\begin{bmatrix}
			c_0^{-2}&1&0&c_0^{-2}\\-Z_0^{-1}&0&0&Z_0^{-1}\\
			0&0&1&0\\1&0&0&1
		\end{bmatrix},
		&\Lambda_x&=\operatorname{diag}(u_0-c_0,u_0,u_0,u_0+c_0),\\[-1mm]
		R_y&=\begin{bmatrix}
			c_0^{-2}&1&0&c_0^{-2}\\0&0&1&0\\
			-Z_0^{-1}&0&0&Z_0^{-1}\\1&0&0&1
		\end{bmatrix},
		&\Lambda_y&=\operatorname{diag}(v_0-c_0,v_0,v_0,v_0+c_0).
	\end{split}
\end{equation}
The characteristic splitting is
\begin{equation}
	A_d^\pm=R_d\Lambda_d^\pm R_d^{-1},
	\qquad \Lambda_d^\pm=\tfrac12(\Lambda_d\pm|\Lambda_d|).
	\label{eq:sdaf-split}
\end{equation}
}
\section{Additional details and plots for Shear Problems}
\label{sec:Addplot}
This section provides additional context and plots pertaining to the shear problems discussed previously.

The Discontinuous Galerkin method used for comparisons employs a modal basis with p+1 Gauss-Legendre quadrature points per dimension. HLLC flux was used between elements and no spatial limiting was used. Time integration was performed using TVD RK3 and a CFL number was set to  0.4/(2p+1)=0.13. The DG method and the present method have the same number of degrees of freedom for a given mesh.

The Continuous Galerkin method is a spectral element scheme based on summation-by-parts and two point fluxes that ensure additional consistency of Kinetic energy, entropy and internal energy. Additional stabilization is added through gradient jump penalty terms. The results use degree 4 polynomials and classical RK4 method; the grids with cells $32 \times 16, 64 \times 32$ and $128 \times 64$ correspond to equivalent AF grids of sizes $128 \times 64$, $256 \times 128$ and {$512 \times 256$} respectively.

\begin{figure}
	\centering
	\includegraphics[width=0.48\textwidth]{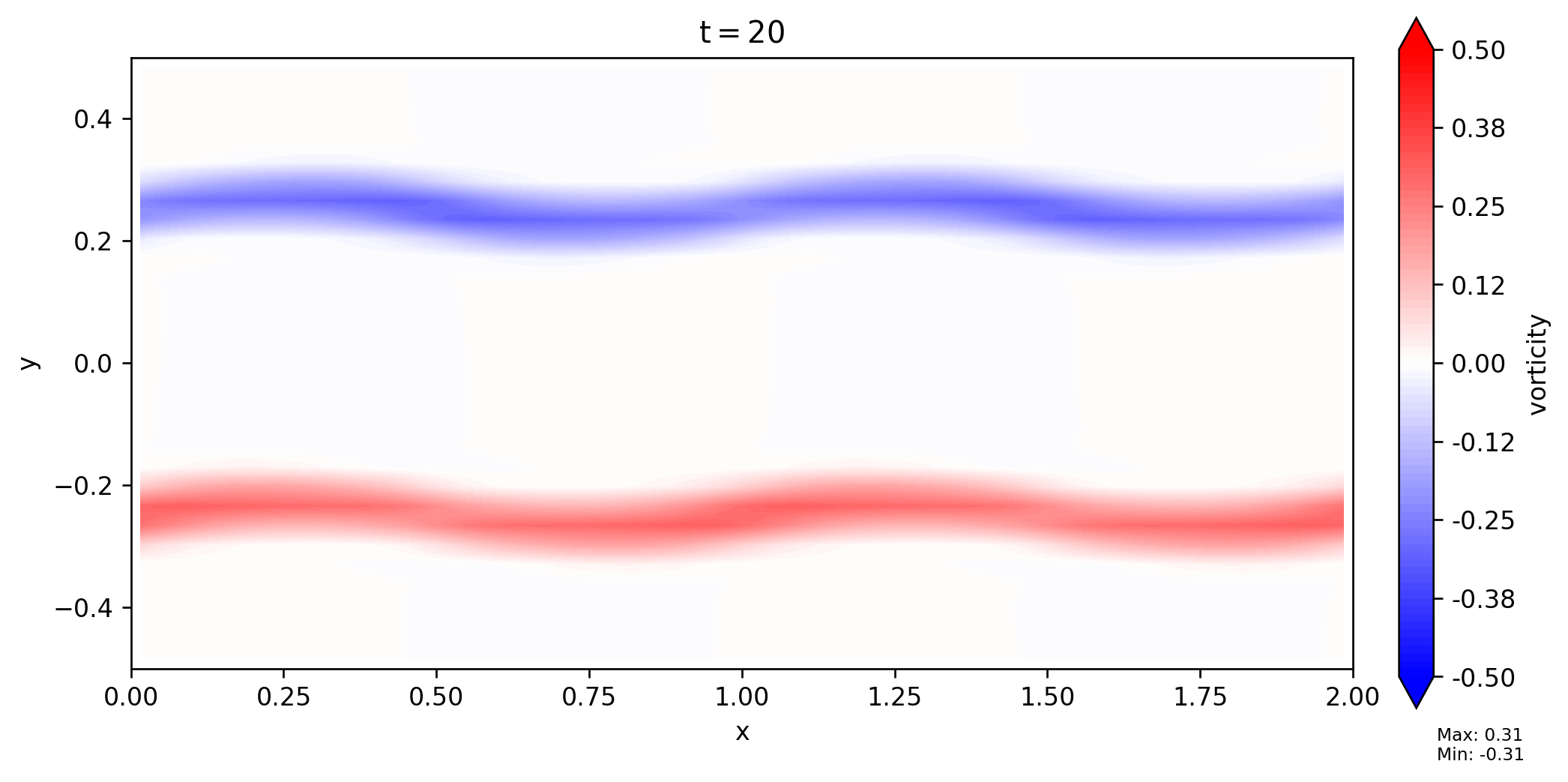}
	\includegraphics[width=0.48\textwidth]{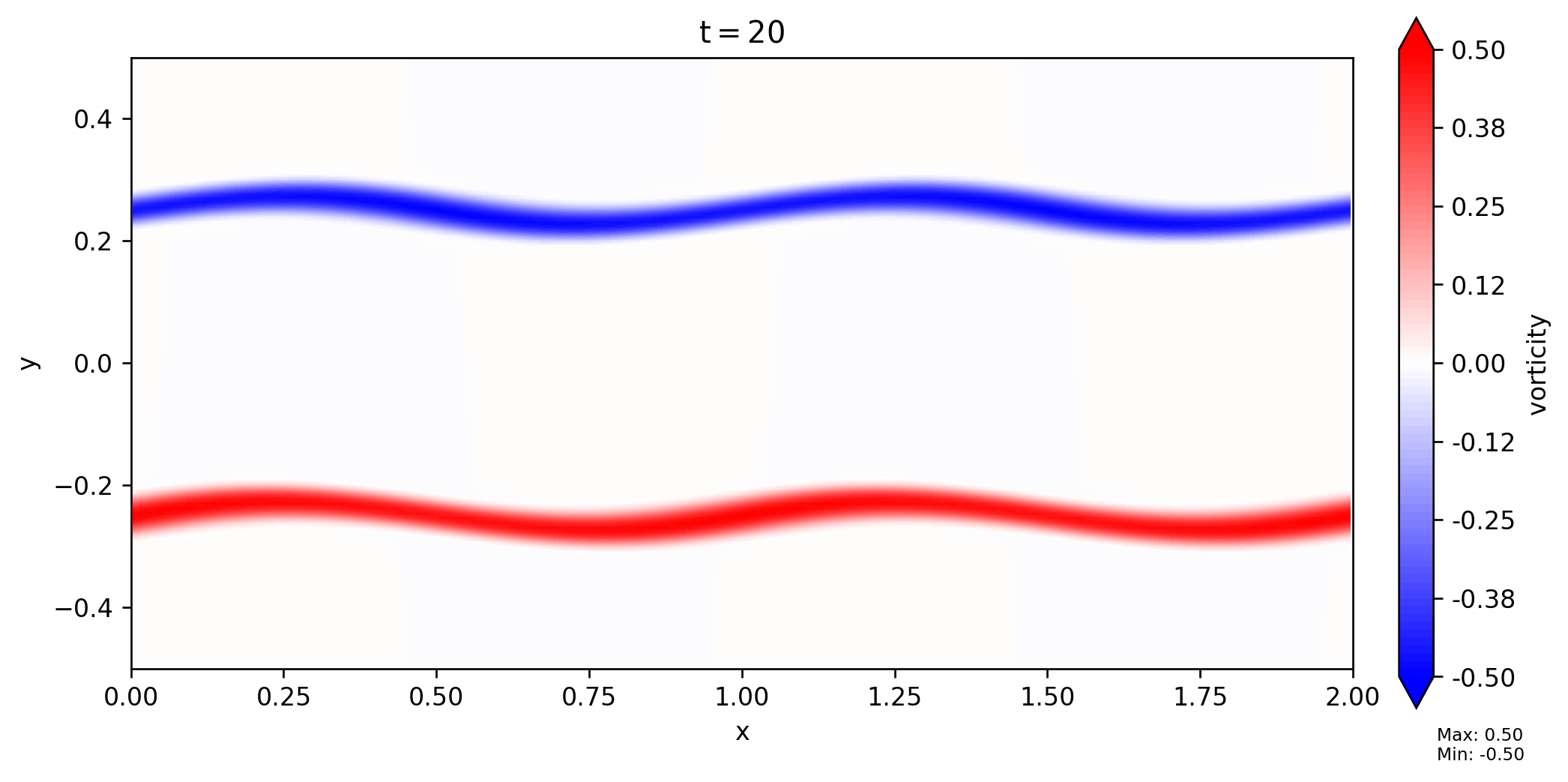}
	\includegraphics[width=0.48\textwidth]{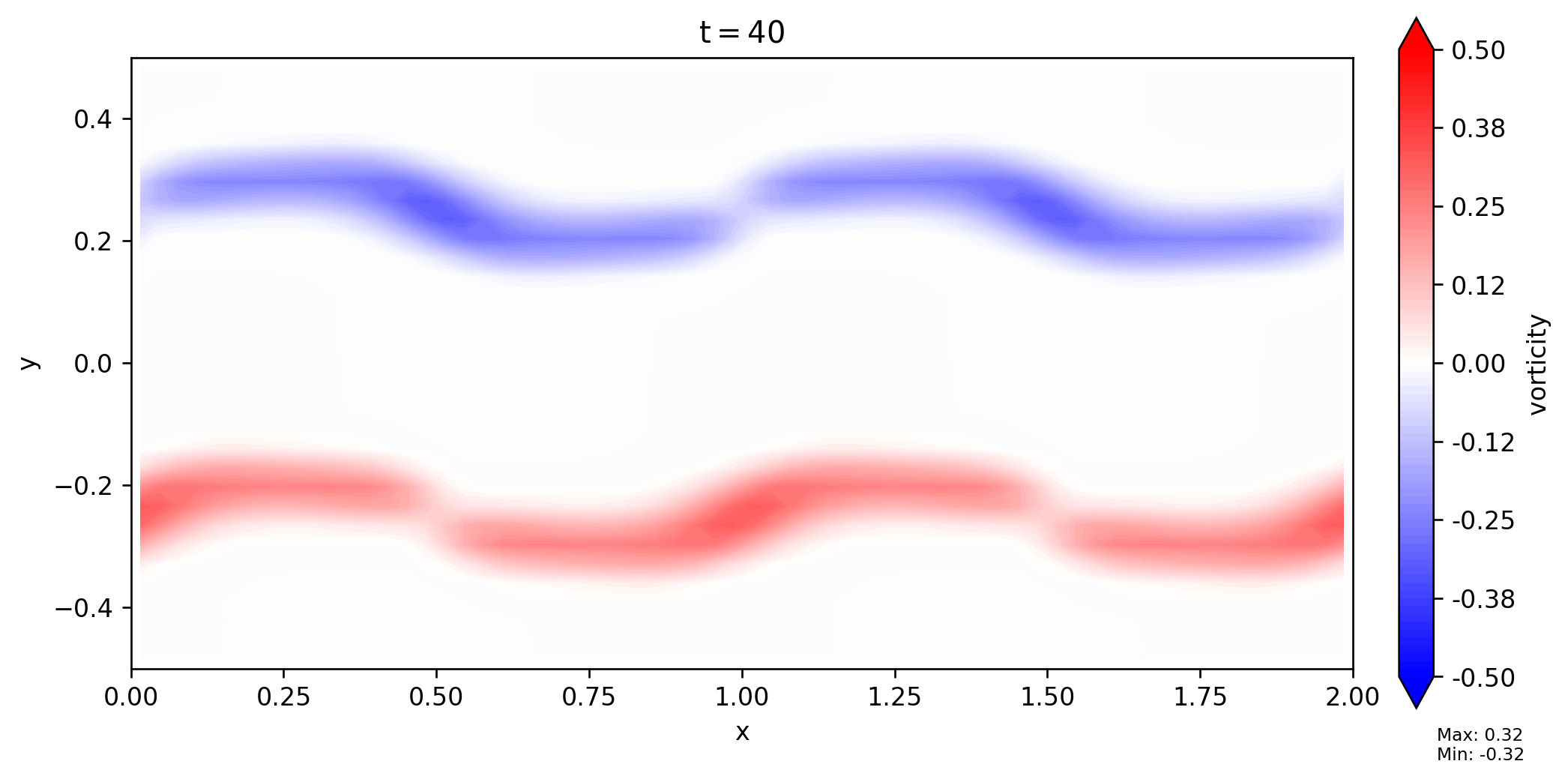}
	\includegraphics[width=0.48\textwidth]{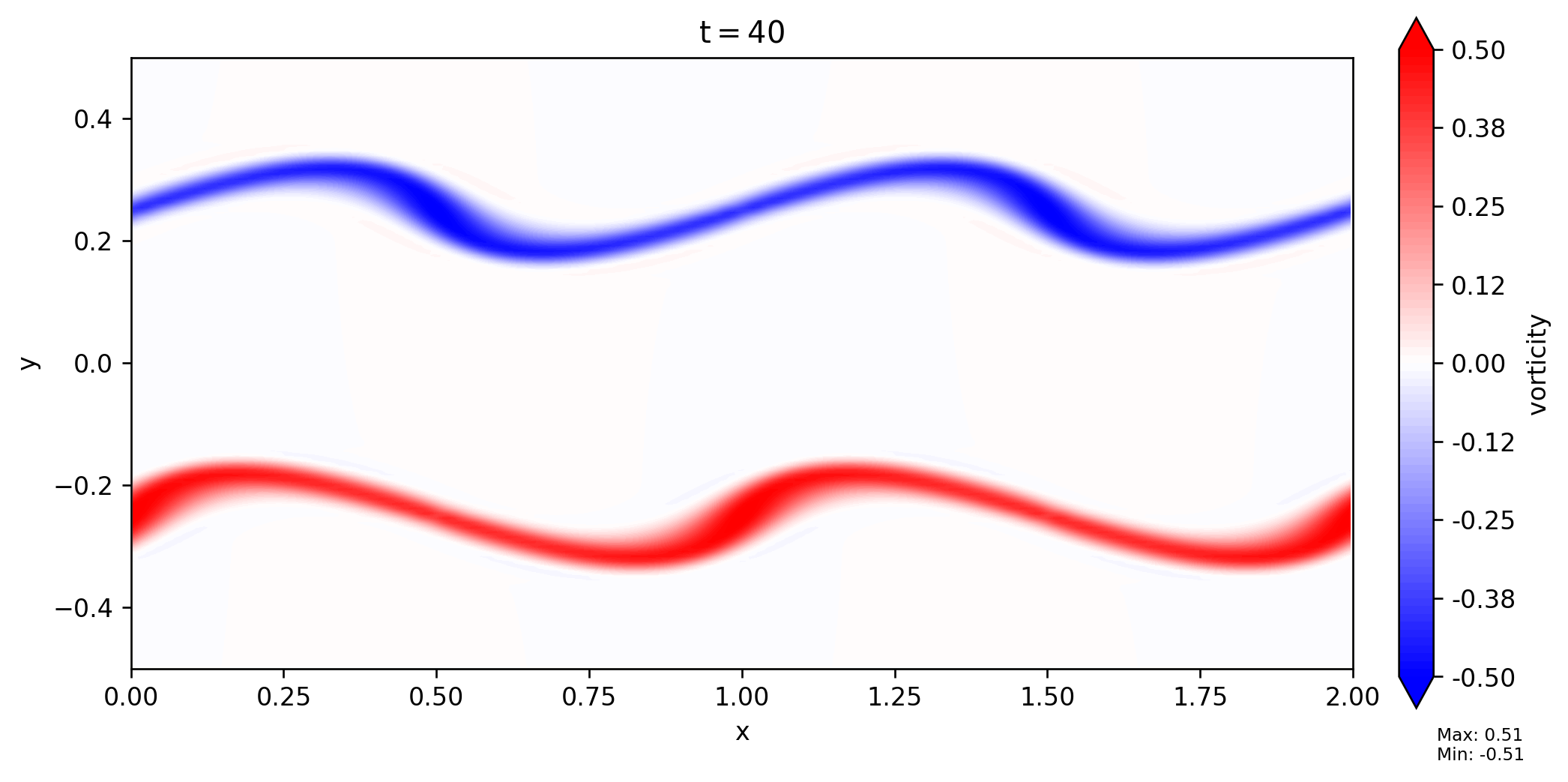}
	\includegraphics[width=0.48\textwidth]{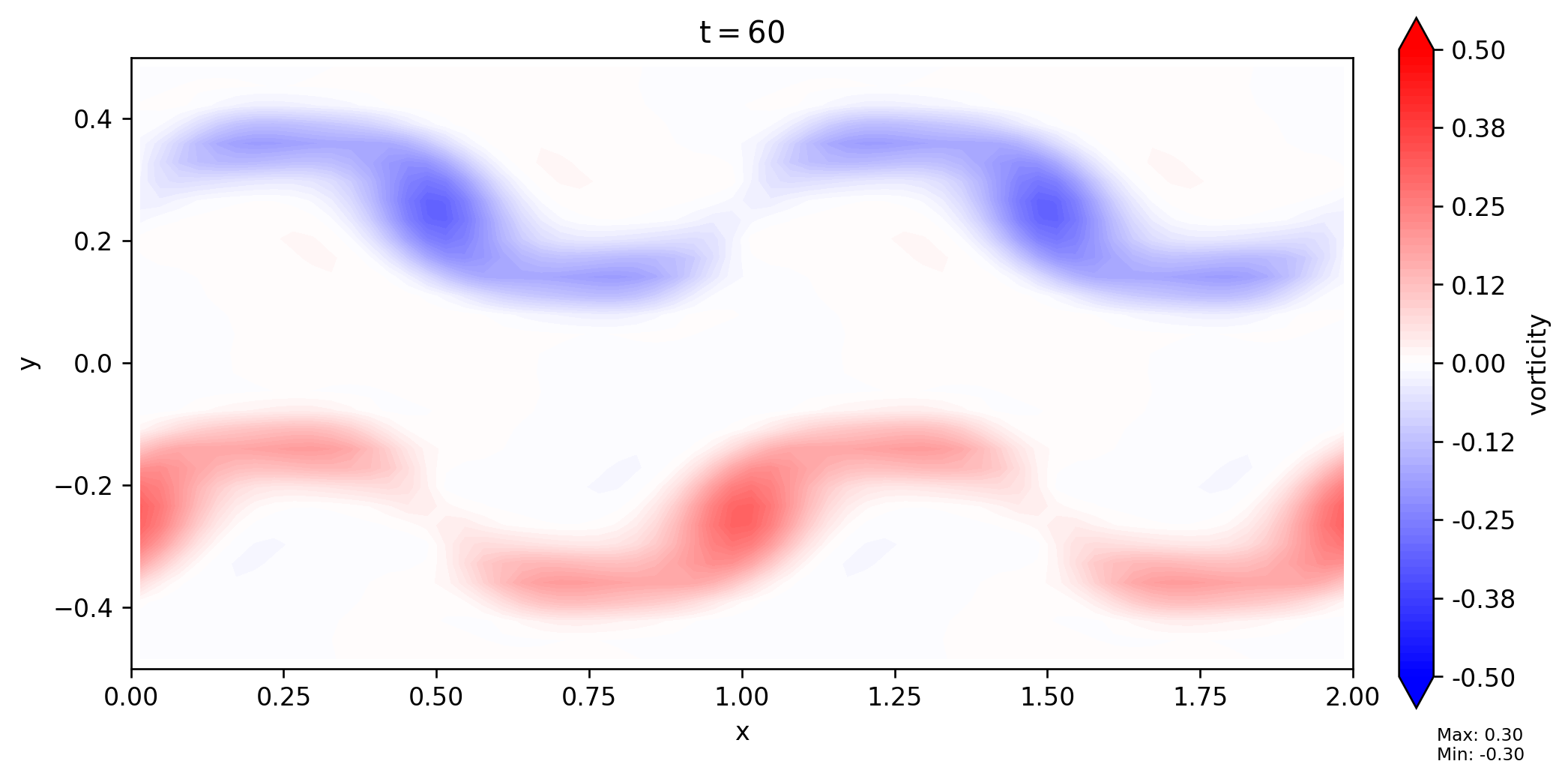}
	\includegraphics[width=0.48\textwidth]{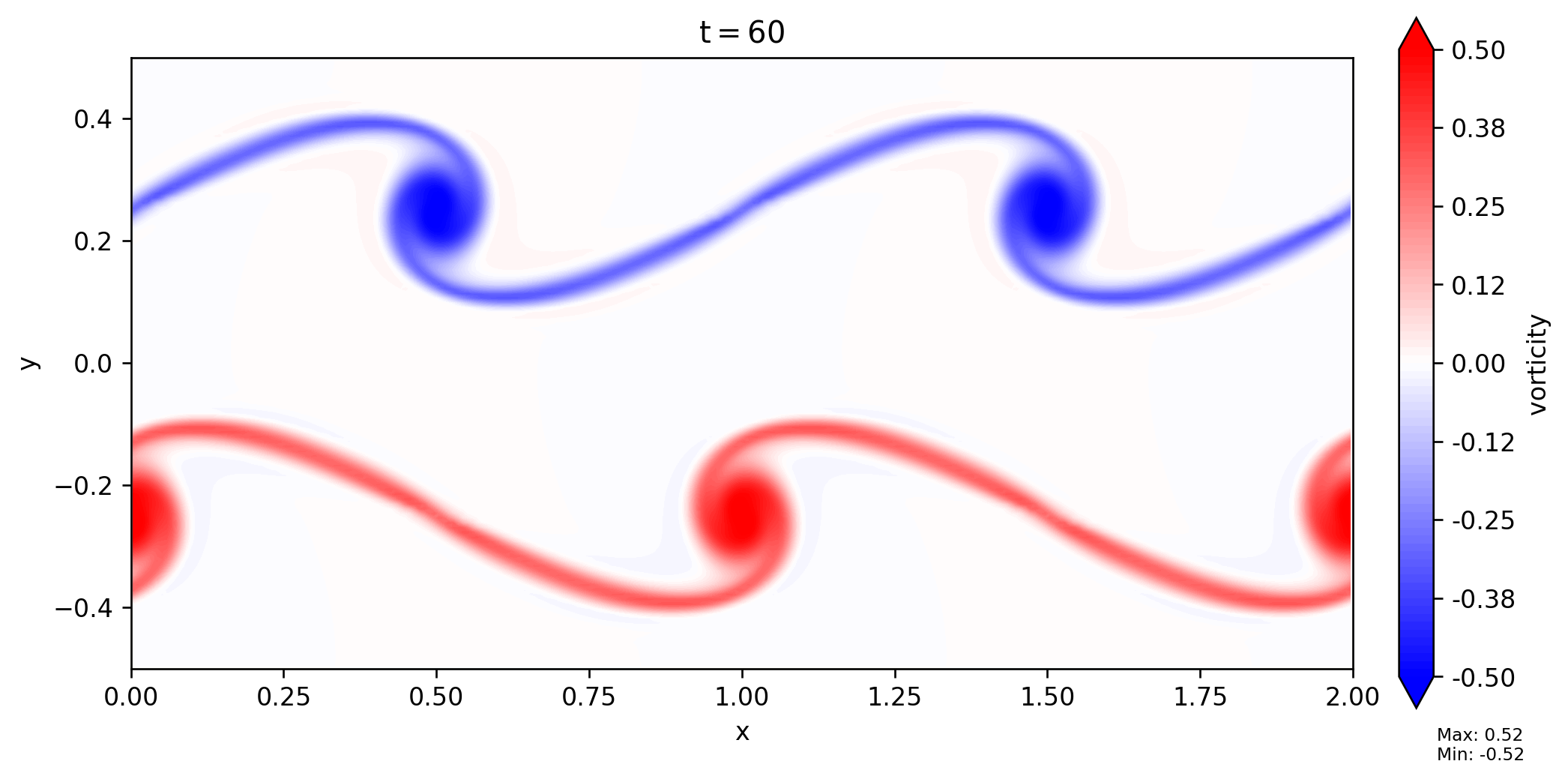}
	\includegraphics[width=0.48\textwidth]{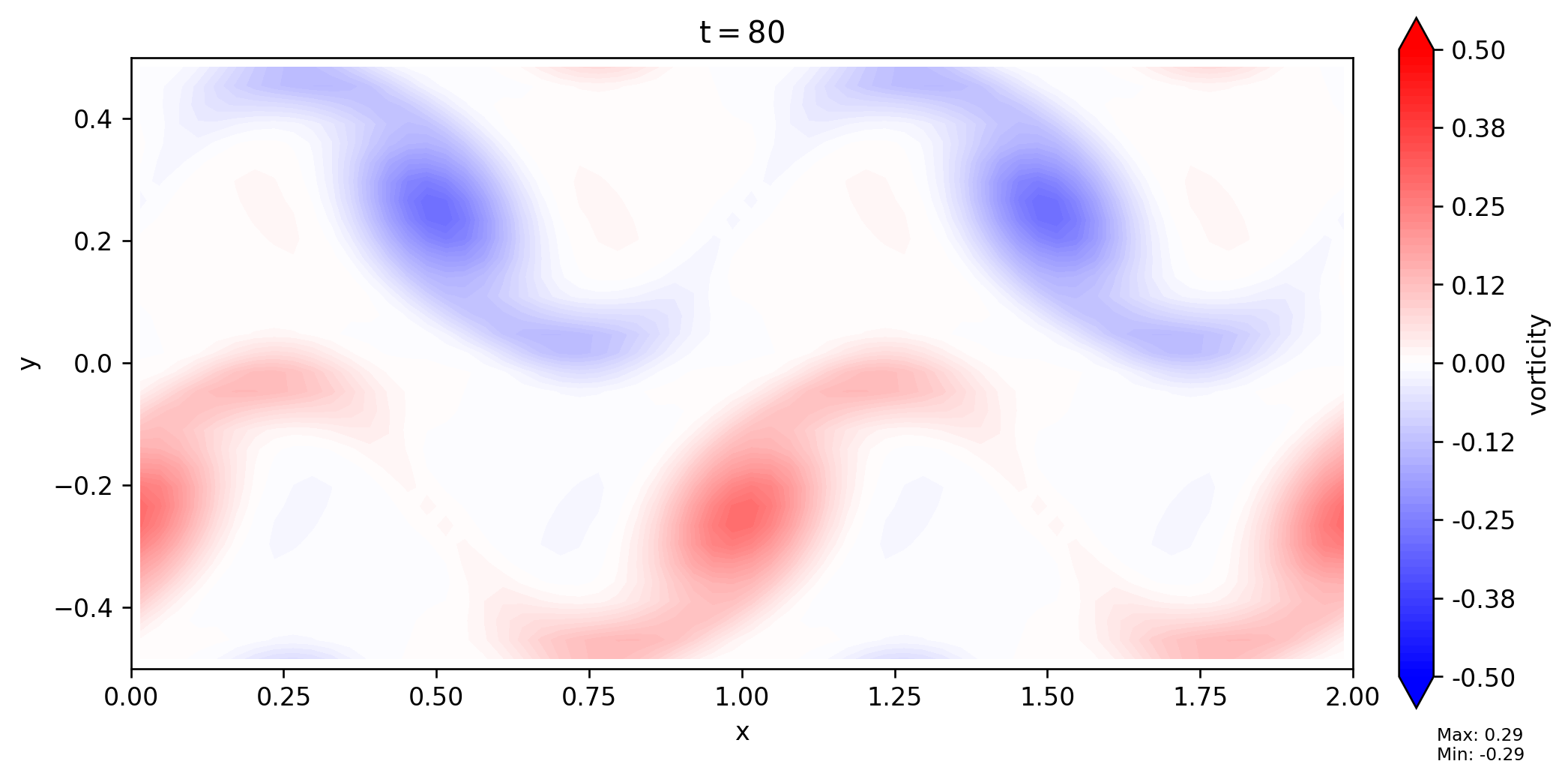}
	\includegraphics[width=0.48\textwidth]{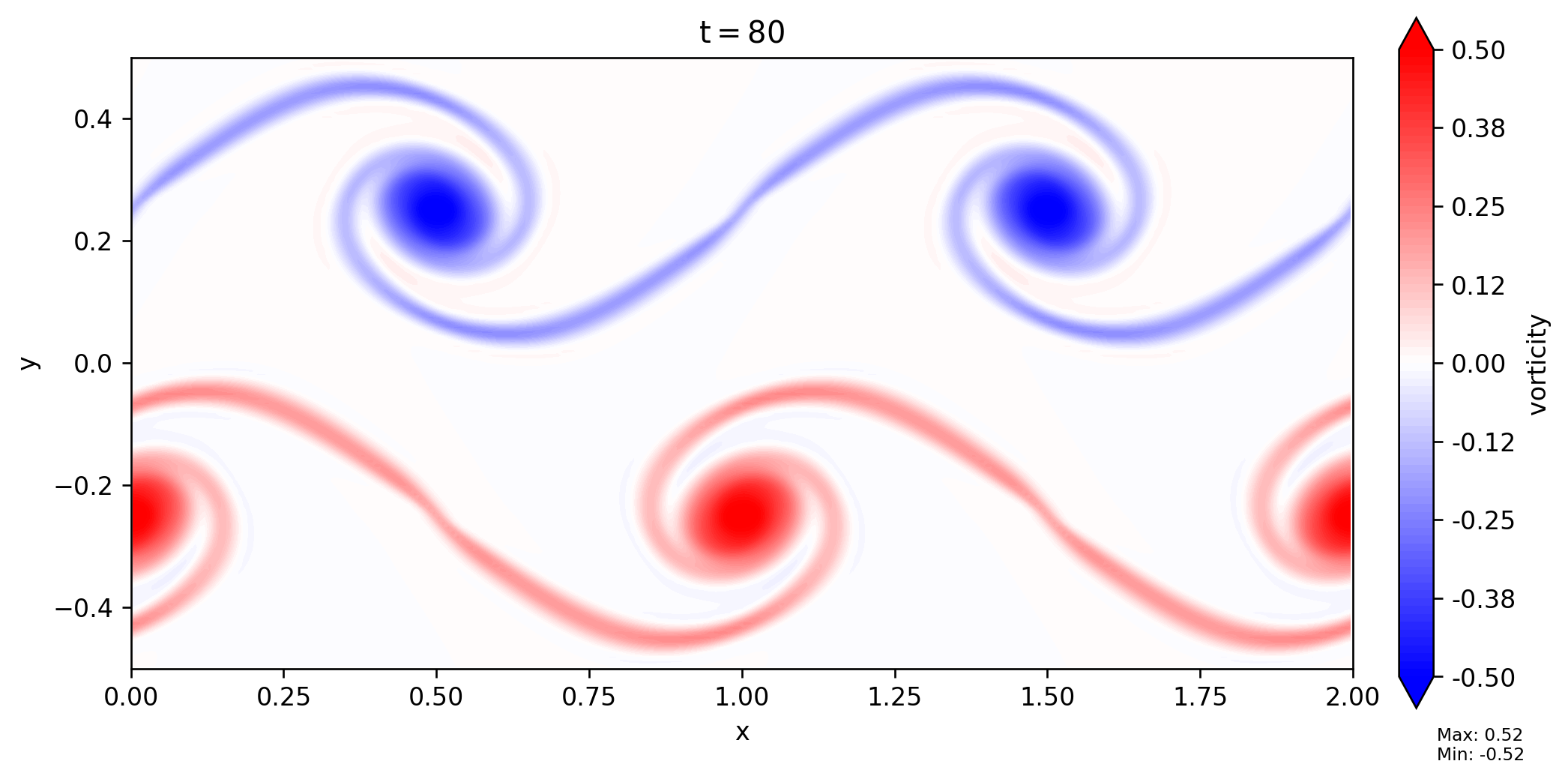}
	\caption{Shear Problem on $64\times 32$ and $256 \times 128$ mesh, Discrete (RB-TAI)}
	\label{fig:shearvort}
\end{figure}

\begin{figure}
	\centering
	\includegraphics[width=0.49\textwidth]{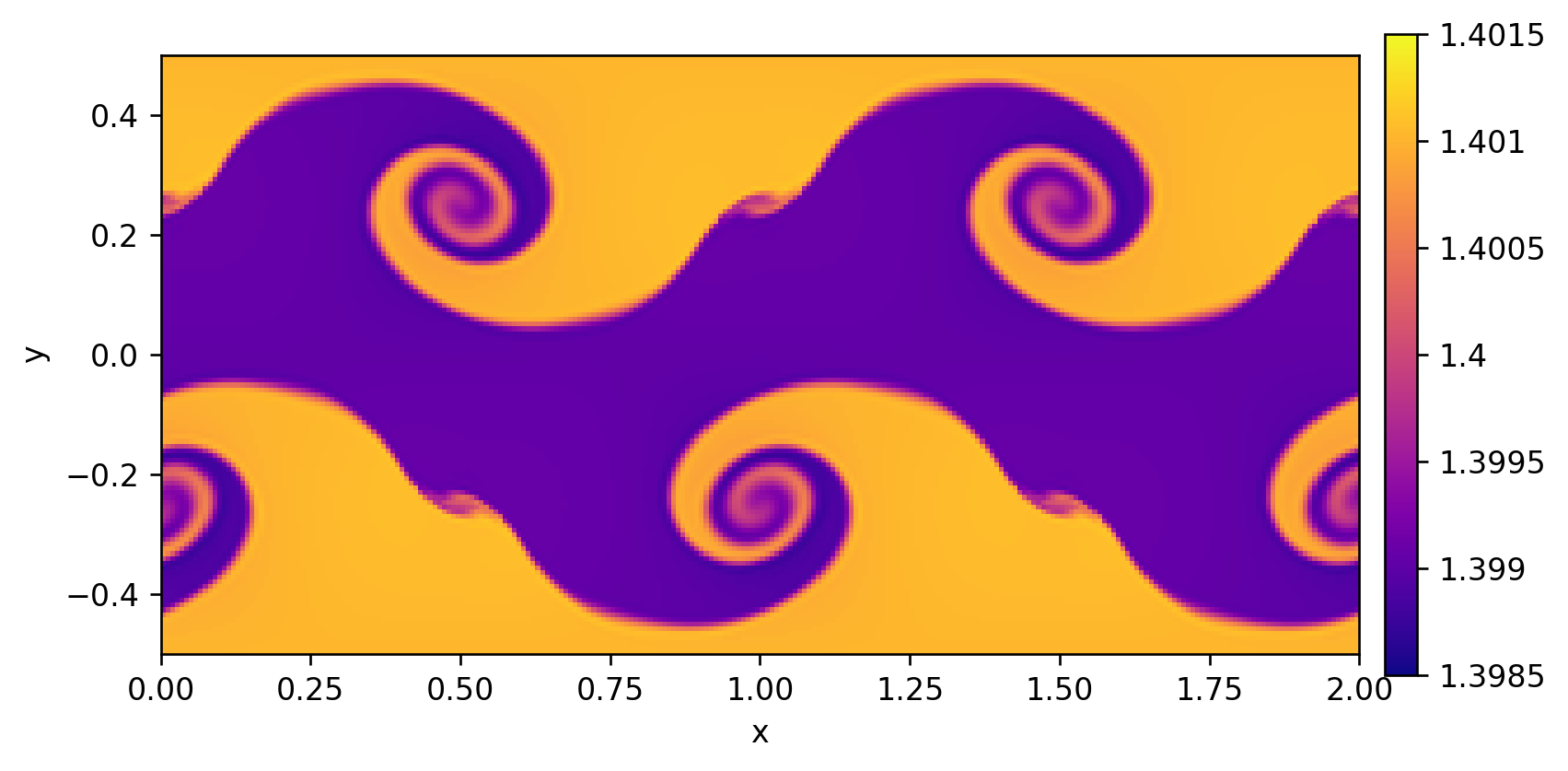}
	\includegraphics[width=0.49\textwidth]{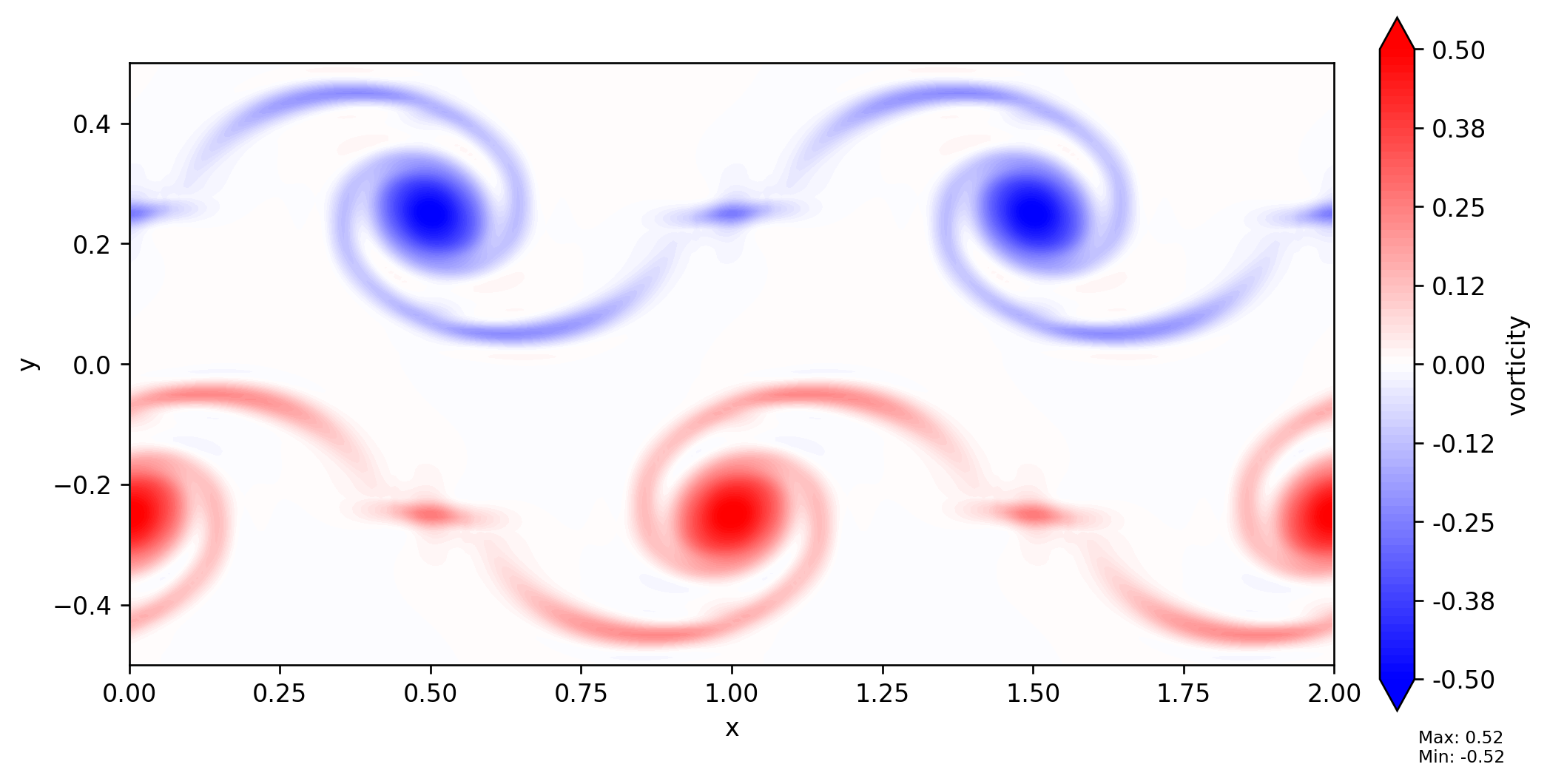}
	\caption{Density and Vorticity contours in shear problem  with P1 Discontinuous Galerkin on $256 \times 128$ grid.}
	\label{fig:shearcomp2}
\end{figure}

\begin{figure}
	\centering
	\includegraphics[trim={0 0 0 1.3cm},clip,width=0.5\textwidth]{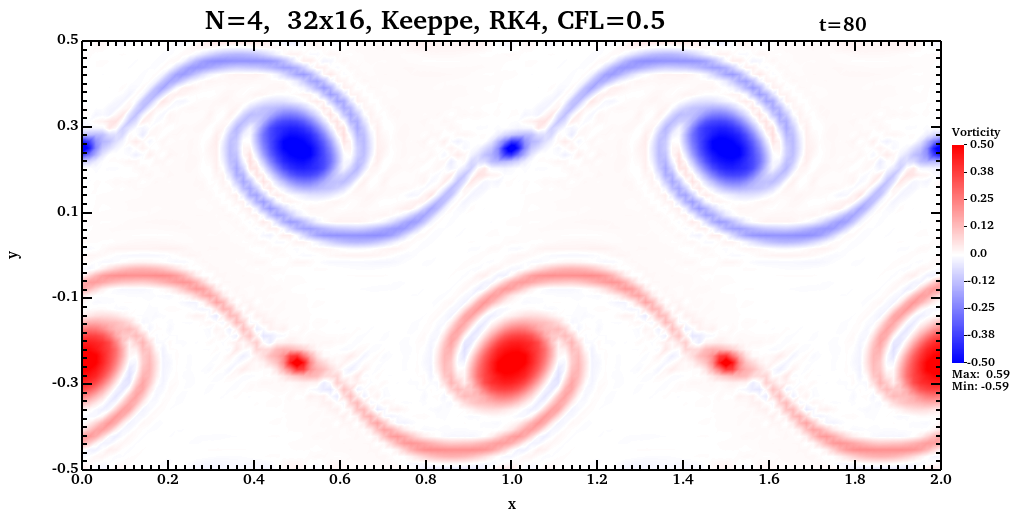}
	\includegraphics[trim={0 0 0 1.3cm},clip,width=0.5\textwidth]{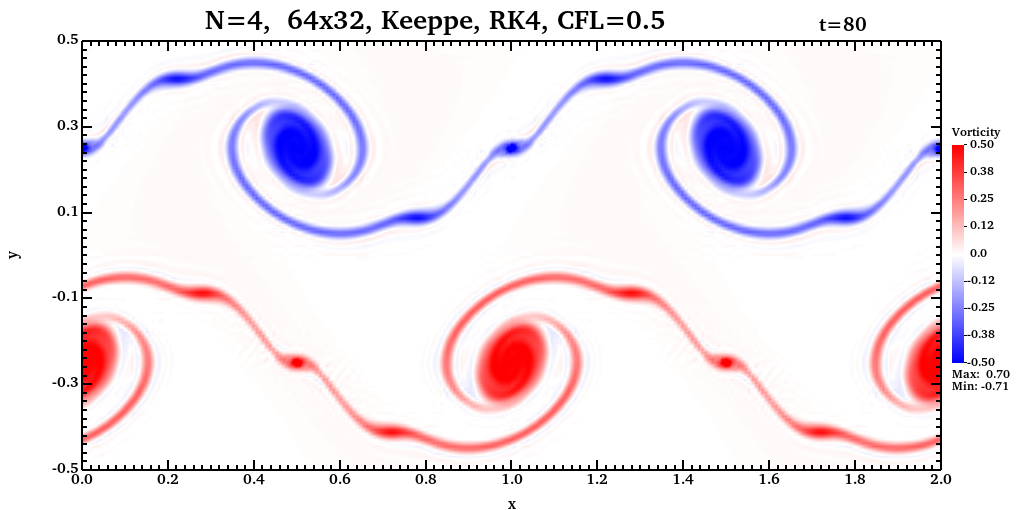}
	\includegraphics[trim={0 0 0 1.3cm},clip,width=0.5\textwidth]{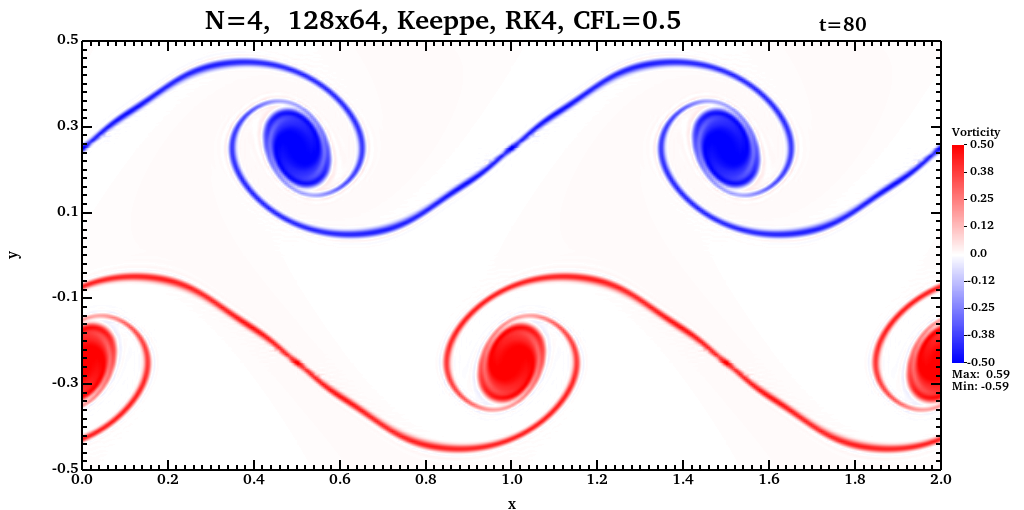}
	\caption{Vorticity contours in shear problem  with 5th order Continuous Galerkin: $32\times 16$; $64\times 32$; $128\times 64$ grids.}
	\label{fig:shearcomp3}
\end{figure}

\begin{figure}
\centering
\includegraphics[width=0.49\textwidth]{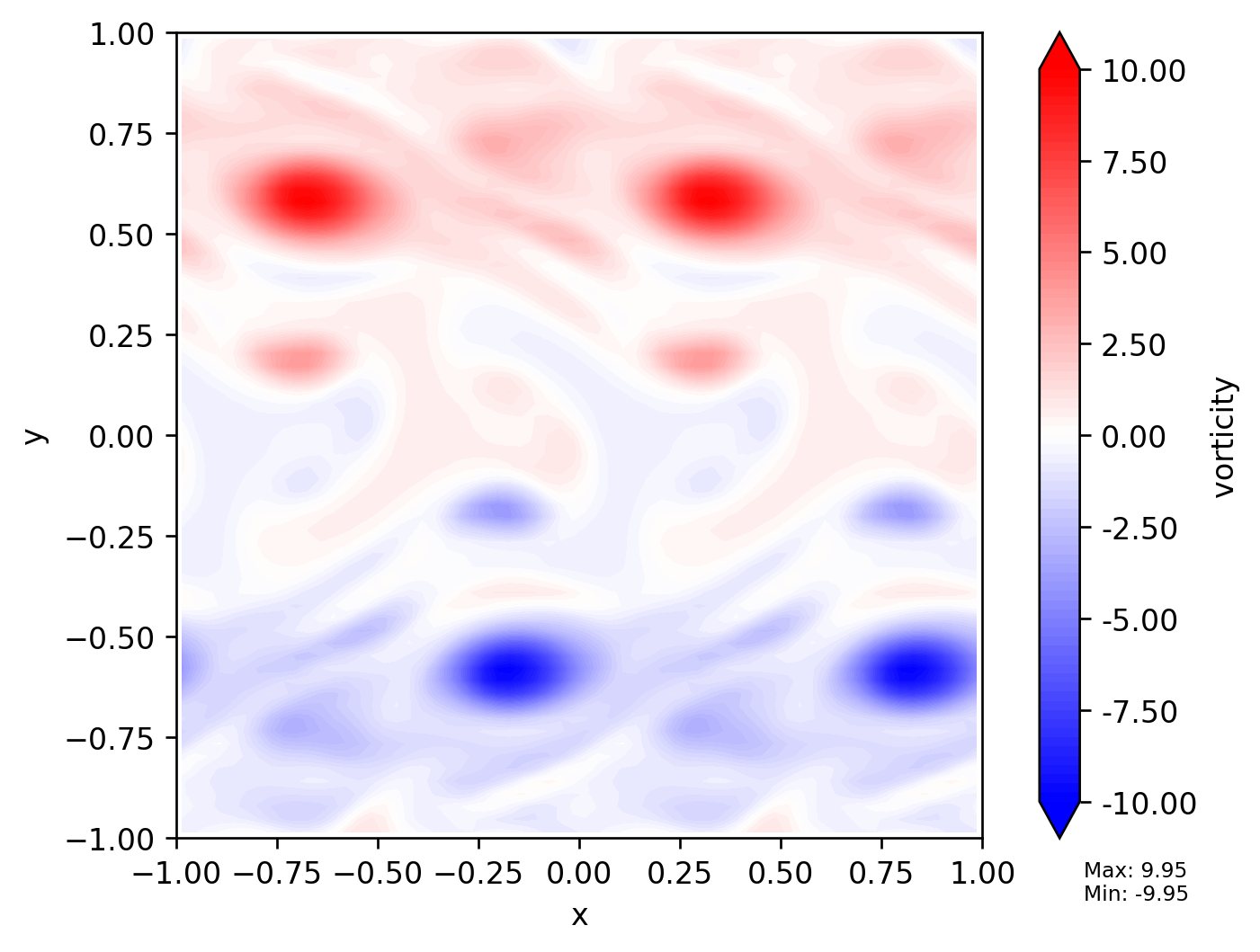}
\includegraphics[width=0.49\textwidth]{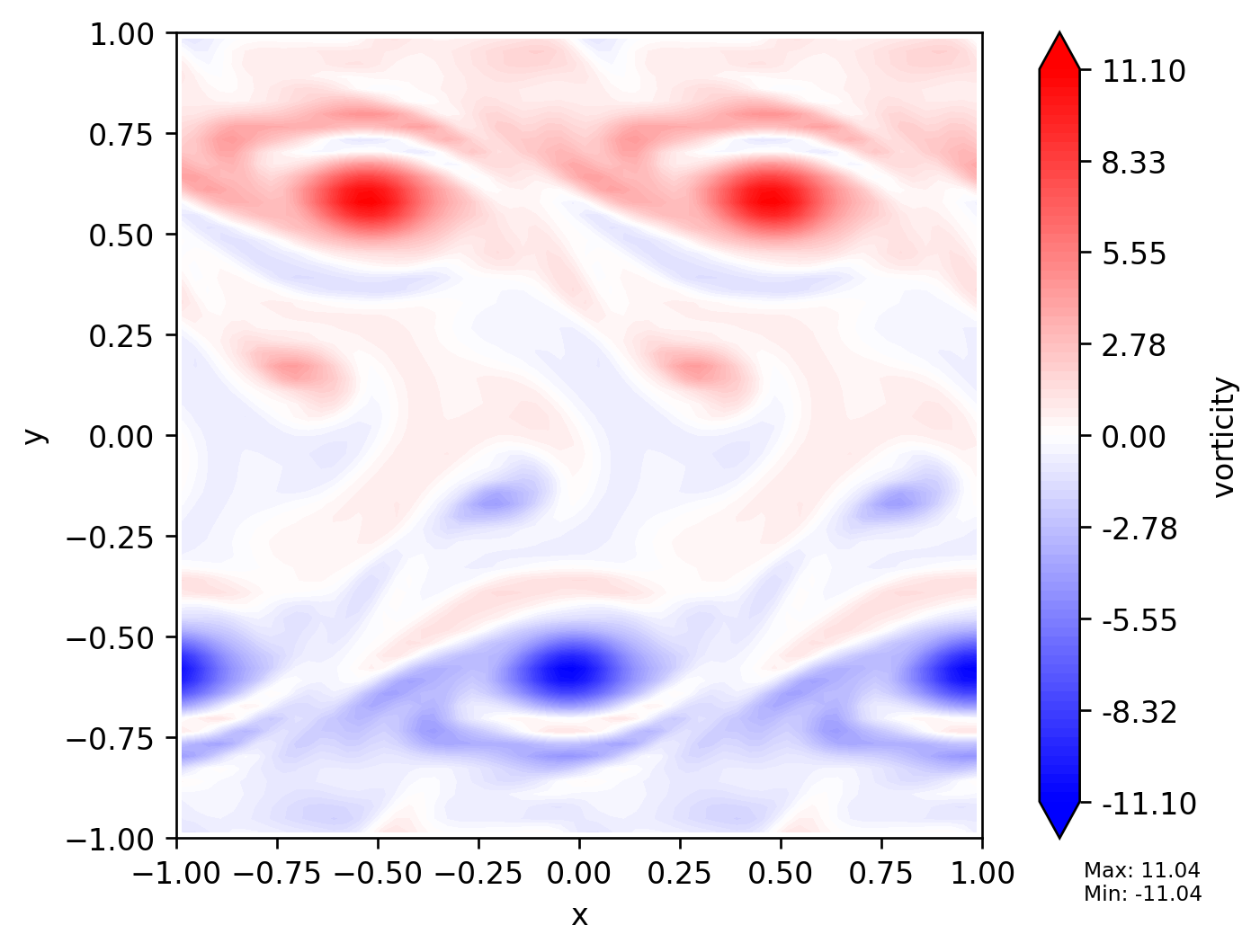}
\includegraphics[width=0.49\textwidth]{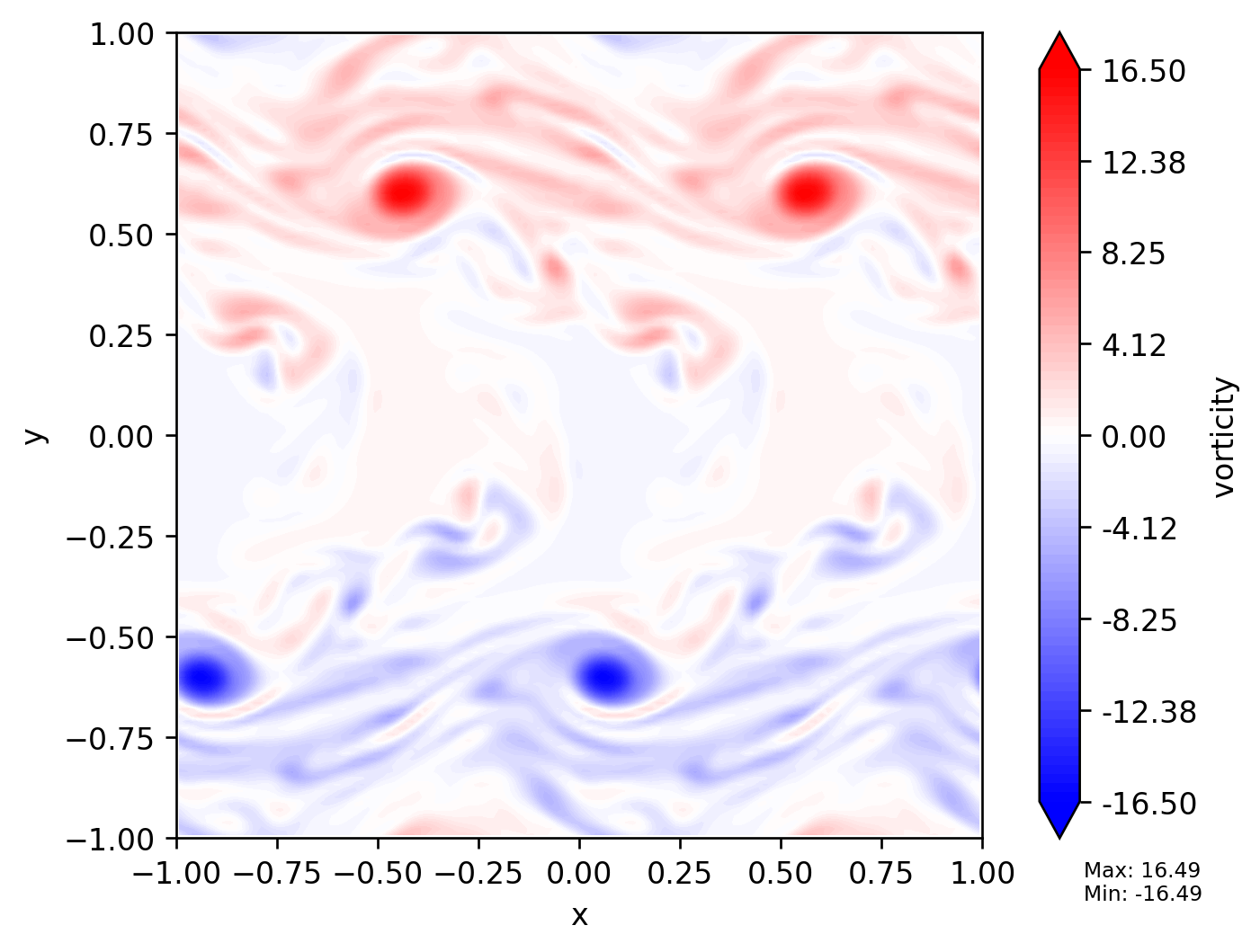}
\includegraphics[width=0.49\textwidth]{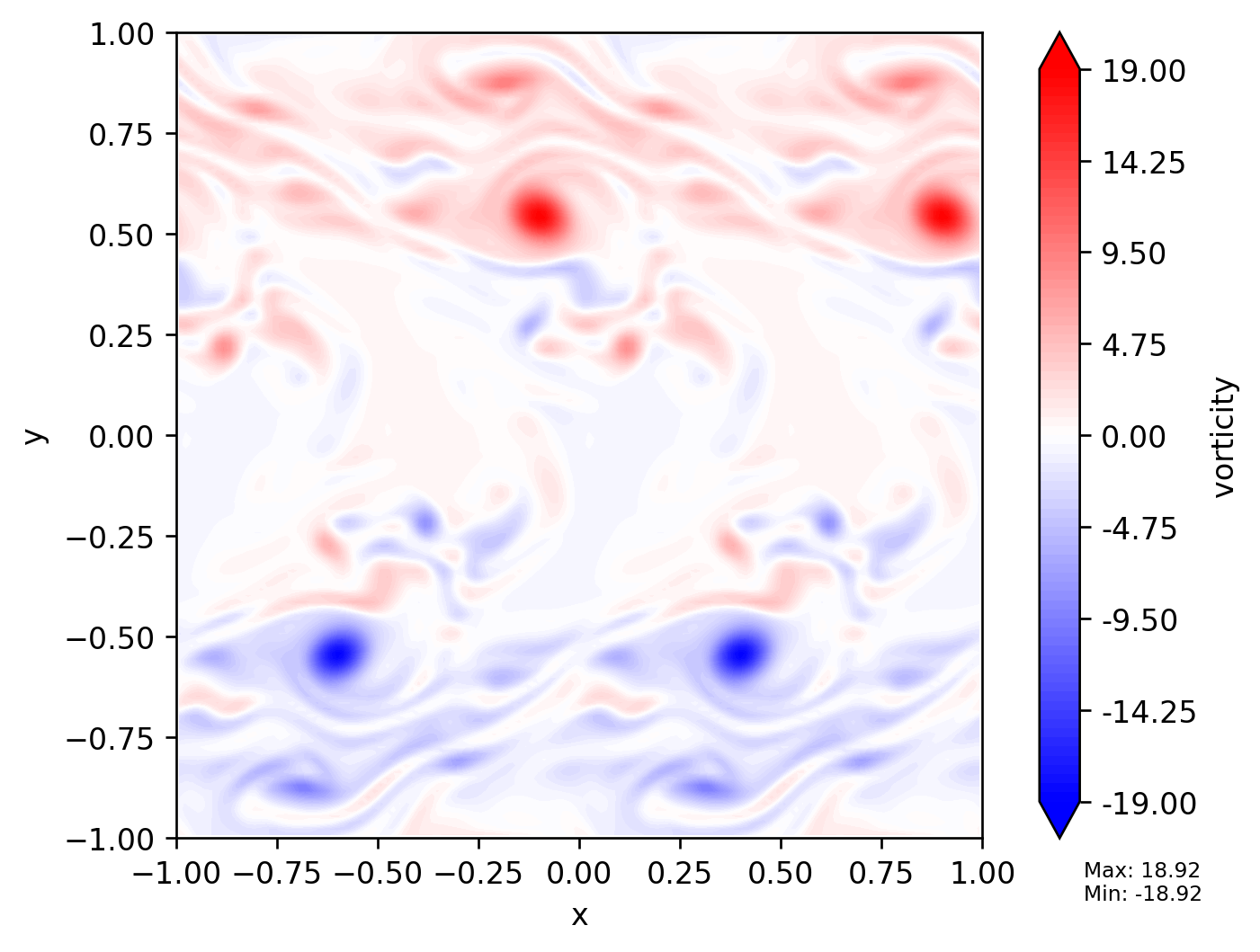}
\includegraphics[width=0.49\textwidth]{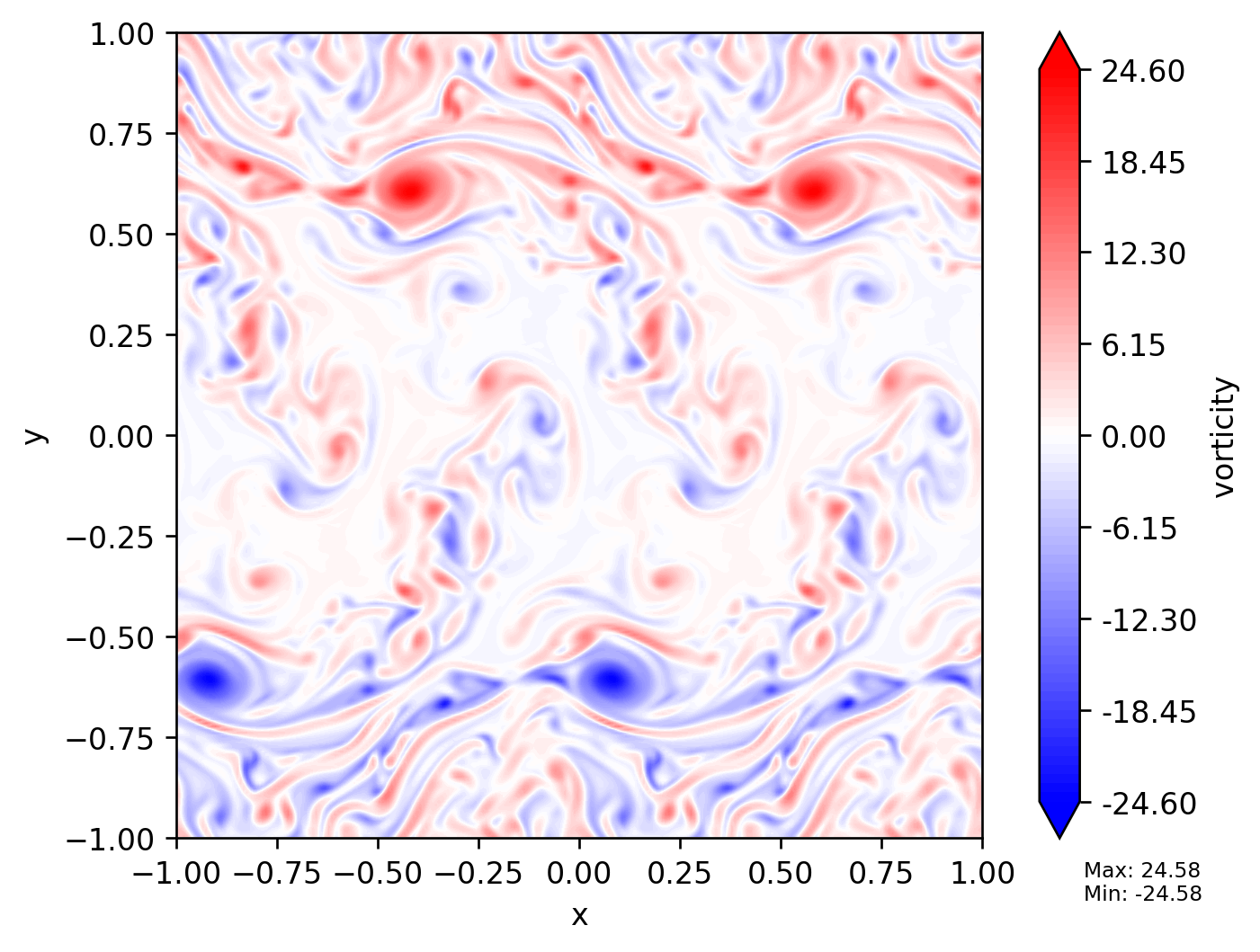}
\includegraphics[width=0.49\textwidth]{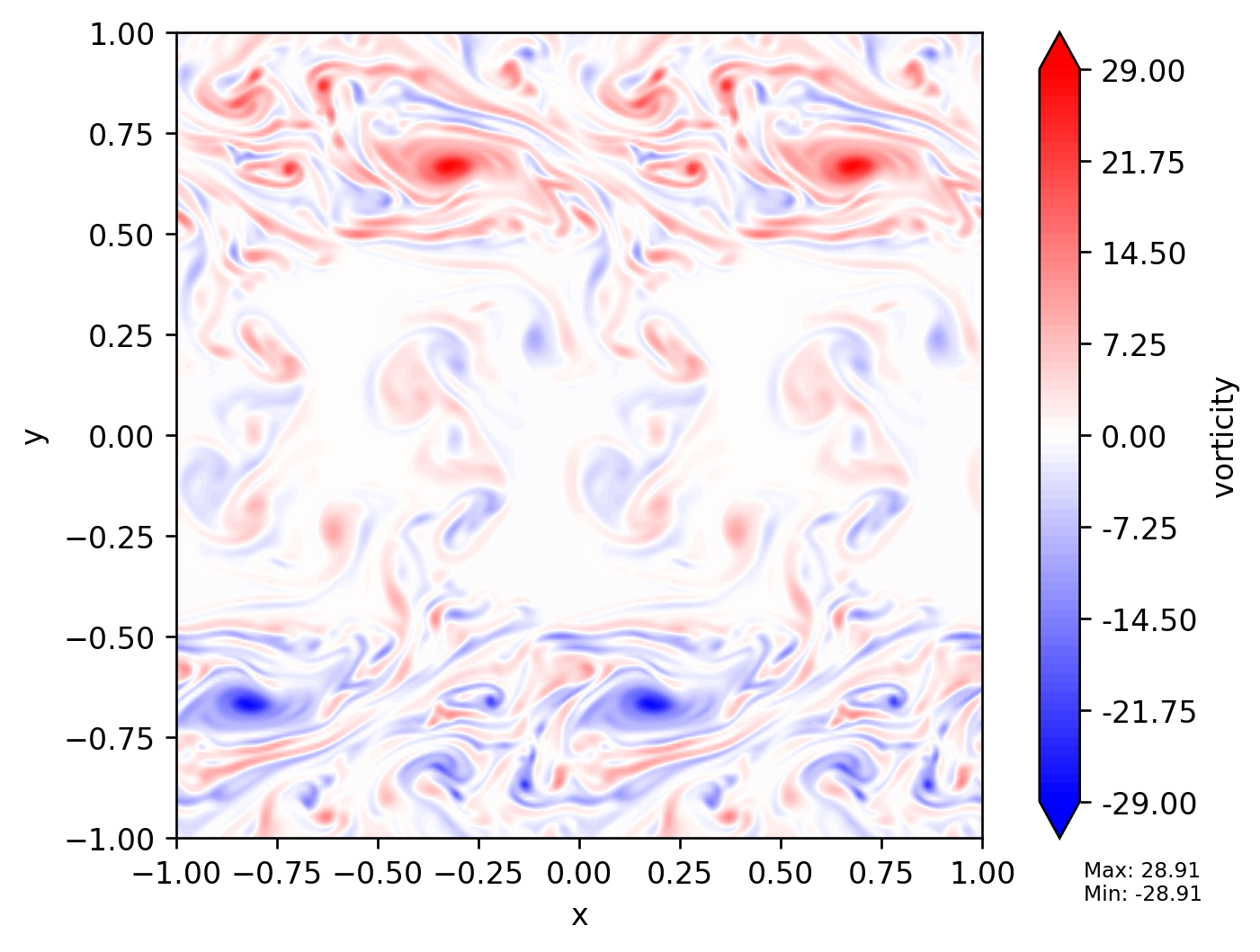}
\caption{Vorticity contours for under-resolved KH problem. Rows:  $64^2,128^2,256^2$ mesh. Columns: Discrete RB, Discrete RB-TAI}
\label{fig:shearnewcomp}
\end{figure}

\bibliographystyle{cas-model2-names}

\bibliography{references}

\end{document}